\documentclass[aap]{imsart}
\usepackage{mathrsfs}
\RequirePackage{amsthm,amsmath,amsfonts,amssymb}
\RequirePackage[numbers,sort&compress]{natbib}
\RequirePackage[colorlinks,citecolor=blue,urlcolor=blue]{hyperref}
\RequirePackage{graphicx}

\startlocaldefs
\numberwithin{equation}{section}
\theoremstyle{plain}

\newtheorem{theorem}{Theorem}[section]
\newtheorem{lemma}[theorem]{Lemma}
\newtheorem{condition}[theorem]{Condition}
\newtheorem{remark}[theorem]{Remark}
\newtheorem{example}[theorem]{Example}
\newtheorem{proposition}[theorem]{Proposition}

\theoremstyle{remark}

\newcommand{\p}{\partial}
\newcommand{\e}{\varepsilon}

\newcommand{\R}{{\mathbb R}}
\newcommand{\Z}{{\mathbb Z}}

\newcommand{\pP}{{\mathbb P}}
\newcommand{\I}{{\mathbb I}}
\newcommand{\E}{{\mathbb E}}

\newcommand{\ty}{\infty}

\newcommand{\mT}{{\mathbb T}}

\newcommand{\aA}{{\cal A}}

\newcommand{\FF}{{\cal F}}

\def\mE{{\mathbb E}}
\def\dif{{\mathord{{\rm d}}}}
\def\cF{{\mathcal F}}
\def\cZ{{\mathcal Z}}
\def\cK{{\mathcal K}}
\def\mP{{\mathbb P}}
\def\mN{{\mathbb N}}
\def\mR{{\mathbb R}}
\def\mZ{{\mathbb Z}}
\def\cB{{\mathcal B}}
\def\cD{{\mathcal D}}
\def\cC{{\mathcal C}}
\def\cM{{\mathcal M}}
\def\cA{{\mathcal A}}
\def\cS{{\mathcal S}}

\def\cR{{\mathcal R}}
\def\cH{{\mathcal H}}
\def\eps{\varepsilon}
\def\cL{{\mathcal L}}

\def\mI{{\mathbb I}}
\def\cO{{\mathcal O}}
\def\cI{{\mathcal I}}
\def\mQ{{\mathbb Q}}

\newcommand{\lag}{\langle}
\newcommand{\rag}{\rangle}

\newcommand{\dd}{{\textup d}}

\def\[{{\big[}}
\def\]{{\big]}}

\endlocaldefs

\begin{document}

\begin{frontmatter}
\title{Ergodic and  mixing properties   of the 2D Navier-Stokes equations with a degenerate  multiplicative  Gaussian noise}
\runtitle{Ergodic and  mixing properties  for   stochastic 2D Navier-Stokes equations }

\begin{aug}
\author[A]{\fnms{Zhao}~\snm{Dong}
\ead[label=e1]{dzhao@amt.ac.cn}},
\author[B]{\fnms{Xuhui}~\snm{Peng}\ead[label=e2]{xhpeng@hunnu.edu.cn}
}
\address[A]{Academy of Mathematics and Systems Science,   University of Chinese Academy of Sciences,  Chinese  Academy of Sciences, Beijing
100190, P.R.China
 \printead[presep={ ,\ }]{e1}}

\address[B]{MOE-LCSM, School of Mathematics and Statistics, Hunan Normal
University, Changsha, Hunan410081, P.R.China
\printead[presep={,\ }]{e2}}
\end{aug}

\begin{abstract}
In this paper, we establish the ergodic and mixing properties of stochastic 2D Navier-Stokes equations driven by a highly degenerate multiplicative Gaussian noise. The noise can appear in as few as four directions, and its intensity depends on the solution. The case of additive Gaussian noise was previously treated by Hairer and Mattingly [\emph{Ann. of Math.}, 164(3):993--1032, 2006]. To derive the ergodic and mixing properties in the present setting, we employ Malliavin calculus to establish the asymptotically strong Feller property. The primary challenge lies in proving the "invertibility" of the Malliavin matrix, which differs fundamentally from the additive case.
\end{abstract}

\begin{keyword}[class=MSC]
\kwd{60H15;	60H07}
\end{keyword}

\begin{keyword}
\kwd{Stochastic Navier--Stokes equations}
\kwd{Malliavin calculus}
\kwd{Ergodicity}
\end{keyword}

\end{frontmatter}

\section{Introduction and main results}

In the physics literature, when discussing the behavior of stochastic fluid dynamics in the turbulent regime, the primary assumptions typically include ergodicity and statistical translational invariance of the stationary state. The theory of turbulence often concerns cases where the noise is degenerate, meaning the driving noise does not directly act on all determining modes of the flow. For example, see \cite{Eyi96,Nov65,Sta88,VKF79}. When the noise is extremely degenerate, i.e., present in only a few directions, many existing papers assume the noise is additive, meaning it does not depend on the solutions. However, it is more natural to assume that the noise depends on the solutions, i.e., that the noise is multiplicative. The aim of this article is to prove the ergodic and mixing properties of 2D Navier-Stokes equations with extremely degenerate multiplicative Gaussian noise, using Malliavin calculus.

Since the main topic of this article is stochastic 2D Navier-Stokes equations, we provide a brief introduction to recent developments in this field.
In \cite{HM-2006,HM-2011},
 Hairer and Mattingly established   the exponential mixing for 2D Navier-Stokes equations on the torus and the sphere, provided that the random perturbation is an additive Gaussian noise containing only several Fourier modes.
 Shirikyan demonstrated the exponential mixing for the 2D Navier-Stokes system with space-time localized noise in \cite{shirikyan-asens2015}.
%
Shirikyan investigated the Navier-Stokes system with noise acting only on the boundary in \cite{Shi2021}.
Peng, Zhai and Zhang
\cite{PZZ24}  established  the ergodicity  for stochastic  2D Navier-Stokes equations driven by
a highly degenerate  pure jump L\'evy noise.
For 2D Navier-Stokes equations driven by multiplicative Gaussian noise, Odasso \cite{Od08}, Ferrario, and Zanella \cite{FZ24}  demonstrated   the exponential mixing property, requiring that the number of noises is sufficiently large.
As far as we know, there is no literature on the ergodicity of stochastic 2D Navier-Stokes equations driven by multiplicative noise when the noise is extremely degenerate, i.e., it appears only in several directions and depends on the solution. This paper provides an affirmative answer to this problem


Not only for the stochastic  2D Navier-Stokes equations,
in many papers concerning    the ergodicity of   stochastic partial equations driven by multiplicative
noise,  they  also  require that
 the noise exists in sufficiently many directions.
For example, this is the case for three-dimensional (3D) stochastic Navier-Stokes equations \cite{Od07}, the 2D stochastic quasi-geostrophic equation \cite{RZZ15}, and the stochastic model of two-dimensional second-grade fluids \cite{WZZ16}, among others.
%

In addition to the literature mentioned above, there are also several other studies that explore the ergodicity  properties of stochastic partial differential equations driven by multiplicative noise. For instance, under specific growth conditions on the coefficients associated with the noise, Brze\'zniak, Ferrario, Maurelli, and Zanella \cite{BFMZ24} demonstrated that all solutions of the stochastic nonlinear Schr\"odinger equation will converge to $0$, with $\delta_0$ being the unique invariant measure. Furthermore, for solutions of stochastic conservation laws with Dirichlet boundary conditions driven by multiplicative noise, Dong, Zhang (Rangrang), and Zhang (Tusheng) \cite{DZZ23} established  the existence of a unique invariant measure and also demonstrate a polynomial mixing property.  Lastly, in the context of stochastic porous media equations with nonlinear multiplicative noise on bounded domains with Dirichlet boundary, Dareiotis, Gess, and Tsatsoulis \cite{DGT20} proved  the existence and uniqueness of invariant measures, accompanied by an optimal rate of mixing.


In this paper, we   study the ergodic  and   exponential mixing properties for the
 2D Navier-Stokes equations with degenerate  multiplicative  Gaussian noises.
 Recall that   the Navier-Stokes equations on torus  $\mT^2=[-\pi,\pi]^2$ are given by
 \begin{eqnarray}
 \label{16-4}
\partial_t u + (u,\cdot \nabla)u = \nu \Delta  u -  \nabla p + \xi,
\text{ div }u = 0,
 \end{eqnarray}
 where
$u=u(x, t)\in \mR^2$ denotes the value of the velocity field at time
 $t$ and position
$x$, $p(x, t)$ denotes the pressure, $\nu>0$ is the viscosity  and $\xi=\xi(x, t)$
 is an external force field acting on
the fluid.
We rewrite the NS system \eqref{16-4}  in the vorticity formulation:
   \begin{equation}
  \label{0.1}
    \p_t w -\nu \Delta w-B( \cK w,w)=\eta ,
  \end{equation}
  where $\eta=\eta(t,x)$ is   an external random force which will
   be described in (\ref{0.4}) below,
  $w=\nabla\wedge u$,   $B(v,w)=-( v \cdot \nabla ) w$,
  and  $\cK$ is the Biot--Savart operator
which
 is defined in Fourier space
by $(\cK w)_k =-iw_kk^{\bot}/|k|^2$,  where $k^\bot=(-k_2,k_1).$
 By $w_k$, it  means
the scalar product of $w$ with $\frac{1}{2\pi} e^{i k \cdot x}$.
  We consider this system  (\ref{0.1}) in the Hilbert space as follows
\begin{equation}\label{0.2}
H=\Big\{w\in L^2(\mT^2,\R):  \int_{\mT^2} w(x)\dif x=0 \Big\}
\end{equation}
   endowed    with the   $L^2$-scalar product $\lag \cdot, \cdot\rag$ and the corresponding $L^2$-norm $\|\cdot\|.$

 In order to describe the noises $\eta$,   we introduce the following notations.
 Denote
  $$
e_k(x)=
\begin{cases}
 \sin \lag k,x\rag& \text{if }k_1>0\text{ or }k_1=0,\, k_2>0, \\
  \cos\lag k,x\rag   & \text{if }k_1<0\text{ or } k_1=0,\, k_2<0,
  \quad k =(k_1,k_2).
 \end{cases}
$$
It is well-known that  $\{e_k: k\in \mZ^2\backslash (0,0)\}$ form a complete orthonormal basis  of $H.$
We   assume that $\eta$ is a  noise of the form
\begin{equation}\label{0.4}
\eta(t,x)=  \p_t \big(\sum_{k\in \cZ_0} q_k(w_t)W_k(t) e_k(x)\big),
\end{equation}
where $\cZ_0 \subset \Z^2_*:=\mZ^2\backslash (0,0)$ is a  finite set,
 $\{ q_k, {k\in \cZ_0}\} $ are continuous  functions  from $H$ to $\mR$,     $W=(W_k)_{k\in \cZ_0}   $
 is a
$|\cZ_0|$-dimensional standard   Brownian motion
on a filtered probability space~$(\Omega,\FF,  \{\FF_t\}_{t\geq 0}, \pP)$    satisfying
        the usual conditions (see   Definition 2.25 in~\cite{KS1991}).
        For the convenient of writing,  we  always denote $|\cZ_0|$ by $d.$

This work is motivated by the paper \cite{HM-2006}, in which Martin Hairer and Jonathan C. Mattingly considered the case of additive noise, where the noise coefficients are constant, i.e., $q_k\equiv$ constant, $\forall k\in \cZ_0.$
Leveraging the asymptotically strong Feller property that they discovered, they proved the uniqueness and existence of an invariant measure for the semigroup generated by the solution to $(\ref{0.1})$, even when the noise is present in only four directions.


To formulate our main result, let us recall that a set $\cZ_0 \subset \Z^2_*$
is a generator if any element of~$\Z^2$
is a finite linear
combination of elements of~$\cZ_0$ with integer coefficients.
Throughout the whole paper,  we    assume that the following  two  conditions ~are ~satisfied.
\begin{condition}
\label{16-5}
 The set  $\cZ_0 \subset \Z^2_*$ in \eqref{0.4} is a
 finite symmetric (i.e.,  $-\cZ_0=\cZ_0$)  generator
that contains at least two non-parallel vectors~$m$
and~$n$
such that $|m|\neq |n|$.
\end{condition}
This is   the   condition under which   the  ergodicity
of the NS system is  established  in~\cite{HM-2006, HM-2011}
in the case of a white-in-time noise and in \cite{KNS-2018} in
 the case of a bounded   noise. The set
$$
\cZ_0=\{(1,0), \, (-1,0), \, (1,1), \, (-1,-1)\}\subset \Z^2_*
$$ is an   example satisfying this condition.

\begin{condition}
\label{16-2}
There exists a constant $\aleph >0$ such that  the followings hold.
\begin{itemize}
\item[(\romannumeral1)]
for any $i\in \cZ_0$ and   $w\in H$,
\begin{eqnarray*}
  |q_i(w)|\in (0,\aleph].
\end{eqnarray*}
\item[(\romannumeral2)]
for  any  $i\in \cZ_0$, the function $q_i(\cdot)$ is  first and second Fr\'echet
differentiable and the following bounds hold
for any $u,v,w\in H$
\begin{eqnarray*}
  && \|Dq_i(w)v\|\leq  \aleph \|v\|,
  \\ &&
  \|D^2q_i(w)(u,v)\|\leq  \aleph \|u\|\|v\|.
\end{eqnarray*}
\end{itemize}
\end{condition}

\begin{remark}
Under the  Condition \ref{16-2}, it well-known that
  the equation (\ref{0.1}) has  a unique solution in the space
  $C([0,T],H)\cap L^2([0,T],H^1)$   for any $T>0,$
  where  $H^1$ is the    usual    Sobolev space of order~$1.$
  For example, please refer to    \cite[Theorem 4.7]{GM10}   etc.
\end{remark}

  Before we state the main theorem, we introduce some  notations.
We denote by $P_t(w_0,\cdot)$ the transition probabilities for equation (\ref{0.1}), i.e,
\begin{eqnarray*}
  P_t(w_0,A)=\mP(w(t)\in A\big|w(0)=w_0)
\end{eqnarray*}
for every Borel set $A\subseteq H.$ For every Borel set $A\subseteq H,$ $\Phi:H\rightarrow \mR$ and probability measure $\mu$ on $H,$ define
\begin{equation*}
P_t \Phi (w_0)=\int_{H} \Phi (w)P_t(w_0,\dif w), \quad  P_t^*\mu(A)=\int_{H}P_t(w_0,A)\mu(\dif w_0).
\end{equation*}
For  $\eta>0$ and  $\Phi:H\rightarrow \mR,$  the norm of $\Phi$  is defined  by
\begin{eqnarray*}
  \|\Phi\|_\eta=\sup_{v\in H}\[e^{-\eta\|v\|^2}(|\Phi(v)|+\|D \Phi(v)\|)\].
\end{eqnarray*}
Let
\begin{eqnarray*}
  \cO_\eta=\{\Phi \in C^1(H), \|\Phi\|_\eta<\infty\}.
\end{eqnarray*}
Our aim is to prove the following theorem.
\begin{theorem}
  Consider the  2D Navier-Stokes equation  (\ref{0.1}) with a
  multiplicative  Gaussian noise
  (\ref{0.4}).
  Under the Condition \ref{16-5} and Condition {\ref{16-2}},
  the followings  hold.
  \begin{itemize}
    \item[(\romannumeral1)]
   There exists a unique invariant measure $\mu_*$
associated to (\ref{0.1}), i.e., $\mu_*$ is a probability measure on $H$
 such that $P_t^*\mu_*=\mu_*$ for every $t\geq 0.$
    \item[(\romannumeral2)]
    There exists a $\eta_0>0$ such that for any $\eta\in (0,\eta_0],\Phi\in \cO_\eta$ and  $t\geq 0,$
one has
    \begin{eqnarray}
    \label{p1214-1}
      \big|\mE  \Phi(w_t)-\int_{H} \Phi(v)\mu_*(\dif v)\big|
      &\leq & Ce^{-\gamma t+\eta \|w_0\|^2}\|\Phi\|_\eta,
    \end{eqnarray}
    where $C=C(\eta),\gamma=\gamma(\eta)$ are positive constants depending on $\eta$ and $\nu,\aleph, d=|\cZ_0|.$
  \end{itemize}
\end{theorem}
\begin{proof}
 In order to prove this theorem,
  we only   need to verify the Assumptions 4--6 in \cite{HM-2008}
 for   $V(w)=\exp{(\eta \|w\|^2)}$ when $\eta$ is small enough,
  i.e,  a type of Lyapunov structure on the solution $w_t,$
a type of gradient inequality for $P_t$,
and  a  weak form of irreducibility.
By \cite[Theorems 3.4,3.6]{HM-2008},
these three assumptions   not only   imply
the existence  and  uniqueness of invariant measure,
but also imply  the  exponential  mixing property (\ref{p1214-1}).
See also \cite[Section 7]{FGRT-2015}.

The verification  of   a type of Lyapunov structure   is  essentially  the same as that in  \cite{HM-2008}.
The proof of  the    weak form of irreducibility  is demonstrated in Section \ref{123}.
The  type of gradient inequality for $P_t$  is often referred as    asymptotically strong Feller property, and we will give  a  proof  of it in Proposition \ref{1-7} below.
\end{proof}

\begin{proposition}\label{1-7}
(Asymptotically strong Feller property)
For  all $\eta,\gamma>0,w_0,\xi\in H$  and  $f:H \rightarrow \mR$ with $\|f\|_\infty+\|D  f\|_\infty\leq 1 $, it holds  that
\begin{eqnarray*}
 \big|D_\xi P_tf(w_0)\big|\leq C_{\eta,\gamma}\exp\{\eta\|w_0\|^2\} (\|f\|_\infty+\|D    f\|_\infty e^{-\gamma t})\|\xi\|,~\forall t\geq 0,
\end{eqnarray*}
where $C_{\eta,\gamma}$  is a   positive constant  depending on $\eta,\gamma$ and $\nu,\aleph, d.$
\end{proposition}
We will give a proof of this proposition in Section \ref{p4-1}.

\begin{example}

 For any $k\in \cZ_0,$ if we set
  \begin{eqnarray*}
    q_k(w)=f_k(\langle w,e_i\rangle),
  \end{eqnarray*}
  where  $i=i(k)\in \mZ^2_*$   and  $f_k:\mR \rightarrow \mR$  is a 2-order differentiable function such that
  for some $\aleph>0,$
  \begin{eqnarray*}
    && |f_k(x)|\in (0,\aleph],
   \quad
    |f_k'(x)|\leq \aleph,
    \quad  |f_k^{(2) }(x)|\leq \aleph,\quad \forall x\in \mR,
  \end{eqnarray*}
  then the Condition \ref{16-2} holds.
\end{example}
\begin{example}
For any $k\in \cZ_0,$ if we set
  \begin{eqnarray*}
    q_k(w)=f_k(\| w\|),
  \end{eqnarray*}
  where $f_k:[0,\infty) \rightarrow \mR$  is a 2-order differentiable function such that
  for some $\aleph>0$
  \begin{eqnarray*}
    && |f_k(x)| \in (0,\aleph],
    \quad
    |f_k'(x)|\leq \aleph   \min\{ |x|,1\},
   \quad  |f_k^{(2) }(x)|\leq \aleph, \quad\forall x\geq 0,
  \end{eqnarray*}
  then the Condition \ref{16-2} holds.
 \end{example}
 \begin{proof}
   By direct calculations, we have
   \begin{align*}
     Dq_k(w)u &= f_k'(\|w\|)\frac{\langle w,u\rangle }{\|w\|},
     \\ D^2q_k(w)(u,v) &=
     f_k{''}(\|w\|)\frac{\langle w,v\rangle }{\|w\|}
     \frac{\langle w,u\rangle }{\|w\|}
     +f_k'(\|w\|)
     \frac{\langle v,u\rangle \|w\|^2-\langle
     w,u\rangle\langle w,v\rangle}{\|w\|^3}.
   \end{align*}
   Then,    we obtain   the desired results.
\end{proof}

Because of  the intensity of the noise depends on the solutions,
 we need to introduce a  set of new  ideas and techniques to establish the exponential mixing property
in comparison to the case of additive  noise (see \cite{HM-2006}  \cite{HM-2011}\cite{FGRT-2015}\cite{PZZ24} etc.).  Now we
explain  them in details.

 Let $\cM_{0,t}$ be the Malliavin operator related to   $w_t$
    and
   $\cS_{\alpha,N}$ be some   subspace of $H$
   (for
rigorous  definitions  of $\cM_{0,t}$ and  $\cS_{\alpha,N}$,    please see (\ref{28-5}) and Proposition \ref{3-10} in this  article).
To obtain the ergodicity via Malliavin calculus,
in the classical literature,
one key  ingredient    is   to show
   \begin{align}
   \label{1818-2}
     \mP(\inf_{\xi \in\cS_{\alpha,N} } \langle \cM_{0,t}\xi,\xi\rangle <\eps)\leq C(\|w_0\|) r(\eps), \forall \eps\in (0,1)\text{ and } w_0\in H,
   \end{align}
   where  $\|w_0 \|$ denotes the $L^2$ norm of $w_0$,   $C$ is some locally  bounded  function  from $[0,\infty)$ to $ [0,\infty)$  and  $r$ is a function on $(0,1)$ with $\lim_{\eps\rightarrow 0}r(\eps)=0$.
      In the
      process of proving  (\ref{1818-2}),
   the properties of  Gaussian polynomials
   (see, e.g.,   \cite[Theorem 7.1]{HM-2011})  play very essential roles in   the estimate of
    probability
   $\mP(\inf_{\xi \in\cS_{\alpha,N} } \langle \cM_{0,t}\xi,\xi\rangle <\eps).$
   Recenty, Peng, Zhai and Zhang \cite{PZZ24}
   established  the ergodicity  for stochastic  2D Navier-Stokes equations driven by
a highly degenerate  pure jump    L\'evy noise
requiring   a little weaker  version of  (\ref{1818-2}).
In a word,  in order to  obtain the ergodicity via Malliavin calculus,
we need to prove (\ref{1818-2}) or  some  little weaker  version of  (\ref{1818-2}).

  For any $\xi \in H,$ by direct calculations, one easily sees that
  \begin{align*}
\langle \cM_{0,t}\xi,\xi \rangle=\sum_{j\in \cZ_0}\int_{0}^tq_j^2(w_r)\langle J_{r,t}e_j,\xi\rangle^2\dif r,
\end{align*}
  where   $J_{s,t}\xi$ is the derivative of $w_t$ with respect to initial value $w_s$ in the direction of $\xi.$
  For rigorous definition of $J_{s,t},$ please to see (\ref{10-1}) below.
  Obviously, $J_{s,t}:H\rightarrow H$ is a linear operator,  let~$K_{s,t}$  be the   adjoint of $J_{s,t}$.
  In all the existing  papers,  in order to prove  (\ref{1818-2}) or some weak form of  (\ref{1818-2}),   they need to use the
 equation that   $\varrho_s=K_{s,t}\xi, s\in [0,t]$ is satisfied.
Since they  are concerned with  the   additive noise,
when we  take the solutions    $w_t, t\geq 0$
 as fixed,
they are
deterministic   differential equations and don't have  It\^o's integral  terms.
In the case of    multiplicative  Gaussian noises,
formally,
there will have a  non-adapted stochastic  integral
in   the
 the equation that    $\varrho_s$ is satisfied.
It is not easy  to     interpret this  non-adapted integral.
Unlike the previous literature,  we directly  use the equation which
$x_t=J_{s,t}\xi, t\in [s,\infty)$ is  satisfied  to prove some weak form of
(\ref{1818-2}),   and  eventually avoid this   problem.

 The rest of the paper is organized as follows.
    In Section \ref{p2-1}, we provide some estimates for the solution $w_t$ and introduce the essential ingredients of the Malliavin  calculus for the solution. Moreover, we give all the necessary estimates associated with the Malliavin matrix.
    Section \ref{p3-1} is devoted to  the proof of the invertibility of the  Malliavin matrix  of  the solution  $w_t$ which plays a key role in the proof of Proposition \ref{1-7}.
The proof  of  Proposition \ref{1-7} is demonstrated   in Section \ref{p4-1}. Some of the technical proofs are put in the Appendix.

\section{Preliminaries}
\label{p2-1}
\subsection{Notation}

In this paper, we use the following  notation.

 \smallskip
\noindent
 Let $\Z^2_*:=\mZ^2\backslash (0,0).$  Let ${H}_N=span\{e_j: j\in \mZ_*^2  \text{ and } |j|\leq N\}.$
${P}_N$  denotes the orthogonal projections from $H$ onto ${H}_N$.
      Define   ${Q}_Nu:=u-{P}_Nu, \forall u\in H.$

\smallskip
\noindent
For $\alpha\in \mR$ and a smooth function $w\in H$, we define the norm $\|w\|_\alpha$ by
\begin{eqnarray*}
\|w\|_\alpha^2=\sum_{k\in \mZ_*^2}|k|^{2\alpha}w_k^2,
\end{eqnarray*}
where $w_k$ denotes the Fourier mode with wavenumber $k.$
 When $\alpha=0,$    we also    denote  this norm $\|\cdot \|_\alpha$ by $\|\cdot \|.$
 $H^m= H^m(\mT^2, \R)\cap H$, where $H^m(\mT^2, \R)$ is the    usual    Sobolev space of order~$m\ge1$. We endow the space  $H^m$ with the Sobolev norm~$\|\cdot \|_m$.

 \smallskip
\noindent
For any $(s_1,s_2,s_3)\in \mR_{+}^3$ with $\sum_{i=1}^3s_i\geq 1$
and $(s_1,s_2,s_3)\neq (1,0,0),(0,1,0),(0,0,1)$, the following relations will be used frequently in this paper(c.f. \cite{CF88}):
\begin{eqnarray}
\label{44-1}
\begin{split}
\langle B(u,v),w\rangle& = -\langle B(u,w),v\rangle, \quad \text{if }\nabla \cdot u=0,
\\
  \big|  \langle B(u,v),w\rangle\big| &\leq
C\|u\|_{s_1}\|v\|_{1+s_2}\|w\|_{s_3},
\\ \|\cK u\|_{\alpha} & = \|u\|_{\alpha-1},
\\  \|w\|_{1/2}^2 &\leq  \|w\|\|w\|_1.
\end{split}
\end{eqnarray}

\smallskip
\noindent
$L^\ty(H)$  is the space of bounded Borel-measurable functions $f:H\to\R$   with the norm $\|f\|_\ty=\sup_{w\in H}|f(w)|$.~$C_b(H)$ is the space of continuous functions. ~$C^1_b(H)$ is the space of functions~$f\in C_b(H)$ that are continuously Fr\'echet differentiable with bounded derivatives.
The Fr\'echet derivative of $f$ at  point  $w$ is denoted by $D f (w)$
and   we usually also  write $D f (w)\xi$ as $D_\xi f(w)$ for any $\xi\in H.$
Taken as a function of $w$, if  $g(w):=D_\xi f(w)$ is also Fr\'echet differentiable, then  for any $\zeta\in H,$
we usually denote  $D_\zeta g(w)$  by $D^2f(w)(\xi,\zeta).$

\smallskip
\noindent
$\cL(X,Y)$ is the space of bounded linear operators from Banach spaces $X$ into Banach space $Y$ endowed with the natural norm $\|\cdot\|_{\cL(X,Y)}$.
If there are no confusions, we always write  the operator  norm  $\|\cdot\|_{\cL(X,Y)}$
as $\|\cdot\|$.

\smallskip
\noindent
Let   $d=|\cZ_0|$ and denote the canonical basis of $\mR^d$
by
$\{\theta_j\}_{j\in \cZ_0}.$


  \smallskip
\noindent
Since there are many constants appearing in the proof, we adopt  the following convention. Without otherwise specified,  the letters $C,C_1,C_2,\cdots$ are  always  used to denote unessential constants that  may
change from line to line  and  implicitly  depend on
 $\nu,\aleph$ and $d.$
 Also, we usually do not explicitly indicate the dependencies on the parameters  $\nu,\aleph$ and $d$ on every occasion.

\subsection{Priori estimates  on the solutions}

The following lemma  gather    some standard  priori estimates for the solution  of the stochastic Navier--Stokes system.

  \begin{lemma}\label{L:5.1}
   Assume $w_t$ is a solution of  \eqref{0.1}  with initial value       $w_0\in H$.    Then we have the following        estimates.
\begin{description}
 \item[(\romannumeral1)]
  There  are positive      numbers   $\eta_0=\eta_0(\nu,d,\aleph), C=C(\nu,d,\aleph)$
      such that
  for any $\eta\in(0,\eta_0]$,
 \begin{eqnarray*}
 \begin{split}
 \nonumber    & \mE \Big[\exp\{\eta \sup_{t\geq s}\big(\|w_t\|^2+\nu \int_s^t \|w_r\|_1^2\dif s-d \aleph^2  (t-s)\big)\}\big |\cF_s\Big]
   \leq    \exp\{\eta \|w_s\|^2\},
   \\
       & \mE\exp\{\eta \|w_t\|^2\} \leq C \exp\{\eta e^{-\nu t}\|w_0\|^2\}.
    \end{split}
\end{eqnarray*}
\item[(\romannumeral2)]  There exists  positive constants $\eta_1,a,\gamma>0$  depending only on $ \nu, d$ and  $\aleph$  such that
\begin{eqnarray*}
  \mE \exp\{\eta \sum_{k=0}^n\|w_k\|^2-\gamma n\}\leq \exp\{a\eta\|w_0\|^2 \}
\end{eqnarray*}
holds for every $n>0,$ every $\eta\in (0,\eta_1]$.
\item[(\romannumeral3)]For any $\eta, p,T>0$ and $k\in \{1,2,3,4,5,6\}$,  one has
  \begin{eqnarray*}
   &&  \mE \sup_{t\in [T/2,T]} \|w_t\|_k^p\leq C(\eta,p,T)\exp\{\eta\|w_0\|^2\},
  \end{eqnarray*}
  where $ C(\eta,p,T)$ is a constant depending on $ \eta,p,T$ and
  $\nu,d,\aleph.$
\end{description}
\end{lemma}
\begin{proof}
The proof of\text{ (\romannumeral1)(\romannumeral2) }are the same as that in \cite[Lemma 4.10 and Lemma A.1]{HM-2006}, so we omit the details.
Since $\sum_{j\in \cZ_0}\int_0^t q_j(w_s)W_j(s)e_j$
is   spatially smooth,
the proof of \text{(\romannumeral3)} is similarly  as that  in
\cite[Proposition 2.4.12]{KS-book}.
The proof is complete.

\end{proof}

\subsection{Elements of Malliavin calculus}
Since the main technique of this  article is  Malliavin calculus,
we give a detailed statement of Malliavin calculus and  demonstrate  some   related estimates
in this subsection.

Let   $d=|\cZ_0|$ and denote the canonical basis of $\mR^d$
by
$\{\theta_j\}_{j\in \cZ_0}.$
For any $w\in H,$ we define the linear operator {$Q=Q(w): \mR^d\rightarrow H$ in the following way: for any $z=\sum_{j\in \cZ_0}z_j\theta_j \in \mR^d$,
  \begin{eqnarray*}
   Qz= Q(w)z=\sum_{j\in \cZ_0 }   q_{j} (w)z_j e_{j}.
  \end{eqnarray*}
  The adjoint of $Q(w)$ is given by  $Q^*(w): H\rightarrow \mR^d$ in the following way:
  \begin{eqnarray*}
   Q^*(w)\xi=\big(q_{j}(w)\langle \xi, e_{j}\rangle \big)_{j\in \cZ_0} \in \mR^d   \text{ for } \xi\in H.
  \end{eqnarray*}

  For any $0\le s\le t$ and $\xi\in  H$, let~$J_{s,t}\xi$  be the   solution of the linearised problem:
\begin{align}
\label{10-1}
\begin{split}
	\dif   J_{s,t}\xi- \nu \Delta J_{s,t}\xi \dif t
-\tilde{B}(w_t,J_{s,t}\xi) \dif t &=
\sum_{j\in \cZ_0}  \big(Dq_j(w_t)J_{s,t}\xi\big)  e_j \dif W_j(t),
  \\  J_{s,s}\xi&=\xi,
  \end{split}
\end{align}
where $\tilde{B}(w,v)=B(\cK w,v)+B(\cK v,w)$.

  For any $0\le t\le T$ and $\xi\in  H$,
 let~$K_{t,T}$  be the   adjoint of $J_{t,T}$, i.e, $K_{t,T}\xi$ is an elemnt in $H$ such that $\langle K_{t,T}\xi,\phi\rangle=\langle \xi,J_{t,T}\phi\rangle, \forall \phi\in H$.
 For  given  $t>0$, let $\Phi=\Phi(t,W)$ be a ${\cal F}_{t}$-measurable random variable.  For $v\in L^2([0,t];\R^d),$
the Malliavin derivative of~$\Phi$ in the direction~$v$ is defined by
$$
\cD^v\Phi(t,W)=\lim_{\eps \to  0}\frac{1}{\e}
  \left(\Phi(t,w_0,W+\eps \int_0^\cdot v\dd  s)-\Phi(t,w_0,W)\right), \quad
$$
  where the limit   holds almost surely (e.g., see the book~\cite{nualart2006} for finite-dimensional setting or the papers~\cite{MP-2006, HM-2006, HM-2011, FGRT-2015}   for   Hilbert space~case).
  Then, $\cD^v w_s$ satisfies the following equation:
\begin{eqnarray*}
\begin{split}
&\dif  \cD^v w_s - \nu \Delta \cD^v w_s \dif s  - \tilde{B}(w_s,\cD^v w_s)\dif s
\\ & =
 \sum_{j\in \cZ_0} \big( Dq_j(w_s) \cD^v w_s\big) e_j  \dif W_j(s)+ Q(w_s)v(s) \dif s.
\end{split}
\end{eqnarray*}
  By the Riesz representation theorem, there is a linear operator $  \cD:L^2(\Omega,   H)\to L^2(\Omega; L^2([0,t ];\R^d)\otimes  H)$ such that
\begin{equation*}
	  \cD^v w_t  =\lag    \cD w,v \rag_{L^2([0,t];\R^d)},~\forall v\in L^2([0,t];\R^d) .
\end{equation*}
 By the formula of constant variation, we have
\begin{equation}\label{2.3}
	  \cD^v w_t  =\int_0^t J_{r,t}Q(w_r)v(r) \dif r.
\end{equation}
Therefore,   we also have
  \begin{eqnarray*}
    \cD_u^jw_t=J_{u,t}  Q(w_u) \theta_j, ~\forall u\in [0,t],
  \end{eqnarray*}
  where  $ \cD_u^j$ denotes the Malliavin derivative with respect to the $j^{th}$ component of the noise at time $u.$
By (\ref{2.3}),  we get
\begin{equation}\label{5-3}
	\cD^v w_t=\cA_{0,t}v,
\end{equation}
where   {{$\cA_{s,t}:L^2([s,t],\mR^d)\to   H$}}  is  the random operator defined  by
\begin{equation}
\label{5-4}
\cA_{s,t}v=\int_{s}^{t}J_{r,t} Q(w_r)v(r) \dd  r, \quad 0\le s\le t<\infty.
\end{equation}
The random operator   {{$\mathcal{A}^*_{s,t}:H\rightarrow L^2([s,t],\mR^d)  $}}     is given by
$$
(\cA_{s,t}^* \xi)(r)= Q^*(w_r) J_{r,t}^*\xi  , \quad  \xi\in  H,\, r\in [s,t].
$$

For any  $0\le  s< t,$     the Malliavin operator    is defined by
$$
\cM_{s,t}=\cA_{s,t}\cA_{s,t}^*:H\rightarrow H
$$
and then
\begin{align}
\label{28-5}
\cM_{s,t}\xi =\int_{s}^tJ_{r,t} Q(w_r)  Q^*(w_r)  K_{r,t}\xi \dif r,~\xi\in H.
\end{align}
It is a
  non-negative self-adjoint operator, so its regularisation~$\cM_{s,t}+\beta\I$
is invertible for any  $\beta>0$,  where    $\I$ is the identity operator.
By direct calculations, for any $T>0$ and $\xi\in H,$  we see that
\begin{eqnarray*}
  \langle \cM_{0,T}\xi, \xi  \rangle
  =\sum_{j\in \cZ_0}\int_0^T q_j^2(w_r) \langle J_{r,T}e_j,\xi\rangle^2\dif r.
\end{eqnarray*}

\begin{lemma}
\label{2-3}
 For any $n\in \mN,\kappa,T>0$ and  $0\leq s \leq T,$   one has
  \begin{eqnarray*}
  \sup_{\|\xi\|\leq 1}\mE\Big[  \Big( \sup_{t\in [s,T]} \|J_{s,t}\xi\|^2+\int_s^T \|J_{s,r}\xi\|_1^2\dif r\Big)^n \big|\cF_s \Big]   \leq    C_{\kappa,n,T}   e^{\kappa \|w_s\|^2 },
  \end{eqnarray*}
     where $C_{\kappa,n,T}$  is  a    constant   depending on
  $\kappa,n,T$ and
 $\nu,\aleph,d$\footnote{Obviously, we only need to prove this lemma for $\kappa$ small enough, i.e.,  $\kappa\leq \kappa_0(\nu,\aleph,d)$. }.
\end{lemma}
\begin{proof}
By (\ref{44-1}) and Young's inequality, it holds that
\begin{eqnarray}
\label{0927-2}
\begin{split}
 &  \big| \langle B(\cK J_{s,t}\xi,w_t ),J_{s,t}\xi \rangle\big|
  \leq C\|w_t\|_1 \|J_{s,t}\xi\|_{1/2}\|J_{s,t}\xi\|
  \\ & \leq  \frac{\nu}{4}\|J_{s,t}\xi\|_1^2+ C\|w_t\|_{1}^{4/3}\| J_{s,t}\xi\|^2,
  \end{split}
\end{eqnarray}
where $C=C(\nu).$  So,  using  It\^o's formula to $\|J_{s,t}\xi\|^{2n}$, one arrives at that \begin{eqnarray}
 \nonumber  && \dif \|J_{s,t}\xi\|^{2n}\leq 2n\|J_{s,t}\xi\|^{2n-2}
  \langle J_{s,t}\xi, \dif J_{s,t}\xi \rangle
 \\ \nonumber  && \quad\quad  +2n(n-1)\sum_{j\in \cZ_0}  \|J_{s,t}\xi\|^{2n-4} |Dq_j(w_t)J_{s,t}\xi|^2  \langle J_{s,t}\xi,   e_j  \rangle^2\dif t
  \\ \nonumber &&\quad\quad +n \|J_{s,t}\xi\|^{2n-2}\sum_{j\in \cZ_0} | Dq_j(w_t)J_{s,t}\xi|^2\dif t
  \\ \nonumber &&\leq  2n\|J_{s,t}\xi\|^{2n-2}\big[-\nu \|J_{s,t}\xi\|^2_1\dif t+C\|w_t\|_{1}^{4/3}\| J_{s,t}\xi\|^2\dif t\big]
  +C_n \|J_{s,t}\xi\|^{2n}\dif t
  \\ \nonumber &&\quad\quad+2n\|J_{s,t}\xi\|^{2n-2} \Big(\sum_{j\in \cZ_0}(Dq_j(w_t)J_{s,t}\xi)  \langle J_{s,t}\xi, e_j\rangle \dif W_j(t)\Big)
  \\ \nonumber &&\leq ( C_n\|w_t\|_{1}^{4/3}+C_n)\| J_{s,t}\xi\|^{2n}\dif t
   - 2n\nu \|J_{s,t}\xi\|^{2n-2} \|J_{s,t}\xi\|^2_1\dif t,
  \\ \label{2-4}
  &&\quad\quad
  +2n\|J_{s,t}\xi\|^{2n-2} \Big(\sum_{j\in \cZ_0}(Dq_j(w_t)J_{s,t}\xi)  \langle J_{s,t}\xi, e_j\rangle \dif W_j(t)\Big),
\end{eqnarray}
where $C_n$ is a constant depending on $n$ and $\nu,\aleph, d.$
Therefore, we obtain
\begin{eqnarray*}
  \mE \big[\|J_{s,t}\xi\|^{2n}e^{-C_n(t-s)-C_n\int_s^t \|w_r\|_{1}^{4/3}\dif r}\big|\cF_s \big] \leq 1.
\end{eqnarray*}
Combining the above with Lemma \ref{L:5.1}, Cauchy-Schwartz inequality   and
\begin{eqnarray}
\label{p1204-1}
C_n \|w_r\|_{1}^{4/3} \leq \frac{\kappa }{2}\|w_r\|_{1}^{2}  + C_{n,\kappa},\quad \forall \kappa>0,
\end{eqnarray}
 one arrives at that
  \begin{eqnarray}
  \label{2-5}
  \sup_{\|\xi\|\leq 1}\mE\big[  \|J_{s,t}\xi\|^{n} \big|\cF_s\big]  \leq    C_{\kappa,n,T}    e^{\kappa \|w_s\|^2},  ~\forall n\in \mN \text{ and } 0\leq s\leq t\leq T.
  \end{eqnarray}

Using (\ref{2-4}) with $n=1,$
 we get
 \begin{eqnarray*}
 &&  \sup_{t\in [s,T]}\|J_{s,t}\xi\|^2 + \nu \int_s^T\|J_{s,r}\xi\|_1^2\dif r\leq \|\xi\|^2+
 C  \int_s^T  (\|w_r\|_{1}^{4/3}+1)\| J_{s,r}\xi\|^{2}\dif r
 \\ && \quad\quad\quad\quad\quad\quad\quad\quad  +\sum_{j\in \cZ_0}
\sup_{t\in [s,T]}
  \int_s^t (Dq_j(w_r)J_{s,r}\xi)
  \langle J_{s,r}\xi, e_j\rangle
 \dif W_j(r).
 \end{eqnarray*}
 Thus,  for any $n\geq 2$  and $\xi$ with $\| \xi\| \leq 1,$
 by (\ref{2-5}),   Lemma \ref{L:5.1},  H\"older's inequality  and
 Burkholder-Davis-Gundy  inequality,  it holds that
 \begin{eqnarray*}
&& \mE\Big[  \big( \sup_{t\in [s,T]} \|J_{s,t}\xi\|^2+\int_s^T \|J_{s,r}\xi\|_1^2\dif r\big)^n \big|\cF_s \Big]
\\ &&\leq  C_n+ C_{n,T}\big( \mE \sup_{t\in [s,T]} \|J_{s,t}\xi\|^{4n}\big|\cF_s\big)^{1/2}\Big(\mE\big(1+\int_s^T \|w_r\|_1^{4/3}\dif r)^{2n })\big|\cF_s \Big)^{1/2}
\\ &&\quad\quad+C_n \mE \Big( \big( \int_s^T \|J_{s,r}\xi\|^{4} \dif r  \big)^{n/2} \big| \cF_s\Big)
\\ && \leq    C_{\kappa,n,T}   e^{\kappa \|w_s\|^2 }.
 \end{eqnarray*}
 The proof is complete.
%
\end{proof}

\begin{lemma}
\label{19-1}
For any $n\geq 1,T>0, s\in [0,T],\kappa>0$ and $H$-valued  random variable $\xi\in \cF_s$, it   holds that
  \begin{eqnarray}
  \label{30-1}
    \mE \big[\sup_{t\in [s,T]}\|J_{s,t}\xi\|_1^{n}\big|\cF_s\big]
    \leq  C_{\kappa,n,T}\|\xi\|_1^{n}\exp\{\kappa \|w_s\|^2 \},
  \end{eqnarray}
   where $C_{\kappa,n,T}$  is  a    constant   depending on
  $\kappa,n,T$ and
 $\nu,\aleph,d.$
  \end{lemma}
  \begin{proof}
  For any $\kappa>0$ and  $n\geq 1,$ let   $\eps=\eps(\kappa,n)=\min\{\frac{\nu}{8},\frac{\kappa }{8n}\}$.
    By (\ref{44-1}),  it holds that
  \begin{eqnarray*}
   &&  \big| \langle  B(\cK w_t,J_{s,t}\xi), \Delta J_{s,t}\xi\rangle\big| \leq
    C \| J_{s,t}\xi\|_2 \|J_{s,t}\xi\|_1  \|w_t\|_{1/2}
    \\ &&\leq  \eps  \| J_{s,t}\xi\|_2^2+C_{\kappa,n} \|J_{s,t}\xi\|_1^2  \|w_t\|_{1/2}^2
    \\ &&\leq  \eps  \| J_{s,t}\xi\|_2^2+\eps  \|J_{s,t}\xi\|_1^2  \|w_t\|_{1}^2+C_{\kappa,n}   \|J_{s,t}\xi\|_1^2  \|w_t\|^2
  \end{eqnarray*}
  and
  \begin{eqnarray*}
   && \big| \langle  B(\cK J_{s,t}\xi,w_t), \Delta J_{s,t}\xi\rangle\big| \leq
    C \| J_{s,t}\xi\|_2 \|J_{s,t}\xi\|_{1/2}  \|w_t\|_{1}
    \\ && \leq  \eps  \| J_{s,t}\xi\|_2^2+C_{\kappa,n}  \|J_{s,t}\xi\|_1   \|J_{s,t}\xi\| \|w_t\|_{1}^2
      \\ && \leq  \eps  \| J_{s,t}\xi\|_2^2+\eps \|J_{s,t}\xi\|_1^2 \|w_t\|_{1}^2 +C_{\kappa,n}    \|J_{s,t}\xi\|^2  \|w_t\|_{1}^2,
  \end{eqnarray*}
    where $C_{\kappa,n}$  is  a    constant   depending on
  $\kappa,n$ and
 $\nu,\aleph,d.$
Using It\^o's formula to
  $\langle J_{s,t}\xi, -\Delta J_{s,t}\xi\rangle,$
    we arrive at that
  \begin{eqnarray*}
    && \dif  \langle J_{s,t}\xi, -\Delta J_{s,t}\xi\rangle
    =2
    \langle  -\Delta J_{s,t}\xi, \dif J_{s,t}\xi    \rangle
  +
    \big(\sum_{j\in \cZ_0}|j|^2 |Dq_j(w_t)J_{s,t}\xi|^2\big) \dif t
    \\ &&\leq -2\nu   \| J_{s,t}\xi \|_2^{2} \dif t
    \\ &&+ 2   \big[2 \eps  \| J_{s,t}\xi\|_2^2+ 2 \eps  \|J_{s,t}\xi\|_1^2  \|w_t\|_{1}^2+C_{\kappa,n}  \|J_{s,t}\xi\|_1^2  \|w_t\|^2+C_{\kappa,n}  \|J_{s,t}\xi\|^2  \|w_t\|_1^2
    \big]\dif t
    \\ &&+ C   \| J_{s,t}\xi \|^{2}\dif t  +2 \big\langle-\Delta J_{s,t}\xi, \sum_{j\in \cZ_0}\big(Dq_j(w_t) J_{s,t}\xi\big)e_j
    \big\rangle \dif W_j(t).
  \end{eqnarray*}
  Therefore, one has
  \begin{eqnarray*}
    && \|J_{s,t}\xi\|_1^2e^{-4 \eps \int_s^t\|w_r\|_1^2\dif r}
    \\ && \leq
    \|\xi\|_1^2+
    C_{\kappa,n} \int_s^t\big[  \|J_{s,r}\xi\|_1^2  \|w_r\|^2+  \|J_{s,r}\xi\|^2  \|w_r\|_1^2+\|J_{s,r}\xi\|^2 \big] \dif r
    \\ && \quad + 2 \int_s^t e^{-4 \eps \int_s^r\|w_u\|_1^2\dif u }
    \big\langle -\Delta J_{s,r}\xi, \sum_{j\in \cZ_0}\big(Dq_j(w_r) J_{s,r}\xi\big)e_j
    \big\rangle \dif W_j(r)
        \\ && \leq
    \|\xi\|_1^2+
    C_{\kappa,n} \sup_{r\in [s,t]}\|w_r\|^2\int_s^t \|J_{s,r}\xi\|_1^2  \dif r
    \\ && \quad +C_{\kappa,n} \sup_{r\in [s,t]} \|J_{s,r}\xi\|^2  \int_s^t  \|w_r\|_1^2\dif r+C_{\kappa,n} \int_s^t \|J_{s,r}\xi\|^2  \dif r
    \\ && \quad + 2 \int_s^t e^{-4 \eps \int_s^r\|w_u\|_1^2\dif u }
    \big\langle -\Delta J_{s,r}\xi, \sum_{j\in \cZ_0}\big(Dq_j(w_r) J_{s,r}\xi\big)e_j
    \big\rangle \dif W_j(r).
  \end{eqnarray*}
  Furthermore, by H\"older's inequality,
 Burkholder-Davis-Gundy  inequality  and  Lemmas  \ref{2-3}, \ref{L:5.1},  it holds that
  \begin{eqnarray*}
    \mE \big[\sup_{t\in [s,s+T]}\|J_{s,t}\xi\|_1^{2n} e^{-4 \eps n  \int_s^t\|w_r\|_1^2\dif r}   \big|\cF_s\big]
    \leq  C_{\kappa,n,T}\|\xi\|_1^{2n}\exp\{\kappa \|w_s\|^2  \}.
  \end{eqnarray*}
  The result of this lemma follows from the above inequality,
Cauchy-Schwartz inequality and Lemma \ref{L:5.1}.
The proof is complete.
  \end{proof}

Recall that $P_N$ is the orthogonal projection from $H$ into $H_N=span\{e_j:  j\in Z_{\ast}^2, |j|\leq N\}$ and $Q_N=I-P_N$.
For any $N\geq 1,t\geq 0$ and $\xi\in H,$  denote $\xi_t^h=Q_NJ_{0,t}\xi,
\xi_t^\iota =P_NJ_{0,t}\xi$ and $\xi_t=J_{0,t}\xi.$
\begin{lemma}
\label{1340-3}
 For any $n\in \mN, \kappa, T>0, t\in  [0,T],\xi\in H$ and $N>\max\{|j|: j\in \cZ_0\},$  one has
 \begin{eqnarray}
 \label{2-6}
   \sup_{\|\xi\|\leq 1}\mE \|\xi_t^h\|^{2n} \leq \exp\{- n\nu N^2 t\}+\frac{C_{\kappa,n,T}}{\sqrt{N} } \exp\{\kappa \|w_0\|^2\},
 \end{eqnarray}
 where $C_{\kappa,n,T}$  is  a    constant   depending on
  $\kappa,n,T$ and
 $\nu,\aleph,d.$  
\end{lemma}
\begin{proof}
Note that $\|\xi_t^h\|^2_1\geq N^2\|\xi_t^h\|^2$  and
\begin{eqnarray*}
  && \big| \langle B(\cK \xi_t,w_t),\xi_t^h\rangle\big|+\big|\langle B(\cK w_t,\xi_t),\xi_t^h\rangle\big| \leq
  C\|\xi_t^h\|_1 \|w_t\|_{1/2}  \|\xi_t\|
  \\ && \leq \frac{\nu}{4}\|\xi_t^h\|_1^2+C \|w_t\|_{1/2}^2  \|\xi_t\|^2.
\end{eqnarray*}
Applying the chain rule
 to $\|\xi_t^h\|^{2n},$  we conclude that
\begin{eqnarray*}
  \dif \|\xi_t^h\|^{2n}\leq  2 n \|\xi_t^h\|^{2n-2}\big[-\nu N^2 \|\xi_t^h\|^2  \dif t+C\|w_t\|_{1/2}^2  \|\xi_t\|^2\dif t\big].
\end{eqnarray*}
Therefore,
\begin{eqnarray}
 \nonumber  && \|\xi_t^h\|^{2n} \leq \exp\{-n\nu N^2 t\}\|\xi\|^2+
  C \int_0^t \exp\{-n \nu N^2 (t-s) \} \|w_s\|_{1/2}^2  \|\xi_s\|^{2n}\dif s.
  \\   \nonumber &&\leq \exp\{-n\nu N^2 t\}\|\xi\|^2
  \\  \nonumber && +C \sup_{s\in [0,T]}\|w_s\| \big(\int_0^t \exp\{-4n \nu N^2 (t-s) \}\dif s \big)^{1/4}\big(\int_0^t\|w_s\|_1^2\dif s\big)^{1/2}\big(\int_0^t \|\xi_s\|^{8n}\dif s\big)^{1/4}.
\end{eqnarray}
Then, with the help of   Lemma \ref{2-3} and Lemma \ref{L:5.1}, the above
inequality  implies the desired result (\ref{2-6}).

\end{proof}

\begin{lemma}
\label{2-1}
 For any $n\in \mN,\kappa>0$  and   $\xi\in H$ with $\xi_0^\iota=P_N\xi=0$, one has
 \begin{eqnarray}
 \label{30-3}
   \sup_{\|\xi\|\leq 1}\mE \|\xi_1^\iota \|^{n}\leq   \frac{C_{\kappa,n}}{N^{1/24} } \exp\{\kappa \|w_0\|^2\},
 \end{eqnarray}
    where $C_{\kappa,n}$  is  a    constant   depending on
  $\kappa,n$ and
 $\nu,\aleph,d.$
  Furthermore, combining the above  inequality with Lemma  \ref{1340-3}, for any  $n\in \mN,\kappa>0,\xi\in H$ and  $N>\max\{|j|: j\in \cZ_0\},$
  we have
    \begin{eqnarray}
    \label{30-2}
  &&  \mE \|J_{0,1}Q_N \xi\|^n  \leq  \frac{C_{\kappa,n}}{N^{1/24} } \exp\{\kappa \|w_0\|^2\}\|\xi\|^n.
  \end{eqnarray}
\end{lemma}
\begin{proof}
With similar arguments as that in (\ref{0927-2}),  and by (\ref{44-1}), we obtain
\begin{eqnarray*}
  \big| \langle B(\cK \xi_t^\iota,w_t),\xi_t^\iota\rangle\big|
  \leq \frac{\nu}{4} \|\xi_t^\iota \|_1^2+C\|w_t\|_1^{4/3}\|\xi_t^\iota\| ^2
\end{eqnarray*}
  and
\begin{eqnarray*}
&&  \big| \langle B(\cK w_t,\xi_t^h),\xi_t^\iota\rangle\big|
 +  \big| \langle B(\cK \xi_t^h, w_t),\xi_t^\iota\rangle\big|
 \\ && \leq  C \|\xi_t^\iota \|_1 \|w_t\|\|\xi_t^h \|_{1/2}
 \leq \frac{\nu}{4} \|\xi_t^\iota \|_1^2+C\|w_t\|^2 \|\xi_t^h \|_{1/2}^2.
\end{eqnarray*}
So,  using  It\^o's formula to $\|\xi_t^\iota\|^{2n}$
and following similar arguments as that in (\ref{2-4}),   one arrives at that
\begin{eqnarray*}
 &&  \dif \|\xi_t^\iota\|^{2n}\leq  2n \|\xi_t^\iota\|^{2n-2}\langle \xi_t^\iota,\dif \xi_t^\iota \rangle+C_n \|\xi_t^\iota\|^{2n-4}
 \|\xi_t\|^4\dif t+C_n \|\xi_t^\iota\|^{2n-2}
 \|\xi_t\|^2\dif t
 \\ &&= 2n \|\xi_t^\iota\|^{2n-2}\langle \xi_t^\iota,\nu \Delta \xi_t+
 B(\cK \xi_t, w_t)+B(\cK w_t,\xi_t )\rangle \dif t
 \\ && \quad\quad +C_n \|\xi_t^\iota\|^{2n-4}
 \|\xi_t\|^4\dif t+C_n \|\xi_t^\iota\|^{2n-2}
 \|\xi_t\|^2\dif t+\dif M_t
 \\ &&\leq C_n \|w_t\|_1^{4/3} \|\xi_t^\iota\|^{2n}\dif t+C_n  \|\xi_t^\iota\|^{2n-2}\|w_t\|^2\|\xi_t^h \|_{1/2}^2\dif t
 \\ && \quad\quad +C_n \|\xi_t^\iota\|^{2n}\dif t+C_n \|\xi_t^h\|^{2n}\dif t+\dif \mathscr{M}_t,
\end{eqnarray*}
where $\mathscr{M}_t= 2n \sum_{j\in \cZ_0 }\int_0^t \|\xi_s^\iota\|^{2n-2}\big(Dq_j(w_s)\xi_s\big)\langle \xi_s^\iota, e_j \rangle \dif W_j(s).$
Therefore, we conclude that
\begin{eqnarray*}
  && \mE
  \big[\|\xi_t^\iota\|^{2n}
  \exp\{-C_n t-C_n\int_0^t\|w_r\|_1^{4/3}\dif r \}\big]
  \\ &&\leq \mE \int_0^t \|\xi_s^h\|^{2n}  \dif s
  +C_n \mE \int_0^t  \|\xi_s^\iota\|^{2n-2}\|w_s\|^2\|\xi_s^h \|_{1/2}^2\dif s
  \\ && \leq \mE \int_0^t \|\xi_s^h\|^{2n}  \dif s
  +C_n \mE \sup_{s\in [0,t]}\|w_s\|^2  \int_0^t \|\xi_s\|^{2n-2}   \|\xi_s^h \|_{1} \|\xi_s^h \|  \dif s
   \\ && \leq \mE \int_0^t \|\xi_s^h\|^{2n}  \dif s
  \\ && \quad +C_n t^{1/6 }  \mE \sup_{s\in [0,t]}\|w_s\|^2  \big(\int_0^t \|\xi_s\|^{12n-12}\dif s \big)^{1/6}  \big(\int_0^t \|\xi_s^h \|_{1}^2 \dif s\big)^{1/2}\big(\int_0^t \|\xi_s\|^{6}\dif s \big)^{1/6}
  \\ &&\leq \mE \int_0^t \|\xi_s^h\|^{2n}  \dif s
  +C_n t^{1/6}\big(\mE \sup_{s\in [0,t]}\|w_s\|^{12}\big)^{1/6}
  \big( \mE  \int_0^t \|\xi_s\|^{12 n-12}\dif s  \big)^{1/6}
  \\ &&\quad\quad\quad\quad\times
  \big(\mE \int_0^t \|\xi_s^h \|_{1}^2 \dif s  \big)^{1/2}
  \big(\mE \int_0^t \|\xi_s^h \|^6 \dif s  \big)^{1/6} .
\end{eqnarray*}
Setting $t=1$ in the above, with the help of Lemmas \ref{L:5.1},  \ref{2-3},
\ref{1340-3},
 we
arrive at
\begin{eqnarray*}
\mE \big[\|\xi_1^\iota\|^{2n}
  \exp\{-C_n -C_n\int_0^1\|w_r\|_1^{4/3}\dif r \}\big]
  \leq \frac{C_{\kappa,n}}{N^{1/12} } \exp\{\kappa \|w_0\|^2\}.
\end{eqnarray*}
 (\ref{30-3})    follows from the above inequality, (\ref{p1204-1}), Lemma \ref{L:5.1}  and
Cauchy-Schwartz inequality.
The proof is complete.

\end{proof}

\begin{lemma}
\label{19-4}
  For any $w_0\in H, 0<a<b<\infty$, smooth function $\xi$ on $\mT^2, s\in [a,b]$ and $n\in \mN,$ we have
  \begin{eqnarray*}
  w\in C([a,b],H^n) \text{ and }
   \mP\Big(\sup_{s<t\leq b} \| J_{s,t}\xi\|_n <\infty\Big)=1,
  \end{eqnarray*}
  where $w=(w_r)_{r\geq 0}$ is the solution to (\ref{0.1}) with initial value $w_0$ and $J_{s,t}\xi, t\geq s$ is the solution to (\ref{10-1}).
\end{lemma}
\begin{proof}
Noticing that $\xi,$ $\sum_{j\in \cZ_0}\int_0^t q_j(w_s)W_j(s)e_j$ and
$\int_s^t \sum_{j\in \cZ_0}  \big(Dq_j(w_r)J_{s,r}\xi\big) e_j \dif W_j(r)$
 are  spatially smooth,
by (\ref{10-1}) and  following the classical proof  for the stochastic 2D Navier--Stokes equations  in  \cite{KS-book}, we arrive at the desired result.

\end{proof}

\begin{lemma}
\label{18-2}
  For any  $T>0$  and   $\lambda\in (0,  T]$,    denote
  \begin{eqnarray*}
D_{ \lambda,T}=\{(t,s)\in \mR^2:\lambda  \leq t\leq s\leq T \}.
\end{eqnarray*}
For any $w_0\in H,\lambda\in (0,T]$ and $k,j\in \mZ^2\setminus \{0,0\}$,  define the following $H$-valued processes  on  $D_{\lambda,T}:$
 \begin{eqnarray*}
 && f_1(u)= J_{s,T}\big(\Delta J_{\lambda,t}e_k\big),
 \\ && f_2(u)=J_{s,T}\tilde B(w_t,J_{\lambda,t}e_k),
 \\ && f_3(u)=J_{s,T}\big(\Delta J_{\lambda,t}\tilde B(w_\lambda,e_k)\big),
 \\ && f_4(u)=J_{s,T}\tilde B(w_t,J_{\lambda,t}\tilde B(w_\lambda,e_k)),
 \\ && f_5(u)=J_{s,T} \tilde B\big(\nu \Delta w_t+B(\cK w_t,w_t),e_k\big),
 \\ &&  f_6(u)=J_{t,s}e_k,\quad
  f_7(u)=J_{t,s} \tilde B(e_k,e_j),
\\ &&  f_8 (u)=J_{t,s} \tilde B(w_t,e_k) ~~\forall u=(t,s)\in D_{ \lambda,T},
 \end{eqnarray*}
 where $w_t$ is the solution to (\ref{0.1}) with initial value $w_0.$
 Then, for  any $T>0,$ there exists  a constant $\vartheta=\vartheta(T)\in [T/2,T)$
 only depending on $T$ and $\nu,d,\aleph $
 such that for any   $w_0\in H,  \lambda\in [\vartheta,T], k,j\in \mZ^2\setminus \{0,0\}$ and $ 1\leq i\leq 8,$ there exists a continuous modification $\tilde f_i$ of $f_i$ in  $H$, i.e., for any $u_0\in D_{ \lambda,T} ,$ it holds that
 \begin{eqnarray*}
   \lim_{u\rightarrow u_0,u\in D_{ \lambda,T} }\| \tilde f_i(u)-\tilde f_i(u_0)\|=0.
 \end{eqnarray*}
 Here,  we say $\tilde f_i$ is a modification of $f_i$
 if for each $u\in D_{ \lambda,T }, $
$
   \mP\big(f_i(u)=\tilde f_i(u)\big)=1.
$

\end{lemma}

\begin{remark}
We believe that   the results  of   Lemma  \ref{18-2} hold for any    $\vartheta\in (0,T]$,
but the current weaker version is enough for our proof.
We prove  this lemma  through the verifications  of  the conditions of Kolmogorov's continuity theorem.
The processes  are very
tedious and don't  need  much technique, so
we    provide  their    details   in Appendix \ref{A}.
\end{remark}

The following lemma gathers some estimates from Section~4.8 in~\cite{HM-2006}  and Lemma~A.6 in~\cite{FGRT-2015}.
\begin{lemma}
\label{L:2.2}For any
  $0\le s<t$ and   $\beta>0$, it holds that
  \begin{gather}
\|\cA_{s,t}^*(\cM_{s,t}+\beta\I)^{-1/2}\|_{\cL({H},L^2([s,t];\R^d))} \le  1,\label{2.8}
   \\
   \|(\cM_{s,t}+\beta\I)^{-1/2}\cA_{s,t}\|_{\cL(L^2([s,t];\R^d),{H})} \le  1,\label{2.9}
   \\  \|(\cM_{s,t}+\beta\I)^{-1/2}\|_{\cL({H},{H})} \le  \beta^{-1/2}.\label{2.10}
  \end{gather}
\end{lemma}

As in  ~\cite{HM-2006}, if $A:\cH_1\rightarrow \cH_2$ is a random linear map between two Hilbert spaces, we denote by
$\cD_s^iA:\cH_1\rightarrow \cH_2 $ the random linear map defined by
\begin{eqnarray*}
  (\cD_s^iA)h=\langle \cD_s(Ah),\theta_i\rangle.
\end{eqnarray*}

\begin{lemma}\label{L:2.3}
The   operators $J_{s,t}$, $\cA_{s,t}$, and $\cA_{s,t}^*$
are Malliavin differentiable, and for any
$\kappa>0$, $m, n\in \mN, r\in[2n,2n+1] $ and $2n\leq s\leq t\leq 2n+1$, the following inequalities~hold
\begin{align}
\label{2-7}
  \mE \big[ \|\cA_{s,t}\|^{m}_{\cL(L^2([s,t];\R^d),{H})}\big|\cF_{2n} \big] & \leq C_{\kappa,m}\exp\{\kappa \|w_{2n}\|^2\},
  \\
  \label{2.11}  \sup_{\xi\in H:\|\xi\|\leq 1}\E  \[ \|\cD_r^iJ_{s,t}\xi\|^{m}\big|\cF_{2n} \]  & \le    C_{\kappa,m}  \exp\{\kappa \|w_{2n}\|^2\},
  \\  \label {2.12}
   \E \[  \|\cD_r^i \cA_{2n,2n+1}\|^{m}_{\cL(L^2([2n,2n+1];\R^d),{H})}\big|\cF_{2n}\]   & \le    C_{\kappa,m} \exp\{\kappa \|w_{2n}\|^2\},
  \\   \label{2.13}  \E  \[   \|\cD_r^i\cA^*_{2n,2n+1}\|^{m}_{\cL({H},L^2([2n,2n+1];\R^d))} \big|\cF_{2n}\] & \le    C_{\kappa,m}  \exp\{\kappa \|w_{2n}\|^2\},
\end{align}
where $C_{\kappa,m}>0$ is a constant only depending on $\kappa,m$ and $\nu,\aleph,d$.

\end{lemma}

\begin{proof}

We only prove this lemma  for the case $n=0,$ and for the general case, the proof is similar.

First, we give a proof of (\ref{2-7}).
For any $v\in L^2([s,t];\R^d)$,
  \begin{eqnarray}
  \nonumber  && \|\cA_{s,t}v\|= \big\|\int_{s}^tJ_{r,t}Q(w_r)v(r)\dif r\big \|
  =  \big \| \int_{s}^tJ_{r,t}\big(\sum_{j\in \cZ_0}q_j(w_r)v_j(r)e_j\big)\dif r \big \|
    \\  \nonumber && =  \big \| \sum_{j\in \cZ_0}\int_{s}^t q_j(w_r)v_j(r)  J_{r,t}e_j  \dif r \big \| \leq  \sum_{j\in \cZ_0} \int_{s}^t |q_j(w_r)|\cdot |v_j(r)| \cdot   \| J_{r,t}e_j\| \dif r
    \\  \label{24-5}  && \leq  C  \sum_{j\in \cZ_0} \big(\int_s^t  |v_j(r)|^2\dif r\big)^{1/2} \cdot  \sum_{j\in \cZ_0}  \big( \int_s^t \| J_{r,t}e_j\|^2\dif r\big)^{1/2}.
  \end{eqnarray}
  By the above inequality and  Lemma \ref{2-3},   the  proof of (\ref{2-7}) is complete. 

 Now, let us   prove (\ref{2.11})   for the case $n=0.$
In the followings,  we divide into  the following   two cases to prove (\ref{2.11}).
Assume that $\xi\in H$ and $\|\xi\|\leq 1.$

Case 1: $r\leq s\leq t.$ In this case,
 $\zeta_t=\cD_r^iJ_{s,t}\xi$
is the solution to the following equation:
\begin{align*}
	& \dif  \zeta_t - \nu  \Delta \zeta_t \dif t
-\tilde{B}(J_{r,t} Q\theta_i,J_{s,t}\xi) \dif t- \tilde{B}( w_t,\zeta_t) \dif t
\\ &=
\sum_{j\in \cZ_0}  \big(Dq_j(w_t)\zeta_t\big)   e_j \dif W_j(t)+\sum_{j\in \cZ_0}  D^2q_j(w_t)(J_{r,t} Q\theta_i, J_{s,t}\xi)  e_j \dif W_j(t),
  \\  &  \zeta_t|_{t=s}=0.
\end{align*}
Using  It\^o's formula to  $\|\zeta_t\|^{2m},$ it holds that
\begin{eqnarray*}
&& \dif \|\zeta_t \|^{2m}=2m\|\zeta_t \|^{2m-2}  \langle
\zeta_t, \dif \zeta_t  \rangle
\\ &&~+m \|\zeta_t \|^{2m-2}\sum_{j\in \cZ_0} \Big( Dq_j(w_t)\zeta_t+ D^2q_j(w_t)(J_{r,t} Q\theta_i, J_{s,t}\xi)\Big)^2 \dif t
\\ &&~+2m(m-1)\|\zeta_t \|^{2m-4} \sum_{j\in \cZ_0}\langle
\zeta_t,e_j\rangle^2  \Big( Dq_j(w_t)\zeta_t+ D^2q_j(w_t)(J_{r,t} Q\theta_i, J_{s,t}\xi)\Big)^2\dif t.
\end{eqnarray*}
  In view of  the above equality,  Lemma \ref{19-1}, Young's inequality
 and
   \begin{eqnarray*}
  &&  \big| \langle \zeta_t, \tilde{B}(J_{r,t} Q\theta_i,J_{s,t}\xi) \rangle
   \big|
   +  \big| \langle \zeta_t, \tilde{B}( w_t,\zeta_t)  \rangle \big|
  \\ &&\leq  C \| \zeta_t\|_1 \| J_{s,t}\xi\| \|J_{r,t} Q\theta_i\|_{1/2}
  +C\|\zeta_t\|_{1/2}\|\zeta_t\|\|w_t\|_1
  \\ &&\leq  \frac{\nu}{4} \| \zeta_t\|_1^2+C\| J_{s,t}\xi\|^2  \|J_{r,t} Q\theta_i\|_{1/2}^2+C\| \zeta_t\|^2\|w_t\|_1^{4/3},
\end{eqnarray*}
we obtain
\begin{eqnarray*}
  && \mE \Big[ \|\zeta_t\|^{2m}e^{-C_m \int_0^t \|w_s\|_1^{4/3}\dif s-C_mt  }
  \Big] \leq C_m \mE  \int_s^t\big( \|J_{r,u} Q\theta_i\|^{4m}_{1/2}+\| J_{s,u}\xi\|^{4m}\big)  \dif u
  \\ &&\leq C_{\kappa,m}\exp\{\kappa \|w_0\|^2\}.
\end{eqnarray*}
With the help of Lemma \ref{L:5.1} and Cauchy-Schwartz inequality, the above   implies
\begin{eqnarray*}
 \mE \|\zeta_t\|^{m}\leq  C_{\kappa,m}\exp\{\kappa \|w_0\|^2\}.
\end{eqnarray*}

Case 2: $s\leq r < t.$ In this case
$\zeta_t=\cD_r^iJ_{s,t}\xi$
satisfies   the following equation:
\begin{align*}
	&  \zeta_t =\int_r^t \nu \Delta \zeta_u   \dif u
+\int_r^t \tilde{B}(J_{r,u}Q \theta_i,J_{s,u}\xi) \dif u+ \int_r^t \tilde{B}( w_u,\zeta_u) \dif u
\\ &~+
\sum_{j\in \cZ_0}  \int_r^t\big( Dq_j(w_u )\zeta_u\big)    e_j \dif W_j(u)+\sum_{j\in \cZ_0} \int_r^t  D^2q_j(w_u)(J_{r,u}Q \theta_i, J_{s,u}\xi)  e_j \dif W_j(u)
\\ &~+\big( Dq_i(w_r)J_{s,r}\xi\big)  e_i.
\end{align*}
With similar   arguments as that in the case 1, one arrives at that
\begin{eqnarray*}
  && \mE \Big[ \|\zeta_t\|^{2m}e^{-C_m \int_0^t \|w_s\|_1^{4/3}\dif s-C_mt  }
  \Big]
  \\ && \leq C_m \mE  \int_r^t \[ \|J_{r,u} Q\theta_i\|^{4m}+\| J_{s,u}\xi\|^{4m}\]  \dif u+C_m \mE \|Dq_i(w_r)J_{s,r}\xi e_i\|^{2m}
  \\ &&\leq C_{\kappa,m}\exp\{\kappa \|w_0\|^2\},
\end{eqnarray*}
which also   implies
$
 \mE \|\zeta_t\|^{m}\leq  C_{\kappa,m}\exp\{\kappa \|w_0\|^2\}.
$

 In the end, we give a proof of (\ref{2.12}).
 Once we have proved  (\ref{2.12}), (\ref{2.13}) follows immediately since $\cA_{s,t}^{*}$ is the adjoint of $\cA_{s,t}. $
For any $v\in L^2([0,1];\R^d)$,
  \begin{eqnarray*}
  &&  \big\|  \cD_r^i \cA_{0,1}v\big\| = \big\| \cD_r^i\int_{0}^1J_{s,1}Q(w_s)v(s)\dif s\big\|
=  \big\| \cD_r^i \int_{0}^1J_{s,1}(\sum_{j\in \cZ_0}q_j(w_s)v_j(s)e_j)\dif s\big\|
    \\ && = \Big\|  \sum_{j\in \cZ_0}\int_{0}^1 q_j(w_s)v_j(s) \big(\cD_r^i J_{s,1}e_j\big) \dif s+\sum_{j\in \cZ_0}\int_{0}^1\big( Dq_j(w_s)\cD_r^iw_s \big)     v_j(s) \big( J_{s,1}e_j\big) \dif s\Big\|
    \\ && \leq   C    \sum_{j\in \cZ_0}  \big(\int_0^1  |v_j(s)|^2\dif s\big)^{1/2} \sum_{j\in \cZ_0}\Big( \int_0^1 \|\cD_r^i J_{s,1}e_j\|^2\dif s+\int_r^1 \|J_{r,s}Q\theta_i\|^2\|J_{s,1}e_j\|^2\dif s \Big)^{1/2}
      \\ && \leq   C  \sum_{j\in \cZ_0} \big(\int_0^1  |v_j(s)|^2\dif s\big)^{1/2} \sum_{j\in \cZ_0}\Big( \int_0^1 \|\cD_r^i J_{s,1}e_j\|^2\dif s+\int_r^1\|J_{r,s}e_i\|^2\|J_{s,1}e_j\|^2\dif s \Big)^{1/2}.
  \end{eqnarray*}
  Hence, by (\ref{2.11}) and Lemma  \ref{2-3}, we obtain the desired result (\ref{2.12}).

\end{proof}

\section{The "invertibility"  of the Malliavin matrix.}
\label{p3-1}

Recall that $$
\langle \cM_{0,T}\xi, \xi  \rangle
  =\sum_{j\in \cZ_0}\int_0^T q_j^2(w_r) \langle J_{r,T}e_j,\xi\rangle^2\dif r,
$$ where $\xi\in H$ and  $\{w_r, r\in [0,T] $ is  the solution to  equation (\ref{0.1}) with initial value $w_0$.  The aim of this section is to prove the following two propositions.
\begin{proposition}
\label{3-10}
 For any   $ w_0\in H, T\in  (0,\infty), \alpha\in (0,1)$ and $N\geq 1,$
   one has
   \begin{eqnarray*}
     \mP(\inf_{\xi \in \cS_{\alpha,N}} \langle  \cM_{0,T}\xi,\xi\rangle=0)=0,
   \end{eqnarray*}
    where $\cS_{\alpha,N}=\{\xi\in H:\|P_N\xi\|\geq \alpha, \|\xi\|=1\}.$
\end{proposition}

In our proof, we will prove a little stronger result than Proposition \ref{3-10} for later use, i.e.,
for any   $ w_0\in H, T\in  (0,\infty), \alpha\in (0,1)$ and $N\geq 1,$
 we have
\begin{eqnarray}
\label{28-1}  \mP\Big(\omega: \inf_{\xi \in \cS_{\alpha,N}  }    \sum_{j\in \cZ_0} \int_{T/2}^T q_j^2(w_r)  \langle \xi,J_{r,T} e_j\rangle^2\dif r  =0  \Big)=0.
\end{eqnarray}

Proposition \ref{3-10}  is not sufficient for the proof of Proposition \ref{1-7}, we need a stronger statement.
For  any  $\alpha\in (0,1], w_0\in H,N\in \mN,\mathfrak{R}>0$ and $\eps>0$, let\footnote{Note that $\cM_{0,t}$ is the Malliavin matrix of $w_t$, the solution of equation (\ref{0.1})  at time $t$ with  initial value $w_0.$ }
\begin{eqnarray}
\label{1-5}
  X^{w_0,\alpha,N}=\inf_{ \xi \in \cS_{\alpha,N}}
  \langle \cM_{0,1}\xi,\xi\rangle.
\end{eqnarray}
and denote
\begin{eqnarray}
\label{35-1}
  r(\eps,\alpha, \mathfrak{R},N)=\sup_{{\|w_0\|\leq  \mathfrak{R}}}\mP(  X^{w_0,\alpha,N}<\eps ).
\end{eqnarray}
Based on Proposition \ref{3-10}  and the dissipative property of Navier-Stokes system,  we obtain
the following  result: 
\begin{proposition}
\label{1-2}
For  any   $\alpha \in (0,1], \mathfrak{R}>0,N\in \mN$,
we have
\begin{eqnarray*}
  \lim_{\eps \rightarrow 0}r(\eps,\alpha,\mathfrak{R},N)=0.
\end{eqnarray*}

\end{proposition}

The proof of this  proposition   is essentially the same as that in \cite[Proposition 3.5]{PZZ24}.
But for  the reader's convenience, we provide details   in    Appendix \ref{B}.

The remainder of this section is devoted to prove   (\ref{28-1}).
First, we demonstrate some notation and  five lemmas, then we give  a proof of Proposition \ref{3-10}.

For  any  $T\in (0,\infty),$ by Lemma \ref{18-2}, there exist $\vartheta=\vartheta(T)  \in [T/2,T)$ and   a set $\dot{\Omega}$ with $\mP(\dot {\Omega})=1 $ such that
\begin{itemize}
  \item[(\romannumeral1)]
  for any $\omega \in \dot{\Omega}$,   $\lambda \in [\vartheta,T]\cap \mQ$
 and $k \in \mZ_*^2$,
the following five      $H$-valued processes
 \begin{eqnarray}
 \label{2626-1}
 \begin{split}
 & f_1(t,s)= J_{s,T}\big(\Delta J_{\lambda,t}e_k\big)
 \\ & f_2(t,s)=J_{s,T}\tilde B(w_t,J_{\lambda,t}e_k),
 \\ & f_3(t,s)=J_{s,T}\big(\Delta J_{\lambda,t}\tilde B(w_\lambda,e_k)\big),
 \\ & f_4(t,s)=J_{s,T}\tilde B(w_t,J_{\lambda,t}\tilde B(w_\lambda,e_k)),
 \\ & f_5(t,s)= J_{s,T}\tilde B\big(\nu \Delta w_t+B( \cK w_t,w_t),e_k\big)
 \end{split}
 \end{eqnarray}
 are continuous  with respect to $(t,s)\in D_{\lambda,T} $ in
 $H$, where
  \begin{eqnarray*}
  D_{\lambda,T}=\{(t,s)\in \mR^2:\lambda \leq t\leq s\leq T \}.
\end{eqnarray*}
  \item[(\romannumeral2)]
    for any $\omega \in \dot{\Omega}$
 and $k,j\in \mZ_*^2$,
 the following three  mappings are also continuous:
\begin{eqnarray}
\label{p1112-3}
\begin{split}
 & s\in [\vartheta,T] \rightarrow J_{s,T} e_k\in H,
  \\ & s\in [\vartheta,T] \rightarrow J_{s,T} \tilde B(e_j,e_k)\in H,
 \\ &
  s\in [\vartheta,T] \rightarrow J_{s,T} \tilde B(w_s,e_k)\in H.
  \end{split}
\end{eqnarray}

\end{itemize}

Before we give a proof of Proposition \ref{3-10}, we demonstrate five lemmas first.

\begin{lemma}
\label{23-1}
 For any $ k\in \mZ^2\setminus \{0,0\}$, there exists a set $\Omega^{1,k}$ with $\mP(\Omega^{1,k})=1$ such that for any $\omega\in \Omega^{1,k}$, one has
  \begin{eqnarray*}
  && \sum_{j\in \cZ_0}\sup_{s\in [\vartheta,T]} \big|\langle J_{s,T}e_j, \xi \rangle\big|(\omega)+ \sup_{s\in [\vartheta,T]} \big|\langle J_{s,T}e_k, \xi \rangle\big|(\omega)=0\text{ and }\xi\in \cS_{\alpha,N}
  \\ && \Longrightarrow
  \sup_{s\in [\vartheta,T]}\big|\langle J_{s,T}(\tilde B(w_s,e_k)),\xi\rangle\big|(\omega)=0.
  \end{eqnarray*}
\end{lemma}
\begin{proof}
  For any $ k\in \mZ^2\setminus \{0,0\}, \lambda,s\in [\vartheta,T)\cap \mQ$  with  $\lambda\leq s$,  there exists a set $\Omega^{1,k,\lambda,s}$ with $\mP(\Omega^{1,k,\lambda,s})=1$ such that
  for any $\omega\in \Omega^{1,k,\lambda,s},$
   \begin{eqnarray}
  \nonumber &&  J_{\lambda,T}e_k-J_{s,T}e_k
=J_{s,T}(J_{\lambda,s}e_k-e_k)
\\    \nonumber && =
 J_{s,T}\Big(\int_{\lambda}^{s} \nu \Delta J_{\lambda,t}e_k \dif t
 +\int_{\lambda}^{s}  \tilde{B}(w_t,J_{\lambda,t}e_k)  \dif t\Big)
\\   \nonumber && \quad\quad + J_{s,T}\Big( \sum_{j\in \cZ_0}  \int_{\lambda}^{s}\big(Dq_j(w_t)J_{\lambda,t}e_k\big) e_j \dif W_j(t) \Big)
 \\   \nonumber && =  \int_{\lambda}^{s} \nu J_{s,T} \big( \Delta J_{\lambda,t}e_k \big) \dif t
 +\int_{\lambda}^{s}  J_{s,T}  \tilde{B}(w_t,J_{\lambda,t}e_k)  \dif t
 \\  \label{26-2}  && \quad\quad +\sum_{j\in \cZ_0}  \int_{\lambda}^{s}\big(Dq_j(w_t)J_{\lambda,t}e_k\big) \dif W_j(t)   J_{s,T}  e_j.
  \end{eqnarray}
Set $$\Omega^{1,k}=\bigcap_{\substack{\lambda,s\in \mQ\cap  [\vartheta,T]\\ \lambda< s}}\Omega^{1,k,\lambda,s}\cap \dot{\Omega}.$$
 For any $\omega\in \Omega^{1,k}$ with
 $$
 \sum_{j\in \cZ_0}\sup_{s\in [\vartheta,T]} \big|\langle J_{s,T}e_j, \xi \rangle\big|(\omega)+ \sup_{s\in [\vartheta,T]} \big|\langle  J_{s,T}e_k, \xi\rangle\big|(\omega)=0
  \text{ for some } \xi\in \cS_{\alpha,N},
  $$
  in the equality (\ref{26-2}),
   first  take  inner product with $\xi$,   divide  both sides by $s-\lambda$, and   then     let   $s\rightarrow \lambda^{+}$ with $s \in \mQ$.   By (\ref{2626-1}) and $\sup_{s\in [\vartheta,T]} \big|\langle \xi, J_{s,T}e_j\rangle\big|(\omega)=0,\forall j\in \cZ_0$,
   we finally arrive at
  \begin{eqnarray}
  \label{p1222-1}
    \langle J_{\lambda,T}(\tilde B(w_{\lambda},e_k)),\xi\rangle=0, \forall \lambda\in [\vartheta,T)\cap \mQ.
  \end{eqnarray}
In view of  (\ref{p1112-3}), the mapping $\lambda\in [\vartheta,T] \rightarrow J_{\lambda,T} \tilde B(w_\lambda,e_k)\in H $  is continuous. Thus, (\ref{p1222-1}) implies the desired result  and the proof is  complete.

\end{proof}

\begin{lemma}
\label{21-1}
Recall that $w_t$ is the solution to (\ref{0.1}) with initial value $w_0.$
  For  $\kappa=\frac{1}{4},$ with probability 1, we have
  \begin{eqnarray*}
    \Upsilon: =\sup_{\vartheta \leq  s < t\leq T}\frac{\|w_t-w_s\|}{|t-s|^{\kappa}}<\infty.
  \end{eqnarray*}
\end{lemma}
\begin{lemma}
\label{21-3}
Define $h(t)=\sqrt{2t \log \log(1/t)},t>0.$
For any $T>0, \varepsilon\in (0,1), j\in \cZ_0$ and $s_1\in [\vartheta,T)$  one has
\begin{eqnarray*}
  \mP\Big(\lim_{n\rightarrow \infty}\frac{\int_{s_1}^{s_1+\varepsilon^n}(q_j(w_t)-q_j(w_{s_1}))\dif W_j(t)}{h(\varepsilon^n)}=0\Big)=1.
  \end{eqnarray*}
\end{lemma}
\begin{proof}
  For any $\delta\in (0,1),$ let
  \begin{eqnarray*}
    \alpha_n=h(\varepsilon^n)(1+\delta)\varepsilon^{-5n/4},
    \beta_n=h(\varepsilon^n)\varepsilon^{n/4}/2.
  \end{eqnarray*}
  Then, by exponential martingale inequality,  we obtain
  \begin{eqnarray*}
    && \mP\left(\sup_{s_2\in [s_1,s_1+1]}\Big(\int_{s_1}^{s_2}(q_j(w_t)-q_j(w_{s_1}))\dif W_j(t)
    \right.
    \\ && \left. \quad\quad \quad\quad\quad\quad\quad\quad -\alpha_n \int_{s_1}^{s_2}|(q_j(w_t)-q_j(w_{s_1}))|^2  \dif t  \Big)\geq \beta_n\right)
    \\ &&\leq \exp\{-\alpha_n\beta_n\}=Kn^{-1-\delta}
  \end{eqnarray*}
  for some constant $K$.
  By  Borel-Cantelli Lemma,
  we find that for almost every $\omega,$ there is a integer $n_0(\omega)$
such that for all $n>n_0(\omega)$ and any $s_2\in [s_1,s_1+1]$,
\begin{eqnarray*}
 \int_{s_1}^{s_2}(q_j(w_t)-q_j(w_{s_1}))\dif W_j(t) -\alpha_n \int_{s_1}^{s_2 }|(q_j(w_t)-q_j(w_{s_1}))|^2 \dif t <\beta_n.
\end{eqnarray*}
Setting   $s_2=s_1+\varepsilon^n$ in the above,   we  obtain
\begin{eqnarray*}
 &&  \int_{s_1}^{s_1+\varepsilon^n }(q_j(w_t)-q_j(w_{s_1}))\dif W_j(t)
  \\ && \leq  \alpha_n \int_{s_1}^{s_1+\varepsilon^n }|(q_j(w_t)-q_j(w_{s_1}))|^2 \dif t+\beta_n
  \\ &&\leq  C \alpha_n \int_{s_1}^{s_1+\varepsilon^n }\|w_t-w_{s_1}\|^2 \dif t+\beta_n
  \\ &&\leq  C  \Upsilon \alpha_n  \int_{s_1}^{s_1+\varepsilon^n }(t-s_1)^{2\kappa} \dif t+\beta_n=\frac{C \Upsilon \alpha_n}{2\kappa+1}\varepsilon^{n(1+2\kappa)}+\beta_n.
\end{eqnarray*}
In the above,  $\kappa=\frac{1}{4}$ and  $ \Upsilon$ is a random variable   given in   Lemma \ref{21-1}.
Thus,
\begin{eqnarray*}
&& \limsup_{n\rightarrow \infty}\frac{\int_{s_1}^{s_1+\varepsilon^n }(q_j(w_t)-q_j(w_{s_1}))\dif W_j(t) }{h(\varepsilon^n)}
\\ &&\leq  \limsup_{n\rightarrow \infty}\Big[
\frac{C \Upsilon \alpha_n}{(2\kappa+1) h(\varepsilon^n)}\varepsilon^{n(1+2\kappa)}+\frac{\beta_n }{h(\varepsilon^n)}\Big]=0.
\end{eqnarray*}
Therefore, it holds that
\begin{eqnarray}
\label{21-2}
  \mP\Big(\limsup_{n\rightarrow \infty}\frac{\int_{s_1}^{s_1+\varepsilon^n }(q_j(w_t)-q_j(w_{s_1}))\dif W_j(t) }{h(\varepsilon^n)}
\leq 0\Big)=1.
\end{eqnarray}

In the  above,  replacing   $(q_j(w_t)-q_j(w_{s_1}))$ by the term
$(-q_j(w_t)+q_j(w_{s_1}))$, and with similar arguments,
we get
\begin{eqnarray*}
&&\mP\Big(\liminf_{n\rightarrow \infty}\frac{\int_{s_1}^{s_1+\varepsilon^n }(q_j(w_t)-q_j(w_{s_1}))\dif W_j(t) }{h(\varepsilon^n)}
\geq  0\Big)
 \\ &&= \mP\Big(\limsup_{n\rightarrow \infty}\frac{\int_{s_1}^{s_1+\varepsilon^n }(-q_j(w_t)+q_j(w_{s_1}))\dif W_j(t) }{h(\varepsilon^n)}
\leq 0\Big)=1.
\end{eqnarray*}
Combining the above equality  with (\ref{21-2}), the proof is complete.
\end{proof}

\begin{lemma}
\label{21-4}
Define $h(t)=\sqrt{2t \log \log(1/t)}, t>0$
and denote $c=\frac{1}{2 \sqrt{d}}.$
For  any $s_1\geq 0$,  there exists a
  set $\tilde  \Omega$ only depending on
  $s_1$ and the $d$ dimensional Brownian $\{W_j(t),t\geq 0\}_{j\in \cZ_0}$ such that  $\mP(\tilde  \Omega)=1$ and  for any $\omega\in \tilde   \Omega$,  $\varepsilon \in (0,\frac{c^2}{(2c+8)^2})\cap \mQ$ and    real  numbers     $a_j,j\in \cZ_0$, we have
   \begin{eqnarray*}
\limsup_{n\rightarrow \infty}\sum_{j\in \cZ_0} \frac{a_j \cdot  \big(W_j(s_1+\varepsilon^n)-W_j(s_1)\big) }{h(\varepsilon^n)}(\omega)
     \geq \frac{c}{2} \sum_{j\in \cZ_0}|a_j|.
   \end{eqnarray*}
\end{lemma}
\begin{proof}
We only prove this lemma for the case $s_1=0$; for the case $s_1>0,$ the proof is similar.

For any  $\varepsilon \in (0,\frac{c^2}{(2c+8)^2}) \cap \mQ, j\in \cZ_0$
 and $n\geq 1, $  define  the events
  \begin{eqnarray*}
    &&  \Gamma_{n,j}^1:=\big\{W_j(\varepsilon^n)-W_j(\varepsilon^{n+1})\geq c(1-\sqrt{\varepsilon})h(\varepsilon^n)\big\},
    \\ && \Gamma_{n,j}^{-1}:=\big\{\big(-W_j(\varepsilon^n)\big)-\big(-W_j(\varepsilon^{n+1})\big)\geq c(1-\sqrt{\varepsilon})h(\varepsilon^n)\big\}.
  \end{eqnarray*}
For any
  $$
 \ell =( \ell_j)_{j\in \cZ_0}\in  \mathbb{L}:=\big\{  \ell=(\ell_j)_{j\in \cZ_0}:     \ell_j\in \{1,-1\}\big\},
  $$
  since  the sequence of events  $(\Gamma_{n,j}^{\ell_j})_{j\in \cZ_0}$  are independent,
  we have
  \begin{eqnarray}
 \nonumber   (\sqrt{2\pi})^d \mP(\text{\scalebox{1.3}{$\cap$}}_{j\in \cZ_0}\Gamma_{n,j}^{\ell_j})
  &=& \left(\int_{a}^\infty e^{-\frac{u^2}{2}}\dif u\right)^d
  \\ \label{p1112-1}  &\geq  & \left(\frac{a}{1+a^2}e^{-\frac{a^2}{2}}\right)^d,
  \end{eqnarray}
  where $a=c(1-\sqrt{\varepsilon})\sqrt{2 \log\log \varepsilon^{-n}}/\sqrt{1-\varepsilon}$.
  For a  proof of the  inequality   (\ref{p1112-1}), one can see the proof  of L\'evy's modulus of continuity  theorem
  in \cite[Theorem (2.7), Section 2, Chapter 1]{RY04}. 
 Observe that   the  right side of (\ref{p1112-1})  is  bigger than
$$
    C_{\varepsilon,d} \cdot  n^{- 2 c^2d (1-\sqrt{\varepsilon})^2 /(1-\eps)}=C_{\varepsilon,d} \cdot  n^{- \frac{1-\sqrt{\eps}}{2(1+\sqrt{\eps})} },~\forall n\geq 1.
$$
  As a result, $\sum_{n=1}^\infty\mP\Big(
   \text{\scalebox{1.2}{$\cap$}}_{j\in \cZ_0}\Gamma_{n,j}^{\ell_j} \Big)=\infty.$
  Thus, for any for any $\ell=(\ell_j)_{j\in \cZ_0}\in \mathbb{L}$ and
  $\varepsilon \in (0,\frac{c^2}{(2c+8)^2}),$ by the Borel-Cantelli lemma,  for almost every $\omega\in \Omega,$  the following events
  \begin{eqnarray*}
 \text{\scalebox{1.4}{$\cap$}}_{j\in \cZ_0}\Gamma_{n,j}^{\ell_j}, \quad n=1,2,3,\cdots
  \end{eqnarray*}
  infinitely often.
  Let
  \begin{eqnarray*}
   && \tilde  \Omega^1:=\Big\{\omega: \big(\text{\scalebox{1.4}{$\cap$}}_{j\in \cZ_0}\Gamma_{n,j}^{\ell_j}\big)_{n=1}^\infty
    \text{ infinitely often for any }  \ell \in \mathbb{L}
    \\ && \quad\quad \quad\quad\quad\quad\quad\quad\quad\quad  \text{ and } \varepsilon \in (0,\frac{c^2}{(2c+8)^2}) \cap \mQ
     \Big\}.
  \end{eqnarray*}
  Then, $\mP(\tilde \Omega^1)=1.$

  On the other hand, by the proof of the   Law of the iterated logarithm
  in \cite[Section 1,Chapter 2]{RY04}, there exist a set $\tilde \Omega^2$ with $\mP(\tilde \Omega^2)=1$ such that  for  every $\omega\in \tilde \Omega^2$ and $\varepsilon \in (0,\frac{c^2}{(2c+8)^2})\cap \mQ,$   the following
 \begin{eqnarray*}
   && -2h(\varepsilon^{n+1})
  \leq  W_{j}(\varepsilon^{n+1})\leq 2h(\varepsilon^{n+1}),  \quad \forall j\in \cZ_0
 \end{eqnarray*}
 holds  from some integer $n_0=n_0(\omega, \varepsilon).$

 With these preparations, we will prove that  for any $\omega\in \tilde \Omega:= \tilde \Omega^1\cap \tilde \Omega^2$, $\varepsilon \in (0,\frac{c^2}{(2c+8)^2})\cap \mQ $ and real numbers  $a_j,j\in \cZ_0,$
 \begin{eqnarray*}
\limsup_{n \rightarrow \infty }\sum_{j\in \cZ_0} \frac{a_j \cdot  W_j(\varepsilon^n) }{h(\varepsilon^n)}
     \geq \frac{c}{2} \sum_{j\in \cZ_0}|a_j|.
 \end{eqnarray*}
First, we define the following sign function
\begin{eqnarray*}
 \text{sign}(x)=
  \left\{
  \begin{split}
   &  1,  \quad x\geq 0
    \\
    & -1, \quad  x<0.
  \end{split}
  \right.
\end{eqnarray*}
Assume that the element   $\omega \in \tilde \Omega $
and
\begin{eqnarray*}
 \text{sign}\big(a_j\big)=\ell_j, \forall j\in \cZ_0.
\end{eqnarray*}
 Then, the  events
$
\cap_{j\in \cZ_0}\Gamma_{n,j}^{\ell_j},n=1,2,3,\cdots,
$
  infinitely often, i.e.,  the following events
  \begin{eqnarray*}
 \ell_j \big(W_j(\varepsilon^n)-W_j(\varepsilon^{n+1})\big) \geq c(1-\sqrt{\varepsilon})h(\varepsilon^n), \quad  \forall j\in \cZ_0
  \end{eqnarray*}
  infinitely often.
  Notice that
    \begin{eqnarray*}
   && -2h(\varepsilon^{n+1})
  \leq  W_{j}(\varepsilon^{n+1})\leq 2h(\varepsilon^{n+1}),~~  \forall j\in \cZ_0
 \end{eqnarray*}
 holds  from some integer $n_0=n_0(\omega).$
 Therefore,  the below sequence of events
 \begin{eqnarray*}
\Big\{\omega:  \ell_j  \cdot  W_j(\varepsilon^n)
&=&  \ell_j  \cdot  \big(W_j(\varepsilon^n)-W_j(\varepsilon^{n+1})\big)+\ell_j W_j(\varepsilon^{n+1})
\\ &\geq &   c(1-\sqrt{\varepsilon})h(\varepsilon^n) -2h(\varepsilon^{n+1}),~ \forall j\in \cZ_0 \Big\}_{n=1}^\infty
 \end{eqnarray*}
 infinitely often.
 Consequently, by  $h(\varepsilon^{n+1})\leq 2\sqrt{\varepsilon}h(\varepsilon^n)$
 from some $n$ on, for any $\omega\in \tilde \Omega$,  one arrives at that

 \begin{eqnarray*}
   && \limsup_{n\rightarrow \infty }\sum_{j\in \cZ_0} \frac{a_j \cdot  W_j(\varepsilon^n) }{h(\varepsilon^n)}
 =\limsup_{n \rightarrow \infty}\sum_{j\in \cZ_0} \frac{|a_j|  \cdot  \text{sign}(a_j) W_j(\varepsilon^n) }{h(\varepsilon^n)}
 \\ && =\limsup_{n\rightarrow \infty }\sum_{j\in \cZ_0} \frac{|a_j|  \cdot  \ell_j  W_j(\varepsilon^n) }{h(\varepsilon^n)}\geq  \liminf_{n\rightarrow \infty} \sum_{j\in \cZ_0 }|a_j| \frac{c(1-\sqrt{\varepsilon})h(\varepsilon^n) -2h(\varepsilon^{n+1}) }{ h(\varepsilon^n)}
 \\ &&
  \\
   &&  \geq \big(c-c\sqrt{\varepsilon}-4\sqrt{\varepsilon}\big)\sum_{j\in \cZ_0 }|a_j|\geq \frac{c}{2}\sum_{j\in \cZ_0 }|a_j| .
 \end{eqnarray*}
The proof is complete.

\end{proof}
\begin{lemma}
\label{23-2}
For any $ k\in \mZ^2\setminus \{0,0\}$, there exists a set $\Omega^{2,k}$ with $\mP(\Omega^{2,k})=1$ such that for any $\omega\in \Omega^{2,k}$, one has
  \begin{eqnarray*}
  &&
  \sum_{j\in \cZ_0}\sup_{s\in [\vartheta,T]} \big|\langle J_{s,T}e_j, \xi \rangle\big|(\omega)=0 ,
   \sup_{s\in [\vartheta,T]}\big|\langle  J_{s,T}e_k, \xi \rangle\big|(\omega)=0,
  \\  && \sup_{s\in [\vartheta,T]}\big|\langle J_{s,T}(\tilde B(w_s,e_k)),\xi\rangle\big|(\omega)=0
   \text{ and }\xi \in \cS_{\alpha,N}
  \\ && \Longrightarrow
   \sup_{s\in [\vartheta,T]}\big|\langle J_{s,T}(\tilde B(e_j,e_k)),\xi\rangle\big|(\omega)=0,~~ \forall  j\in \cZ_0.
  \end{eqnarray*}
\end{lemma}

\begin{proof}
Assume that
  \begin{eqnarray}
  \label{20-2}
  \begin{split}
  &
  \sum_{j\in \cZ_0}\sup_{s\in [\vartheta,T]} \big|\langle J_{s,T}e_j, \xi \rangle\big|(\omega)=0 ,
   \sup_{s\in [\vartheta,T]}\big|\langle  J_{s,T}e_k, \xi \rangle\big|(\omega)=0,
  \\  &\sup_{s\in [\vartheta,T]}\big|\langle J_{s,T}(\tilde B(w_s,e_k)),\xi\rangle\big|(\omega)=0
   \text{ and }\xi \in \cS_{\alpha,N}.
   \end{split}
  \end{eqnarray}
%
Then,   for  any $\lambda,s\in [\vartheta,T)\cap \mQ$  with  $\lambda\leq s$,  it holds that
  \begin{eqnarray}
   \nonumber  0&=& \langle J_{\lambda,T}(\tilde B(w_{\lambda},e_k)),\xi\rangle
    -\langle J_{s,T}(\tilde B(w_{s},e_k)),\xi\rangle
    \\ \nonumber &=&   \big\langle J_{\lambda,T}(\tilde B(w_{\lambda},e_k))- J_{s,T}(\tilde B(w_{\lambda},e_k)),\xi\big\rangle
    \\ \nonumber &&+ \big\langle J_{s,T}(\tilde B(w_{\lambda},e_k))- J_{s,T}(\tilde B(w_{s},e_k)),\xi \big\rangle
    \\ \label{21-5}  &:=&\langle I_1(\lambda,s),\xi\rangle-\langle I_2(\lambda,s),\xi\rangle.
  \end{eqnarray}
  Define $h(t)=\sqrt{2t \log \log(1/t)}, f(t)=\nu \Delta w_t+ B( \cK w_t,w_t)$ for  $t>0$ and  denote $\zeta=\tilde B(w_{\lambda},e_k).$   By direct calculations,  for any $k\in \mZ^2\setminus \{0,0\}$ and $\lambda<s$, there exists a set $\hat \Omega^{k,\lambda,s}$ with $\mP(\hat \Omega^{k,\lambda,s})=1$ such that for any $\omega\in \hat \Omega^{k,\lambda,s},$
it holds that
  \begin{eqnarray*}
    && I_1(\lambda,s)
    \\  & &=
  J_{s,T}\big(J_{\lambda,s}\zeta-\zeta\big)
    \\ & &=   J_{s,T} \Big(\int_{\lambda}^{s} \nu \Delta J_{\lambda,t}\zeta \dif t
 +\int_{\lambda}^{s}  \tilde{B}(w_t,J_{\lambda,t}\zeta )  \dif t
+  \sum_{j\in \cZ_0}  \int_{\lambda}^{s}\big(Dq_j(w_t)J_{\lambda,t}\zeta \big)  \dif W_j(t)  e_j  \Big)
 \\ & &=    \int_{\lambda}^{s} \nu J_{s,T}\big(\Delta J_{\lambda,t}\zeta \big) \dif t
 +\int_{\lambda}^{s}  J_{s,T} \tilde{B}(w_t,J_{\lambda,t}\zeta )  \dif t
\\ && \quad +\sum_{j\in \cZ_0}  \int_{\lambda}^{s}\big(Dq_j(w_t)J_{\lambda,t}\zeta \big)  \dif W_j(t)  J_{s,T} e_j
  \end{eqnarray*}
  and
  \begin{eqnarray}
\nonumber   && I_2(\lambda,s)=
   J_{s,T}\tilde B(w_{s},e_k)-J_{s,T}\tilde B(w_{\lambda},e_k)
   \\ \nonumber && =  J_{s,T} \tilde B\big(\int_{\lambda}^{s}f(t)\dif t,e_k\big)
   +\sum_{j\in \cZ_0}\int_{\lambda}^{s}q_j(w_t)\dif W_j(t)  J_{s,T}\tilde B\big( e_j,e_k\big)
   \\ \nonumber  &&= J_{s,T}  \int_{\lambda}^{s}  \tilde  B\big( f(t),e_k\big)\dif t
   +\sum_{j\in \cZ_0}\int_{\lambda}^{s}\big(q_j(w_t)-q_j(w_{\lambda})\big)\dif W_j(t)   J_{s,T}\tilde B\big( e_j,e_k\big)
   \\ \nonumber &&\quad +\sum_{j\in \cZ_0} q_j(w_{\lambda}) \big(W_j(s)-W_j(\lambda)\big)  J_{s,T}\tilde B\big( e_j,e_k\big)
    \\  \nonumber &&=  \int_{\lambda}^{s} J_{s,T}  \tilde B\big( f(t),e_k\big)\dif t
   +\sum_{j\in \cZ_0}\int_{\lambda}^{s}\big(q_j(w_t)-q_j(w_{\lambda})\big)\dif W_j(t)   J_{s,T}\tilde B\big( e_j,e_k\big)
   \\  \nonumber &&\quad +\sum_{j\in \cZ_0} q_j(w_{\lambda}) \big(W_j(s)-W_j(\lambda)\big)  J_{\lambda,T}\tilde B\big( e_j,e_k\big)
   \\ \nonumber &&\quad +\sum_{j\in \cZ_0} q_j(w_{\lambda}) \big(W_j(s)-W_j(\lambda)\big)  \big(J_{s,T}\tilde B( e_j,e_k)-J_{\lambda,T}\tilde B( e_j,e_k)\big)
   \\ \label{30-4} &&:=I_{21}(\lambda,s)+I_{22}(\lambda,s)+I_{23}(\lambda,s)
   +I_{24}(\lambda,s).
  \end{eqnarray}

  By Lemmas  \ref{21-3},\ref{21-4} and
  L\'evy's modulus of
  continuity theorem(c.f. \cite[Theorem (2.7), Section 2, Chapter1]{RY04}),
  there exists a set $\bar \Omega$ with $\mP(\bar{\Omega})=1$ and for any $\omega\in\bar \Omega,\lambda \in   [\vartheta,T)\cap \mQ,\varepsilon \in (0,\frac{c^2}{(2c+8)^2})\cap \mQ$
  and real numbers $(a_j)_{j\in \cZ_0}$, it holds that
  \begin{eqnarray}
  \label{22-5}
  \begin{split}
 &  \lim_{n\rightarrow \infty}\frac{\int_{\lambda}^{\lambda+\varepsilon^n}(q_j(w_t)-q_j(w_{\lambda}))\dif W_j(t)}{h(\varepsilon^n)}=0,\quad \forall j\in \cZ_0,
  \\ &
  \limsup_{n \rightarrow \infty}\sum_{j\in \cZ_0} \frac{a_j \cdot  (W_j(\lambda+\varepsilon^n)-W_j(\lambda)) }{h(\varepsilon^n)}(\omega)
     \geq \frac{c}{2} \sum_{j\in \cZ_0}|a_j|,
     \\ &
     \limsup_{t\rightarrow \lambda^{+}} \frac{| W_j(t)-W_j(\lambda)| }{h(t-\lambda)}(\omega)
     \leq 1,\quad \forall j\in \cZ_0,
     \end{split}
  \end{eqnarray}
  where $c=\frac{1}{\sqrt{2d}}.$
 Let
 $$
 \Omega^{2,k}=\text{\scalebox{1.5}{$\cap$}}_{\substack{\lambda,s\in \mQ
\cap
   [\vartheta,T) \\ \lambda< s} }\hat \Omega^{k,\lambda,s}
   ~\text{\scalebox{1.5}{$\cap$}}~
    \dot{\Omega}
    ~\text{\scalebox{1.5}{$\cap$}}~
     \bar\Omega.
     $$

 Firstly, by (\ref{2626-1}) and the assumption (\ref{20-2}),  for any  $\omega\in \Omega^{2,k}$ and $\lambda\in [\vartheta,T)\cap \mQ$,  it holds that
  \begin{eqnarray}
  \label{24-4}
    \lim_{s \rightarrow \lambda^{+},s\in \mQ }\frac{\langle I_1(\lambda,s),\xi\rangle (\omega) }{h(s-\lambda)}=0.
  \end{eqnarray}
Secondly,  we consider the term $I_{21}(\lambda,s).$
Obviously, by (\ref{2626-1}), for any $\omega\in \Omega^{2,k}$  and  $\lambda\in [\vartheta,T)\cap \mQ,$
  it also holds that
  \begin{eqnarray}
  \label{24-1}
   \lim_{s \rightarrow \lambda^{+},s\in \mQ } \frac{\langle I_{21}(\lambda,s),\xi\rangle}{h(s-\lambda)}=0.
  \end{eqnarray}
Thirdly, we consider the term $I_{22}(\lambda,s).$
 By (\ref{22-5}) and (\ref{p1112-3}), for any $\omega\in \Omega^{2,k}$,  $ \varepsilon \in (0,\frac{c^2}{(2c+8)^2})\cap \mQ$
 and    $\lambda\in [\vartheta,T)\cap \mQ,$ it holds that
  \begin{eqnarray}
  \label{24-2}
   \lim_{n\rightarrow \infty  } \frac{\langle I_{22}(\lambda,s),\xi\rangle}{h(s-\lambda)}\Big|_{s=\lambda+\varepsilon^n}=0.
  \end{eqnarray}
Fourthly, we consider  the term  $I_{23}(\lambda,s).$
 By (\ref{22-5}), for any $\omega\in \Omega^{2,k}$,
 $\varepsilon \in (0,\frac{c^2}{(2c+8)^2})\cap \mQ$
 and    $\lambda\in [\vartheta,T)\cap \mQ,$   one arrives at that
 \begin{eqnarray}
 \label{24-3}
 \begin{split}
 &  \limsup_{n \rightarrow \infty }\frac{\langle I_{23}(\lambda,s),\xi\rangle }{h(s-\lambda)}\Big|_{s=\lambda+\varepsilon^n}
 \\ & \geq \frac{c}{2}
  \sum_{j\in \cZ_0 } \Big|  q_j(w_{\lambda})  \big\langle J_{\lambda,T}\tilde B\big( e_j,e_k)\big),\xi\big\rangle\Big|.
  \end{split}
 \end{eqnarray}
 In the end, we consider the term $I_{24}(\lambda,s).$
  By (\ref{p1112-3}),  the map  $\lambda \in [\vartheta,T] \rightarrow   J_{\lambda,T}\tilde B( e_j,e_k)\in H$ is continuous for any $\omega\in \Omega^{2,k}$.
Therefore,  for any   $\omega\in \Omega^{2,k},\varepsilon \in (0,\frac{c^2}{(2c+8)^2})\cap \mQ$,$\lambda\in [\vartheta,T)\cap \mQ$  and $\gamma>0,$
 there exists a $N=N(\omega,\varepsilon,\lambda,\gamma)$ such that  for any $n\geq N$
 \begin{eqnarray*}
   \|J_{\lambda+\varepsilon^n,T}\tilde B( e_j,e_k)-J_{\lambda,T}\tilde B( e_j,e_k)\|
   \leq  \gamma.
 \end{eqnarray*}
 Therefore, by (\ref{22-5})
 \begin{eqnarray*}
  && \limsup_{n \rightarrow \infty }\frac{\big| \langle I_{24}(\lambda,s),\xi\rangle\big| }{h(s-\lambda)}\Big|_{s=\lambda+\varepsilon^n}
  \\ && \leq \gamma \aleph  \sum_{j\in \cZ_0} \limsup_{n \rightarrow \infty }
  \frac{ |W_j(s)-W_j(\lambda)|  }{h(\varepsilon^n)}\Big|_{s=\lambda+\varepsilon^n}
    \leq  \gamma d \aleph.
 \end{eqnarray*}
 Since the $\gamma>0$ is arbitrary, we get
 \begin{eqnarray}
 \label{29-1}
 \lim_{n \rightarrow \infty }\frac{\langle I_{24}(\lambda,s),\xi\rangle }{h(s-\lambda)}\Big|_{s=\lambda+\varepsilon^n}= 0.
 \end{eqnarray}

 Combining the above estimates (\ref{24-4})--(\ref{29-1}) with  (\ref{21-5}) and (\ref{30-4}),
 for any $\omega\in \Omega^{2,k},\varepsilon \in (0,\frac{c^2}{(2c+8)^2})\cap \mQ$
 and    $\lambda\in [\vartheta,T)\cap \mQ,$
 we obtain
 \begin{eqnarray*}
 0&=&  \liminf_{s\rightarrow \lambda^{+},s\in \mQ}\frac{ \langle I_{1}(\lambda,s),\xi\rangle -\langle I_{21}(\lambda, s),\xi\rangle -\sum_{i=2}^4 \langle I_{2i}(\lambda,s),\xi\rangle }{h(s-\lambda)}
  \\ & = &  -\limsup_{s\rightarrow \lambda^{+},s\in \mQ}\frac{  \sum_{i=2}^4 \langle I_{2i}(\lambda,s),\xi\rangle }{h(s-\lambda)}.
  \end{eqnarray*}
  The above implies that
  \begin{eqnarray*}
  &&0=\limsup_{s\rightarrow \lambda^{+},s\in \mQ}\frac{  \sum_{i=2}^4 \langle I_{2i}(\lambda,s),\xi\rangle }{h(s-\lambda)}
 \\ && \geq
  \lim_{n \rightarrow \infty }\frac{ \langle  I_{22}(\lambda,s),\xi\rangle }{h(s-\lambda)} \Big|_{s=\lambda+\varepsilon^n } +\limsup_{n \rightarrow \infty }\frac{\langle I_{23}(\lambda,s) ,\xi\rangle }{h(s-\lambda)}\Big|_{s=\lambda+\varepsilon^n }
  \\ && \quad\quad +\lim_{n \rightarrow \infty }\frac{\langle I_{24}(\lambda,s) ,\xi\rangle }{h(s-\lambda)}\Big|_{s=\lambda+\varepsilon^n }
  \\ && \geq  \frac{c}{2}
  \sum_{j\in \cZ_0 } \Big|  q_j(w_{\lambda})  \big\langle J_{\lambda,T}\tilde B\big( e_j,e_k)\big),\xi\big\rangle\Big|.
 \end{eqnarray*}
 By (\ref{p1112-3}), the map $\lambda \in [\vartheta,T]\rightarrow J_{\lambda,T}\tilde B( e_j,e_k) \in H $ is continuous. So
 the above  inequality  imply  the desired result. The proof is complete.
\end{proof}

\textbf{We are now in the position to prove (\ref{28-1}).}
Set
\begin{eqnarray}
\nonumber   &&
  \mathcal{L}\! =\! \Big\{\omega: \inf_{\xi \in \cS_{\alpha,N}  } \! \sum_{j\in \cZ_0}\! \int_{\vartheta}^T q_j^2(w_s)  \langle \xi,J_{s,T} e_j\rangle^2\dif s  =0  \Big\}
 \\  \nonumber     && \quad\quad\quad\quad
  \text{\scalebox{1.44}{$\cap$} }
   \dot \Omega  ~\text{\scalebox{1.44}{$\cap$}}~
  \cap_{k\in \mZ_*^2} \Omega^{1,k}
 ~\text{\scalebox{1.44}{$\cap$}}
  \cap_{k\in \mZ_*^2} \Omega^{2,k},
\end{eqnarray}
where $\vartheta=\vartheta(T)\in [\frac{T}{2},T)>0$ is given by (\ref{2626-1}),
$\dot \Omega$ is given in (\ref{2626-1})  and (\ref{p1112-3}),
$\Omega^{1,k}$ and  $\Omega^{2,k}$ are given by Lemma \ref{23-1} and Lemma  \ref{23-2}, respectively.
In the following, we will prove $\mathcal{L}=\varnothing$, completing the proof of  (\ref{28-1}) and  Proposition \ref{3-10} .

Let  $\omega\in \cL,$
then, for some $\xi \in \cS_{\alpha,N}$,    it holds that
\begin{eqnarray}
\label{p1204-6}
  \langle \xi, J_{s,T}e_j\rangle=0, \quad \forall s\in [\vartheta,T],j\in \cZ_0.
\end{eqnarray}
Checking through the    arguments in Lemma \ref{23-1} and Lemma \ref{23-2},  for any $ k\in \mZ_*^2$ and $\omega\in \cL$,  we actually have
\begin{eqnarray}
\label{1657}
\begin{split}
 & \langle \xi, J_{s,T} e_k \rangle = 0, ~\forall s\in [\vartheta,T],
  (\ref{p1204-6}) \text{ and }     \xi \in \cS_{\alpha,N}
   \\ &  \Longrightarrow
  \langle \xi,J_{s,T}\tilde B(e_k,e_j)\rangle= 0, ~~\forall j\in \cZ_0\text{ and }s\in [\vartheta,T].
  \end{split}
\end{eqnarray}
Define the set $\cZ_n\subseteq \mZ_*^2$ recursively:
\begin{eqnarray*}
  \cZ_n=\{k+j \big| j\in \cZ_0,k \in \cZ_{n-1} \text{ with } \langle k^{\perp},j\rangle\neq 0, |k|\neq |j| \},
\end{eqnarray*}
where $i^{\perp}=(i_2,-i_1).$
Assume that we have proved
\begin{eqnarray*}
 \langle \xi, J_{s,T} e_j \rangle = 0, ~~\forall j\in \cZ_{n-1} \text{ and }s\in [\vartheta,T].
\end{eqnarray*}
Then, by (\ref{1657}), it follows that
\begin{eqnarray*}
 \langle \xi, J_{s,T}\tilde B(e_j,e_i) \rangle= 0,~~ \forall j\in \cZ_{n-1}, i\in \cZ_0 \text{ and }s\in [\vartheta,T].
\end{eqnarray*}
It is easy to verify that $\cZ_m$ is symmetric for any $m\geq 0,$ i.e.
$\cZ_m=-\cZ_m.$
Also by the definition of $\cZ_n,$  one can see  that
\begin{eqnarray*}
\{ e_j, j\in \cZ_n \} \subseteq  \text{span}  \{ \tilde B(e_i,e_j):  j\in \cZ_0, i\in \cZ_{n-1}\}.
\end{eqnarray*}
Hence,
\begin{eqnarray*}
 \langle \xi, J_{s,T} e_j \rangle = 0,~~ \forall j\in \cZ_{n} \text{ and }s\in [\vartheta,T].
\end{eqnarray*}
By this  recursion,
\begin{eqnarray*}
 \langle \xi, J_{s,T} e_j \rangle = 0,~~ \forall j\in \cup_{n=1}^\infty \cZ_{n}=\mZ_*^2 \text{ and }t\in [\vartheta,T].
\end{eqnarray*}
(Here, we have used  \cite[Proposition 4.4]{HM-2006}.)
Let $s=T$ to get  $\xi=0$,  which contradicts $\xi\in \cS_{\alpha,N}$ .
Therefore, $\mathcal{L}=\varnothing$.

The proof of Proposition \ref{3-10}  is complete.

\section{Proof of Proposition   \ref{1-7}}
\label{p4-1}

Let us take $f\in C_b^1(H)$ and $\xi\in H$ with $\|\xi\|=1.$  Compute the derivative of $ \mE_{w_0} f(w_{t})$ with respect to $w_0$ in the direction $\xi$:
\begin{eqnarray}
\label{3.2}
  && D_\xi  \mE_{w_0} f(w_{t})= \mE D
  f(w_{t})J_{0,t }\xi.
\end{eqnarray}
Similar to   the papers~\cite{HM-2006, HM-2011},
the main   ideas of the  proof of  the  asymptotic strong Feller property is~to approximate the perturbation
  $J_{0,t}\xi$ caused by the variation of the initial condition with a variation, $\cA_{0,t}v=\cD^v w_t $, of the noise by an appropriate   process $v$. Denote by $\rho_t$ the  residual error between~$J_{0,t}\xi$    and~$\cA_{0,t}v$:
\begin{align*}
  \rho_t=J_{0,t}\xi-\cA_{0,t}v.
\end{align*}

We  define the  perturbation $v$ to be $0$
 on all intervals of the type $[{n+1},{n+2}],
 n\in 2\mN,$ and by some  $v_{{n},{n+1}} \in L^2([n,n+1],H), n\in 2\mN$
on the remaining intervals.
For fixed  $n\in 2\mN,$ define
the infinitesimal variation:
\begin{eqnarray}
\label{p-1}
\begin{split}
  & v_{{n},{n+1}}(r)=\cA^*_{n,n+1}
  \big(\cM_{n,n+1}+\beta\mI \big)^{-1}J_{n,n+1}\rho_{n},
  ~ r\in [{n},{n+1}]
 \\
 & v_{{n+1},{n+2}  }(r)=0,~ r\in [{n+1},{n+2}].
  \end{split}
\end{eqnarray}
where $\rho_{n} $ is the residual of the infinitesimal
 displacement at time $n$,  i.e., $\rho_n=J_{0,n}\xi-\cA_{0,n}v.$  Obviously,    $\rho_0=J_{0,0}\xi-\cA_{0,0}v=\xi.$
 The  $\beta$  in (\ref{p-1}) is a positive  constant will be decided later.
Here and after, we use  $v_{a,b}$  to denote the function $v$ restricted on the interval $[a,b]$.

\subsection{The  choice of  $\beta$  and estimates  on  $\rho_n.$}

Similar to \cite{HM-2006},
for any $\beta\in(0,1),$
  if we define the direction $v$  according to  (\ref{p-1}),
    then
\begin{eqnarray}
\label{26-3}
\rho_{n+2 }&=&{J}_{n+1,n+2}\big[
\beta
\big(  \cM_{n,n+1}+\beta \mI\big)^{-1}\big] {J}_{n,n+1}
\rho_{n}, \quad \forall n\in 2\mN.
\end{eqnarray}

Let  $A_\eps=A_{\eps,w_0,\alpha,N}:=\{X^{w_0,\alpha,N}\geq \eps  \},$
where  the random variable $X^{w_0,\alpha,N}$ is defined in (\ref{1-5}).
To provide an estimate for
$\|\rho_{{n}}\|,$ we start with some preparations.
\begin{lemma}(c.f. \cite[Lemma 5.14]{HM-2011})
\label{3-3}
  For any positive constants  $\beta>0,\eps\in (0,1),\alpha\in(0,1], N\in \mN$
  and $\phi\in H,$  the following inequality holds with probability $1$:
  \begin{eqnarray}
  \label{3-1}
 \beta \|P_N\big(\beta\mI+\cM_{0,1}\big)^{-1}\phi\|
\leq   \| \phi\|
\left(\alpha \vee \sqrt{\beta/\eps }\right) I_{A_\eps} +\|\phi\|I_{A_\eps^c }.
\end{eqnarray}
\end{lemma}
\begin{proof}
On the event  $A_\eps^c$, the inequality (\ref{3-1})  obviously  holds.
On the event $A_\eps$, this inequality is proved in \cite[Lemma 5.14]{HM-2011}, so we omit the details.

\end{proof}

Let $\cR_{n,{n+1}}^{\beta}=\beta(\cM_{n,n+1}+\beta\mI)^{-1}.$
We have  the following estimate for $\cR_{n,{n+1}}^\beta.$
\begin{lemma}
\label{3-9}
    For any $\kappa>0,\delta\in (0,1],  p\geq 1$ and $N\in \mN$, there exists a $\beta>0$ which depends on $\kappa,\delta,p,N$ and $\nu,\aleph,d$
  such  that
    \begin{eqnarray}
    \label{3-4}
\mE\[\|P_N \cR_{{n},{n+1}}^{\beta}\|^p \big| \cF_{n} \]
\leq \delta e^{\kappa \|w_{{n}}\|^2},\quad \forall n\in \mN.
\end{eqnarray}
\end{lemma}
\begin{proof}
   We here give a proof for the case $n=0$ and $p\geq 2.$ The other cases can be proved similarly.
   Let  $\mathfrak{R}=\mathfrak{R}_{\delta,\kappa}$ be a positive  constant  such that
 $
  \exp\{\kappa \mathfrak{R}^2  \}\geq \frac{1}{\delta}.
$
  We divide into the following two cases to prove (\ref{3-4}).

 \textbf{Case 1:} $\|w_0\|\geq \mathfrak{R}$.
  In this case,
  \begin{eqnarray*}
\mE\[\|P_N \cR_{0,1}^{\beta} \|^p \big| \cF_0 \]\leq
1\leq \delta e^{\kappa \|w_{0}\|^2}.
  \end{eqnarray*}

 \textbf{Case 2:}
 $\|w_0\| \leq \mathfrak{R}$.
 For any  positive constants $\eps,\beta$ and $\alpha\in (0,1],$ by Lemma \ref{3-3}, we have
  \begin{eqnarray*}
 \mE\[\|P_N \cR_{0,1}^\beta \|^p \big| \cF_0 \] & \leq&
 C_p \big(\alpha \vee \sqrt{\frac{\beta}{ \eps }}\big)^p+
 C_p r(\eps,\alpha, \mathfrak{R},N),
  \end{eqnarray*}
  where $C_p$ is a constant only depending on $p,$
  and  $r(\eps,\alpha, \mathfrak{R},N)$ is defined in (\ref{35-1}).
 Choose now  $\alpha=\alpha({p})$ small enough such that
 \begin{eqnarray*}
 C_p \alpha^p \leq \frac{\delta}{2}.
 \end{eqnarray*}
 By Proposition \ref{1-2}, $\lim_{\eps \rightarrow 0}  r(\eps,\alpha, \mathfrak{R},N)=0$.
Pick a small constant  $\eps$ such that
 \begin{eqnarray*}
C_p r(\eps,\alpha, \mathfrak{R},N)\leq \frac{\delta}{2}.
 \end{eqnarray*}
 Finally, we choose  $\beta$ small enough so that
 \begin{align*}
  C_p( \sqrt{\beta/\eps})^p <\frac{\delta}{2}.
 \end{align*}
Putting the above steps together, we see that
$
  \mE\[\|P_N \cR_{0,1}^{\beta} \|^p \big|\cF_0  \]  \leq \delta
e^{\kappa \|w_{0}\|^2} .
$
\end{proof}

With regard to   $\rho_n$, we have the following estimate.
\begin{lemma}
\label{30-6}
   For  any   $p\geq 1, \kappa>0$   and  $\delta\in (0,1)$,
 there exist    a    constant   $\beta>0 $  which  depends    on   $p,\kappa,\delta$ and $\nu,\aleph,d$
 such that
  if we define the direction $v$  according to  (\ref{p-1}),
    then
  the following holds
  \begin{eqnarray*}
    \mE \[\|\rho_{{n+2}}\|^{p}   \big| \cF_{n}\]\leq \delta  \exp\{\kappa \|w_{n}\|^2 \}\|\rho_{{n}}\|^{p}, \quad \forall n\in 2\mN.
  \end{eqnarray*}
\end{lemma}
\begin{proof}

By (\ref{26-3}),
for any   $n\in 2\mN, \beta>0$ and $N\in \mN,$   one easily sees that
\begin{eqnarray}
 \nonumber  && \rho_{{n+2}}
 = J_{{n+1},{n+2}}
 \cR_{n,{n+1}}^\beta J_{n,{n+1}}\rho_{n}
  \\  \nonumber && = J_{{n+1},{n+2}}Q_N\cR_{n,{n+1}}^\beta J_{n,{n+1}}\rho_{n}
  +J_{{n+1},{n+2}}P_N\cR_{n,{n+1}}^\beta J_{n,{n+1}}\rho_{n}
  \\  \label{4-1} && :=  \rho_{{n+2}}^{(1)}+ \rho_{{n+2}}^{(2)}.
\end{eqnarray}

  First, we consider the term  $\mE \[\|\rho_{{n+2}}^{(1)}\|^{p}\big| \cF_{n} \].
  $
  By (\ref{30-2}),
  for any $\kappa, p, \delta'>0,$  we can set $N$ big enough such that
  \begin{eqnarray}
  \label{29-2}
  \begin{split}
   & \mE \big[\| J_{{n+1},{n+2}}Q_N\cR_{n,{n+1}}^\beta J_{n,{n+1}}\rho_{n}\|^p \big| \cF_{n+1}\big]
   \\   & \leq \delta' \exp\{\frac{\kappa}{2} \|w_{n+1}\|^2\}\|J_{n,{n+1}}\rho_{n}\|^p.
   \end{split}
  \end{eqnarray}
  So by  Lemma \ref{2-3} and Lemma \ref{L:5.1}, we get
    \begin{eqnarray}
  \nonumber && \mE \big[\| J_{{n+1},{n+2}}Q_N\cR_{n,{n+1}}^\beta J_{n,{n+1}}\rho_{n}\|^p \big| \cF_{n}\big]
  \\ \nonumber && \leq \delta' \mE \big[  \exp\{\frac{\kappa}{2} \|w_{n+1}\|^2\}\|J_{n,{n+1}}\rho_{n}\|^p \big|
  \cF_n \big]
  \\ \nonumber && \leq \delta'\big( \mE \exp\{\kappa  \|w_{n+1}\|^2\} \big|\cF_n \big)^{1/2} \cdot C_{\kappa,p}\exp\{\frac{\kappa}{2}\|w_{n}\|^2 \}\|\rho_{{n}}\|^{p}
 \\ \label{27-1}  && \leq C_{\kappa,p}\delta' \exp\{\kappa \|w_{n}\|^2  \}\|\rho_{{n}}\|^{p}.
  \end{eqnarray}

  Second, we  consider the term  $\mE \[\|\rho_{{n+2}}^{(2)}\|^{p} \big| \cF_{n} \]$.
  Obviously,
  it holds that
  \begin{eqnarray*}
  &&  \mE \Big[ \| J_{{n+1},{n+2}}P_N\cR_{n,{n+1}}^\beta J_{n,{n+1}}\rho_{n}\|^p \big|\cF_{n+1}\Big]
   \\ &&\leq  C_{\kappa,p}  \exp\{\frac{\kappa}{3} \|w_{n+1}\|^2\}\|P_N\cR_{n,{n+1}}^\beta J_{n,{n+1}}\rho_{n}\|^p.
  \end{eqnarray*}
Then,  for any $\kappa, \delta',p>0$ and the $N$ given by  (\ref{29-2}),
 by Lemma \ref{3-9},
  there exist   a     constant   $\beta>0 $
  such that
  \begin{align}
    \nonumber & \mE \big[ \exp\{\frac{\kappa}{3}   \|w_{n+1}\|^2\}\|P_N\cR_{n,{n+1}}^\beta J_{n,{n+1}}\rho_{n}\|^p\big| \cF_n\big]
    \\ \nonumber &\leq  \big(\mE  \exp\{\kappa
      \|w_{n+1}\|^2\}\big| \cF_n\big)^{1/3}
    \big( \mE \|P_N\cR_{n,{n+1}}^\beta\|^{3p}\big|\cF_n\big)^{1/3}
    \big(\mE \|J_{n,{n+1}}\rho_{n}\|^{3p}\big| \cF_n\big)^{1/3}
    \\ \label{27-2}&\leq C_{\kappa,p}\delta' \exp\{\kappa \|w_{n}\|^2 \}\|\rho_{{n}}\|^{p}
  \end{align}

  By (\ref{27-1}) and (\ref{27-2}), for any $\kappa, \delta'>0,$
  there exists    a     constant    $\beta>0 $
  such that
  \begin{eqnarray*}
&& \mE \[\|\rho_{{n+2}}\|^{p}   \big| \cF_{n}\]\leq
C_p \mE \[\|\rho_{{n+2}}^{(1)}\|^{p}   \big| \cF_{n}\]
+C_p \mE \[\|\rho_{{n+2}}^{(2)}\|^{p}   \big| \cF_{n}\]
\\ &&\leq C_{\kappa,p} \delta' \exp\{\kappa \|w_{n}\|^2 \}\|\rho_{{n}}\|^{p}.
  \end{eqnarray*}
 Set $\delta'=\frac{\delta}{C_{p,\kappa}}.$ The proof is complete.

\end{proof}

The constant  $\beta$ in (\ref{p-1}) is decided through  the following Lemma.
\begin{lemma}
\label{L:3.3}
  For any $\gamma>0$ and $\kappa>0$,  there exists  a      constant  $\beta=\beta(\gamma,\kappa, \nu, \aleph,d)>0$
  such that   we define the direction $v$  according to  (\ref{p-1}),
    then
  the following
  \begin{equation}
  \label{3.6}
    \E_{w_0}\|\rho_t\|^8 \le  C_{\gamma,\kappa}  \exp\left\{\kappa \|w_0\|^2-\gamma t \right\}
  \end{equation}
  holds  for  for any $w_0\in  H$ and $t\ge0$,  where
 $C_{\gamma,\kappa}$ is a constant depending   $\gamma,\kappa$ and $\nu,\aleph,d.$
\end{lemma}
\begin{proof}
Following the lines in the proof of \cite[(3.29)]{FGRT-2015}, also with the help of Lemma  \ref{30-6} and Lemma \ref{L:5.1},  for any $\gamma,\kappa>0,$ there is a      $\beta=\beta(\gamma,\kappa, \nu, \aleph,d)>0$   such that
  \begin{equation}
 \label{3.7}
    \E_{w_0}\|\rho_n\|^{32} \le  C_{\gamma,\kappa} \exp\{\kappa \|w_0\|^2- \gamma  n\} \quad \text{for any $w_0\in  H$ and $n\in 2\mN$}.
  \end{equation}From the construction it follows that
 $$
 \rho_t=\begin{cases} J_{n,t}\rho_n-\cA_{n,t}v_{n,t}, & \text{for }t \in [n,n+1], \\ J_{n+1,t}\rho_{n+1}, & \text{for }t\in [n+1,n+2] \end{cases}	
 $$ for any  $n\in 2\mN$.
  Using \eqref{p-1} and inequalities \eqref{2.8} and \eqref{2.10}, we get
    \begin{equation} \label{3.8}
    \|v_{n,n+1}\|_{L^2([n,n+1];\R^d)} \le \beta^{-1/2}\|J_{n,n+1} \rho_n\|.
\end{equation}
Hence,   for any $t\in [n,n+1]$,
   \begin{align}
    \|\rho_t\|&\le  \|J_{n,t}\rho_n\|+\|\aA_{n,t}v_{n,t}\|  \nonumber
    \\  \nonumber  &\le  \|J_{n,t}\rho_n\|+\|\aA_{n,t}\|_{\cL(L^2([n,t];\R^d), H)} \|v_{n,t}\|_{L^2([n,n+1];\R^d)}
    \\  \nonumber  &\le  \|J_{n,t}\rho_n\|+\|\aA_{n,t}\|_{\cL(L^2([n,t];\R^d), H)} \|v_{n,n+1}\|_{L^2([n,n+1];\R^d)}
    \\ \label{3.9}
     &\le C\Big(\|J_{n,t}\rho_n\| +\beta^{-1/2}  \|J_{n,n+1}\rho_n\| \cdot \sum_{j\in \cZ_0}\big( \int_n^t \| J_{r,t}e_j\|^2\dif r\big)^{1/2}  \Big).
  \end{align}
 In the above ,   we have   used  \eqref{24-5} and \eqref{3.8}. For any
     $t\in [n+1,n+2],$
     it holds that
\begin{eqnarray}
\label{24-6}
\|\rho_t\|\le  \|J_{n+1,t}\rho_{n+1}\| .
\end{eqnarray}
  Combining (\ref{3.9})(\ref{24-6})  with  inequalities      \eqref{3.7}, Lemmas \ref{2-3}, \ref{L:5.1}  and   H\"older's  inequality, also noticing that    fact that the    $\kappa>0$ and $\gamma>0$  in (\ref{3.7}) are arbitrary,   we arrive at~\eqref{3.6}.

\end{proof}

\subsection{The control of $\int_0^t v(s)\dif W(s)$}

The aim of this subsection is to prove the following lemma.
\begin{lemma}
\label{L:3.4}
For any $\gamma>0$ and $\kappa>0$,
 let $\beta=\beta(\gamma,\kappa, \nu, \aleph,d)>0$
be a  constant
chosen
according to    Lemma \ref{L:3.3}. Then, the following
  \begin{equation}
  \label{3.10}
  \E_{w_0}  \Big| \int_{2n}^t v(s)\dd W(s)\Big|^2 \le   C_{\gamma,\kappa}\exp\{\kappa \|w_0\|^2-\gamma  n/4 \}
  \end{equation}
    holds  for any $n\geq 0$,   $t\in [2n,2n+1]$ and $w_0\in  H$.
  In the above,   $C_{\gamma,\kappa}$ is a constant depending   $\gamma,\kappa$ and $\nu,\aleph,d.$
\end{lemma}
\begin{proof}  In this proof, we consider the endpoint case $t=2n+1$; the case $t\in [2n,2n+1)$ is treated in a similar way.
  Using  the generalised It$\hat{\text{o}}$ isometry (see Section~1.3 in~\cite{nualart2006}) and the fact that $v(t)=0$ for $t\in [2n+1,2n+2]$ (see \eqref{p-1}), we~obtain
  \begin{align}
  \nonumber &   \E_{w_0}  \left| \int_{2n}^{2n+1} v(s)\dd  W(s)\right|^2
  \\ & =  \E_{w_0}  \int_{2n}^{2n+1} |v(s)|_{\R^d}^2 \dd  s   +  \E_{w_0} \int_{2n}^{2n+1}\!\!\!\int_{2n}^{2n+1} \text{Tr}(\cD_sv(r)\cD_rv(s))\dd  s \dd  r\nonumber
   \\  \nonumber   & \le    \E_{w_0}  \int_{2n}^{2n+1} |v(s)|_{\R^d}^2 \dd  s +  \E_{w_0}  \int_{2n}^{2n+1}\!\!\!\int_{2n}^{2n+1} |\cD_rv_{2n,2n+1}(s)|_{\R^d\times \R^d}^2\dd  s\dd  r\nonumber\\  & =L_1+L_2.
  \label{3.11}
  \end{align}
  We estimate $L_1$ by using
       \eqref{3.8}, \eqref{3.6},  Lemma \ref{2-3} and Lemma \ref{L:5.1}:
     \begin{align}
  \nonumber &  \E_{w_0}  \int_{2n}^{2n+1} |v(s)|_{\R^d}^2 \dd  s
   \le     \beta^{-1}  \E_{w_0}  \|J_{2n,2n+1}\rho_{2n}\|^2
   \\ \nonumber
  &
  \leq  C_\kappa   \beta^{-1}   \E_{w_0} \big[  \exp\{\kappa \|w_{2n}\|^2/2 \}\|\rho_{2n}\|^2\big]
 \\ \nonumber  &  \le
  C_\kappa  \beta^{-1}  \exp\{\kappa \|w_0 \|^2/2 \} \left(  \E_{w_0}  \|\rho_{2n}\|^4\right)^{1/2}
   \\
     \label{3.12}
   &\le   C_{\gamma,\kappa}\exp\{\kappa \|w_0\|^2-\gamma   n/4 \}.
  \end{align}
   To estimate $L_2$,   we use the explicit form of   $\cD_r v.$
  Notice  that,   for any $r\in[2n,2n+1]$ and~$i=1,\dots, d$,
\begin{align*}
   \cD_r^iv_{2n,2n+1}
 &= \cD_r^i(\cA_{2n,2n+1}^*) (\cM_{2n,2n+1}+\beta \I)^{-1}J_{2n,2n+1}\rho_{2n}
    \\ &\quad + \cA_{2n,2n+1}^*(\cM_{2n,2n+1}+\beta\I)^{-1}
 \\ & \quad\quad\times \Big( \cD_r^i(\cA_{2n,2n+1})\cA_{2n,2n+1}^*
  +\cA_{2n,2n+1}\cD_r^i(\cA_{2n,2n+1}^*)\Big) \\ & \quad\quad\times
 (\cM_{2n,2n+1}+\beta\I)^{-1}
 J_{2n,2n+1} \rho_{2n}\\ & \quad +\cA_{2n,2n+1}^* (\cM_{2n,2n+1}+\beta\I)^{-1}\cD_r^i(J_{2n,2n+1}) \rho_{2n}.
\end{align*}
By inequalities \eqref{2.8}-\eqref{2.10}, we have
\begin{align*}
   & \|\cD_r^iv_{2n,2n+1}\|_{L^2([2n,2n+1];\R^d)}
  \\ & \le   \beta^{-1}\| \cD_r^i(\cA_{2n,2n+1})\|_{\cL(L^2([2n,2n+1];\R^d), {H})}
  \|J_{2n,2n+1} \rho_{2n}\|
    \\ &\quad +2\beta^{-1}\| \cD_r^i(\cA_{2n,2n+1}^*)\|_{\cL({H},L^2([2n,2n+1];\R^d))}
  \|J_{2n,2n+1} \rho_{2n}\|  \\ &\quad
  + \beta^{-1/2} \| \cD_r^i(J_{2n,2n+1})  \rho_{2n}\|,
\end{align*}
where
$\|\cD_r^iv_{2n,2n+1}\|_{L^2([2n,2n+1];\R^d)}=\big(\int_{2n}^{2n+1}
|\cD_r^iv_{2n,2n+1}(s)|_{\R^d)}^2\dif s\big)^{1/2}. $
Inequalities  \eqref{2.11}-\eqref{2.13}, \eqref{3.6} and Lemmas \ref{L:5.1},  \ref{2-3}, imply  that
\begin{align}
 \nonumber  & \E_{w_0}  \int_{2n}^{2n+1}\!\!\!\int_{2n}^{2n+1} |\cD_rv_{2n,2n+1}(s)|_{\R^d\times \R^d}^2\dd  s\dd  r
  \\  &  \le  C_\kappa  \beta^{-2}  \exp\{\kappa \|w_0\|^2/2 \} \left(\E \|\rho_{2n}\|^8\right)^{1/4}
 \nonumber \\ \label{3.13} &\le   C_{\gamma,\kappa}  \exp\{\kappa \|w_0\|^2-n\gamma/4 \}.
\end{align}
Combining   estimates \eqref{3.12}-\eqref{3.13} with \eqref{3.11},  the proof is complete.
\end{proof}

\subsection{Proof of Proposition  \ref{1-7}}
For any $\kappa, \gamma>0$, we set the value of  $\beta$ according to   Lemma  \ref{L:3.3}.
Then, by Lemma \ref{L:3.3} and Lemma \ref{L:3.4},
for any $\xi\in H$ with $\|\xi\|=1,$
it holds that
\begin{eqnarray*}
&& \big| D_\xi P_tf(w_0)\big|=\big| \mE D f(w_t) J_{0,t}\xi \big|
= \big| \mE D  f(w_t) \rho_t +  \mE D  f(w_t) \cD^vw_t   \big|
\\ &&\leq   \|D  f\|_\infty \mE \|\rho_t\| +
\big| \mE  f(w_t)\int_0^t v(s)\dif W(s)  \big|
\\ &&\leq  C_{\kappa,\gamma} \exp\{\kappa \|w_0\|^2 \}\|D f\|_\infty e^{-\gamma t/8}
\\ && \quad\quad \quad\quad +\|f\|_\infty  \sum_{n=0}^{\lfloor t \rfloor-1}\mE \big| \int_n^{n+1} v(s)\dif W(s) \big|+\|f\|_\infty \mE \big| \int_{\lfloor t \rfloor}^{t} v(s)\dif W(s) \big|
\\ &&\leq   C_{\kappa,\gamma}  \exp\{\kappa \|w_0\|^2 \}\|D f\|_\infty e^{-\gamma t/8}
+\|f\|_\infty  \sum_{n=0}^\infty \exp\{\kappa \|w_0\|^2-\gamma  n/16   \}
\\ &&\leq  C_{\kappa,\gamma} \exp\{\kappa \|w_0\|^2 \}(\|D f\|_\infty e^{-\gamma t/8}+\|f\|_\infty).
\end{eqnarray*}
 Since $\gamma>0$  is   arbitrary,  we complete  the proof.

\section{Proof of irreducibility}
\label{123}

First, we demonstrate a lemma which plays a key role in the proof of
irreducibility.
\begin{lemma}
\label{16-1}
Under the Condition \ref{16-2},  for any $T,\eps>0$  and $w_0\in H,$   one has
\begin{eqnarray*}
\mP(\sup_{t\in [0,T]}\sum_{i \in \cZ_0 } \big|\int_0^tq_i(w_r)\dif W_i(r)\big|<\eps)\geq p_0>0.
\end{eqnarray*}
where $p_0=p_0({T,\eps,\aleph,d })$ is a positive  constant  only depending on  $T,\eps,\aleph,d.$
\end{lemma}
\begin{proof}
Denote $M_i(t)=\int_0^tq_i(w_r)\dif W_i(r), i\in \cZ_0$.
By \cite[Proposition 16.8]{Kal97}, there exist
independent Brownian motions $B^i=(B^i_t)_{t\geq 0},i\in \cZ_0$ such that
\begin{eqnarray*}
  M_i(t)=B^i_{[M_i](t)},
\end{eqnarray*}
where
\begin{eqnarray*}
[M_i](t) =\int_0^t |q_i(w_r)|^2\dif r\leq \aleph^2 t.
\end{eqnarray*}
Thus, we   arrive at that
\begin{eqnarray*}
&& \mP(\sup_{t\in [0,T]}\sum_{i\in \cZ_0}\big|\int_0^tq_i(w_r)\dif W_i(r)\big|<\eps)
\\ &&=\mP(\sup_{t\in [0,T]}\sum_{i\in \cZ_0 } |M_i(t)|<\eps)
\\ &&=\mP(\sup_{t\in [0,T]}\sum_{i\in \cZ_0 } |B^i_{[M_i](t)}|<\eps)
\\ &&\geq \mP(\sup_{t\in [0,\aleph^2 T]}\sum_{i\in \cZ_0} |B^i_{t}|<\eps)
>0.
\end{eqnarray*}
The above inequality implies the desired result and the proof is complete.

\end{proof}

\begin{proposition}(Weak  form of irreducibility)
\label{17-2}
 Under the Condition \ref{16-2},     for any $\cC,\gamma>0$, there exists a $T=T(\cC,\gamma)>0$ such that
  \begin{eqnarray*}
    \inf_{\|w_0\|\leq \cC}P_{T}(w_0, \cB_\gamma)>0,
  \end{eqnarray*}
  where $\cB_\gamma=\{w\in H,\|w\|\leq \gamma\}.$
\end{proposition}
\begin{proof}
The proof of this proposition is almost the same as that in
 \cite[Lemma 3.1]{EM01}.
 For the convenient  of reading, we still give details here.
 Define $v_t=w_t-\hat f_t$,  where $\hat f_t=\sum_{i\in \cZ_0}\int_0^t
 q_i(w_r)e_i \dif  W_i(r)
  $. Then, $v=v_t$ satisfies
  \begin{eqnarray*}
    \frac{\partial v}{\partial t}=\nu \Delta (v+\hat f)
    +B(\cK  w,w)=\nu \Delta (v+\hat f)
    +B(\cK  w,v+\hat f).
  \end{eqnarray*}
  Taking the $L^2-$inner product of this equation with $v$ produces
  \begin{eqnarray}
  \nonumber   \frac{1}{2}\frac{d}{dt}\|v\|^2 &=& -\nu \|\nabla v\|^2+
    \langle \nu \Delta \hat f,v\rangle+
    \langle B(\cK w, \hat f), v\rangle
    \\ \nonumber
    &\leq &  -\nu \|\nabla v\|^2 +C \|v\| \|\Delta \hat f\|
    +C\|v\|\|\Delta \hat f\| \|\cK w \|_1
     \\  \nonumber
    &= &  -\nu \|\nabla v\|^2 +C \|v\| \|\Delta \hat f\|
    +C\|v\|\|\Delta \hat f\| \|\cK (v+\hat f) \|_1
    \\ \nonumber &\leq & -\nu \|\nabla v\|^2 +C \|v\| \|\Delta \hat f\|
    +C\|v\|\|\Delta \hat f\| \|v \|
    +C\|v\|\|\Delta \hat f\| \|\hat f\|
     \\ \nonumber  &\leq & -\frac{\nu}{2} \|\nabla v\|^2 +\frac{C}{\nu } \|\Delta \hat f\|^2
    +C\|v\|^2 \|\Delta \hat f\|
    +\frac{C}{\nu } \|\Delta \hat f\|^2  \|\hat f\|^2
    \\  \label{1231-1}  &\leq &
   -\frac{\nu}{2} \| v\|^2
    +C\|v\|^2 \|\Delta \hat f\|
    +\frac{C}{\nu } \[\|\Delta \hat f\|^2  \|\hat f\|^2
    +\|\Delta \hat f\|^2
    \].
  \end{eqnarray}
  For    any $\delta\in (0,1) $ and $T>0$,  define
  \begin{eqnarray*}
    \Omega'(\delta,T)&=& \Big\{
    g=(g_s)_{s\in [0,T]}\in C([0,T];H):\sup_{s\in [0,T]}\|\Delta g_s\|\leq \min \{\delta,
    \frac{\nu}{4C}\}\Big\},
  \end{eqnarray*}
  where the constant $C$ is taken from    (\ref{1231-1}).
  If $\hat f\in  \Omega'(\delta,T)$, then by (\ref{1231-1}), one has
  \begin{eqnarray*}
    \|v_t\|^2 \leq  \|v_0\|^2e^{-\frac{\nu}{2} t}
    + \frac{4C }{\nu ^2 }\min\big(\delta,\frac{\nu}{4C}\big)^2.
  \end{eqnarray*}
  In the above, we have used the assumption $\delta\in (0,1).$
  Since $\|w_0\|\leq \cC,$ there exists a  $T$
  and a $\delta$
  such that
  \begin{eqnarray*}
    \|v_T\|\leq \frac{\gamma}{2} \text{ and } \delta \leq \min\{\frac{\gamma}{2},1 \}.
  \end{eqnarray*}
  Putting everything together, for appropriate $T$
  and $\delta$, one has
  \begin{eqnarray*}
    \|w_0\|\leq \cC \text{ and  }
    \hat f\in  \Omega'(\delta,T) \Rightarrow
    \|w_T\|\leq \|v_T\|+\|\hat f_T\|\leq \gamma.
  \end{eqnarray*}
  Combining this fact  with Lemma \ref{16-1}, we complete the proof.
\end{proof}

\begin{appendix}
\section{Proof of Lemma  \ref{18-2}}
\label{A}

First, we list some useful lemmas; then, we demonstrate a proof of Lemma \ref{18-2}.

\begin{lemma}
\label{25-2}
There exists a $C>0$ such that
  \begin{eqnarray}
  \label{25-1}
  \begin{split}
   & \|\tilde B(u,v)\|\leq C\|u\|_1\|v\|_1,
    \\
 &  \big| \langle B(\cK u,v),(-\Delta)^2w\rangle \big|\leq   C  \| u\|_{1/2} \|v\|_3\|w\|_2+C\|u\|_{5/2} \|v\|_1\|w\|_2,
   \\   &  \big| \langle B(\cK u,v),(-\Delta)^3w\rangle \big|
     \leq  C  \| u\|_{1/2} \|v\|_4\|w\|_3+C \| u\|_{7/2} \|v\|_1\|w\|_3.
\end{split}
  \end{eqnarray}
\end{lemma}
\begin{proof}
  (\ref{25-1}) is directly obtained by (\ref{44-1}).
  In the proof of this lemma, for any $\alpha=(\alpha_1,\alpha_2)\in \mN\times \mN,$
 $ D^\alpha:=\partial_{x_1}^{\alpha_1}\partial_{x_2}^{\alpha_2}.$
By interpolating inequality(c.f. \cite[Property 1.1.4]{KS-book}) and  (\ref{44-1}), we conclude that
  \begin{eqnarray*}
  && \langle B(\cK u,v),(-\Delta)^2 w\rangle=\sum_{|\alpha|=2}C_\alpha
  \langle D^\alpha B(\cK u,v),D^\alpha w\rangle
  \\ &&=\sum_{|\alpha|=2}C_\alpha \sum_{|\beta|=2,\beta\leq \alpha}
  \langle  B(D^{\alpha-\beta} \cK u,D^\beta v),D^\alpha w\rangle
  +\sum_{|\alpha|=2}C_\alpha \sum_{|\beta|=1,\beta\leq \alpha}
  \langle  B(D^{\alpha-\beta} \cK u,D^\beta v),D^\alpha w\rangle
  \\ && \quad +\sum_{|\alpha|=2}C_\alpha \sum_{|\beta|=0}
  \langle  B(D^{\alpha-\beta} \cK u,D^\beta v),D^\alpha w\rangle
  \\ &&\leq C\|\cK u\|_{3/2} \|v\|_3\|w\|_2+C \|\cK u\|_{5/2} \|v\|_2\|w\|_2
  +C\|\cK u\|_{7/2} \|v\|_1\|w\|_2
    \\ &&\leq C\| u\|_{1/2} \|v\|_3\|w\|_2+C \| u\|_{3/2} \|v\|_2\|w\|_2
  +C\|u\|_{5/2} \|v\|_1\|w\|_2
  \\ &&\leq   C\| u\|_{1/2} \|v\|_3\|w\|_2+C \| u\|_{1/2}^{1/2} \| u\|_{5/2}^{1/2} \|v\|_1^{1/2}\|v\|_3^{1/2}\|w\|_2
  +C\|u\|_{5/2} \|v\|_1\|w\|_2
  \\ &&\leq C  \| u\|_{1/2} \|v\|_3\|w\|_2+C\|u\|_{5/2} \|v\|_1\|w\|_2.
  \end{eqnarray*}

By interpolating inequality and  (\ref{44-1}), we also get  
  \begin{eqnarray*}
  && \langle B(\cK u,v),(-\Delta)^3 w\rangle=\sum_{|\alpha|=3}C_\alpha
  \langle D^\alpha B(\cK u,v),D^\alpha w\rangle
  \\ &&=\sum_{|\alpha|=3}C_\alpha \sum_{|\beta|=3,\beta\leq \alpha}
  \langle  B(D^{\alpha-\beta} \cK u,D^\beta v),D^\alpha w\rangle
  +\sum_{|\alpha|=3}C_\alpha \sum_{|\beta|=2,\beta\leq \alpha}
  \langle  B(D^{\alpha-\beta} \cK u,D^\beta v),D^\alpha w\rangle
  \\ && \quad +\sum_{|\alpha|=3}C_\alpha \sum_{|\beta|=1,\beta\leq \alpha}
  \langle  B(D^{\alpha-\beta} \cK u,D^\beta v),D^\alpha w\rangle
  +\sum_{|\alpha|=3}C_\alpha \sum_{|\beta|=0}
  \langle  B(D^{\alpha-\beta} \cK u,D^\beta v),D^\alpha w\rangle
  \\ &&\leq C \|\cK u\|_{3/2} \|v\|_4\|w\|_3
  +C\|\cK u\|_{5/2} \|v\|_3\|w\|_3
  \\ &&\quad\quad+\|\cK u\|_{7/2} \|v\|_2\|w\|_3+\|\cK u\|_{9/2} \|v\|_1\|w\|_3
    \\ &&\leq C \| u\|_{1/2} \|v\|_4\|w\|_3
  +C\| u\|_{3/2} \|v\|_3\|w\|_3
  \\ &&\quad\quad+C\| u\|_{5/2} \|v\|_2\|w\|_3+C \| u\|_{7/2} \|v\|_1\|w\|_3
  \\ &&\leq  C \| u\|_{1/2} \|v\|_4\|w\|_3
  +C\| u\|_{1/2}^{2/3} \| u\|_{7/2}^{1/3}\|v\|_1^{1/3}\|v\|_{4}^{2/3}\|w\|_3
  \\ &&\quad\quad+C\| u\|_{1/2}^{1/3}\|u\|_{7/2}^{2/3} C\|v\|_1^{2/3}\|v\|_4^{1/3}\|w\|_3+C \| u\|_{7/2} \|v\|_1\|w\|_3
  \\ &&\leq C  \| u\|_{1/2} \|v\|_4\|w\|_3+C \| u\|_{7/2} \|v\|_1\|w\|_3.
  \end{eqnarray*}
  The proof is complete.

\end{proof}

\begin{lemma}
\label{19-2}
For any  $ T>0,$   there exists  a constant $\vartheta_1=\vartheta_1 (T)\in [T/2,T)$
 only depending on $T$ and $\nu,d,\aleph $
 such that
for any $1\leq n\leq 200, s\in [\vartheta_1,T]$ and $H$-valued random variable $\xi\in \cF_s$, it   holds that
  \begin{eqnarray}
  \label{20-1}
   \begin{split}
   & \mE \big[\sup_{t\in [s,T]}\|J_{s,t}\xi\|_2^{2n}\big|\cF_s\big]
    \\
    &
    \leq  C_{T}\|\xi\|_2^{2n}\exp\{\eta_0 \|w_s\|^2/2  \}\mE\Big[\sup_{t\in [s,T]} (1+   \|w_t\|_{3}^{4n})\big|\cF_s\Big],
    \end{split}
    \end{eqnarray}
    and
    \begin{eqnarray}
    \label{22-2}
    \begin{split}
     &\mE \big[\sup_{t\in [s,T]}\|J_{s,t}\xi\|_3^{2n}\big|\cF_s\big]
    \\ &    \leq   C_{T}\|\xi\|_3^{2n}\exp\{\eta_0 \|w_s\|^2 /2 \}\mE\Big[\sup_{t\in [s,T]}(1+   \|w_t\|_{4}^{4n})\big|\cF_s\Big],
    \end{split}
  \end{eqnarray}
  where $C_{T}$ is a constant   depending on $T,\nu,d,\aleph$ and $\eta_0=\eta_0(\nu,d,\aleph)$ is  given by  Lemma \ref{L:5.1}.
  \end{lemma}
  \begin{proof}
   By (\ref{25-1}),  it holds that
  \begin{eqnarray*}
   && \big|  \langle  B(\cK w_t,J_{s,t}\xi), (-\Delta)^2  J_{s,t}\xi\rangle\big|
   \\ && \leq
    C \| J_{s,t}\xi\|_2 \|J_{s,t}\xi\|_3  \|w_t\|_{1/2}
    +C\| J_{s,t}\xi\|_2 \|J_{s,t}\xi\|_1  \|w_t\|_{5/2}
    \\ &&\leq  \frac{\nu }{4}  \| J_{s,t}\xi\|_3^2+C \|J_{s,t}\xi\|_2^2  \|w_t\|_{1/2}^2  +C \|J_{s,t}\xi\|_1^2   \|w_t\|_{5/2}^2
  \end{eqnarray*}
  and
  \begin{eqnarray*}
   &&\big|  \langle  B(\cK J_{s,t}\xi,w_t), (-\Delta)^2  J_{s,t}\xi\rangle\big|
   \\ && \leq
    C \| J_{s,t}\xi\|_2 \|J_{s,t}\xi\|_{5/2}  \|w_t\|_{1}
    +C  \| J_{s,t}\xi\|_2 \|J_{s,t}\xi\|_{1/2}  \|w_t\|_{3}
    \\ && \leq  \frac{\nu }{4}  \| J_{s,t}\xi\|_3^2+C \|J_{s,t}\xi\|_2^2    \|w_t\|_{1}^{4/3}+ C \|J_{s,t}\xi\|_{1/2}^2   \|w_t\|_{3}^2.
  \end{eqnarray*}
  Using It\^o's formula to
  $\langle J_{s,t}\xi, (-\Delta)^2  J_{s,t}\xi\rangle$,  we arrive at that
   \begin{eqnarray*}
    && \dif  \langle J_{s,t}\xi, (-\Delta)^2  J_{s,t}\xi\rangle
    =2
    \langle  (-\Delta)^2  J_{s,t}\xi, \dif J_{s,t}\xi    \rangle
  +
  \sum_{j\in \cZ_0}|j|^4  |Dq_j(w_t)J_{s,t}\xi|^2  \dif t
    \\ &&\leq
    C \|J_{s,t}\xi\|_2^2\big[ \|w_t\|_{1/2}^2+ \|w_t\|_1^{4/3}\big]\dif t+ C \|J_{s,t}\xi\|_1^2   \|w_t\|_{3}^2\dif t
    \\ &&   \quad +2 \big\langle (-\Delta)^2  J_{s,t}\xi, \sum_{j\in \cZ_0}\big(Dq_j(w_t) J_{s,t}\xi\big)e_j
    \big\rangle \dif W_j(t).
  \end{eqnarray*}
   Therefore, one has
  \begin{eqnarray}
    \nonumber && \|J_{s,t}\xi\|_2^2
    \leq C e^{C \int_s^t(\|w_r\|_{1/2}^2+ \|w_r\|_1^{4/3}) \dif r}
    \|\xi\|_2^2
    \\ && \nonumber  +C e^{C \int_s^t(\|w_r\|_{1/2}^2+ \|w_r\|_1^{4/3}) \dif r}
     \int_s^t  \|J_{s,r}\xi\|_1^2   \|w_r\|_{3}^2
      \dif r
    \\    \nonumber && + 2
    e^{C \int_s^t(\|w_u\|_{1/2}^2+ \|w_u\|_1^{4/3}) \dif u}
    \\  \nonumber && \quad \times  \int_s^t
     e^{-C \int_s^r(\|w_u\|_{1/2}^2+ \|w_u\|_1^{4/3}) \dif u}
    \sum_{j\in \cZ_0}|j|^4 \big(Dq_j(w_r) J_{s,r}\xi\big)
    \big\langle   J_{s,r}\xi, e_j
    \big\rangle \dif W_j(r).
  \end{eqnarray}
  In the above inequality, first  take power of order $n$ and then take conditional  expectation with respect to $\cF_s$.
   For  some constants $C_1=C_1(\nu,d,\aleph), C_T=C(T,\nu,d,\aleph)$ and  any $1\leq n\leq 200,$   by  Cauchy-Schwartz's   inequality and  Burkholder-Davis-Gundy  inequality,  it holds that
\begin{eqnarray*}
  && \mE \Big[ \sup_{t\in [s,T]}\|J_{s,t}\xi\|_2^{2n}\big|\cF_s\Big]
  \leq  C_1  \mE\Big[  e^{C_1  \int_s^T \big(\|w_r\|_{1/2}^2+ \|w_r\|_1^{4/3}\big) \dif r}\big|\cF_s\Big]\|\xi\|_2^{2n}
  \\ && \quad+ C_1 \mE\Big[  e^{C_1  \int_s^T \big(\|w_r\|_{1/2}^2+ \|w_r\|_1^{4/3}\big) \dif r} \big(\int_s^T \|J_{s_1,r}\xi \|_1^2\dif r\big)^n (1+\sup_{r\in [s,T]}\|w_r\|_3^{2n})\big|\cF_s\Big]
  \\ && \quad +C_1 \Big( \mE  e^{C_1  \int_s^T \big(\|w_r\|_{1/2}^2+ \|w_r\|_1^{4/3}\big) \dif r}\Big|\cF_s \Big)^{1/2}
   \Big( \mE \int_s^T \|J_{s,r}\xi \|^{4n}\dif r \Big|\cF_s \Big)^{1/2}.
\end{eqnarray*}
 Thus, by   Lemma \ref{2-3}  and   H\"older's inequality, for  some constant
 $C_T=C(T,\nu,d,\aleph)$ and   any $1\leq n\leq 200,$
  we arrive at
  \begin{eqnarray}
  \label{31-1}
  \begin{split}
  & \mE \Big[ \sup_{t\in [s,T]}\|J_{s,t}\xi\|_2^{2n}\big|\cF_s\Big]
  \leq C_T   \Big(   \mE \big[ e^{4 C_1  \int_s^T \big(\|w_r\|_{1/2}^2+ \|w_r\|_1^{4/3}\big) \dif r}\big|\cF_s\big] \Big)^{1/4 }
  \\   & \quad\quad \times \Big(\|\xi\|_2^{2n}
  +\exp\{\eta_0 \|w_s\|^2/4 \}  \|\xi\|^{2n} \mE \[(1+\sup_{r\in [s,T]}\|w_r\|_3^{4n})\big|\cF_s \] \Big).
  \end{split}
  \end{eqnarray}
Observe that for any $\vartheta_1\in [T/2,T)$ and
$\vartheta_1\leq s \leq T,$ it holds that
  \begin{eqnarray*}
  \begin{split}
   & \exp\big\{4 C_1  \int_s^T(\|w_r\|_{1/2}^2+ \|w_r\|_1^{4/3})\dif r \big\}
   \\ &\leq \exp\big\{\frac{\eta_0}{2}  \int_s^T \|w_r\|_1^2\dif r+  C_2   (T-\vartheta_1)+C_3   \sup_{r\in [s,T]}\|w_r\|^2 (T-\vartheta_1) \big\},
   \end{split}
  \end{eqnarray*}
  where $C_2=C_2(\nu,d,\aleph),C_3=C_3(\nu,d,\aleph).$
  Setting     $\vartheta_1$ close to $T$ enough(for example, $T-\vartheta_1 =\min\{ \frac{\eta_0}{2 C_3},\frac{T}{2}\}),$
  combining the above inequality with (\ref{31-1}) and Lemma \ref{L:5.1}, we get the desired result (\ref{20-1}).

For the proof of (\ref{22-2}),
by (\ref{25-1}),  we obtain
  \begin{eqnarray*}
   &&  \big| \langle  B(\cK w_t,J_{s,t}\xi), (-\Delta)^3  J_{s,t}\xi\rangle\big|
   \\ && \leq
    C \| J_{s,t}\xi\|_3 \|J_{s,t}\xi\|_4  \|w_t\|_{1/2}
    +C\| J_{s,t}\xi\|_3 \|J_{s,t}\xi\|_1  \|w_t\|_{7/2}
    \\ &&\leq  \frac{\nu }{4}  \| J_{s,t}\xi\|_4^2+C \|J_{s,t}\xi\|_3^2  \|w_t\|_{1/2}^2 +C  \|J_{s,t}\xi\|_1^2   \|w_t\|_{7/2}^2
  \end{eqnarray*}
  and
  \begin{eqnarray*}
   && \big|  \langle  B(\cK J_{s,t}\xi,w_t), (-\Delta)^3  J_{s,t}\xi\rangle\big|
   \\ && \leq
    C \| J_{s,t}\xi\|_3 \|J_{s,t}\xi\|_{7/2}  \|w_t\|_{1}
    +C  \| J_{s,t}\xi\|_3 \|J_{s,t}\xi\|_{1/2}  \|w_t\|_{4}
    \\ && \leq  \frac{\nu}{4}  \| J_{s,t}\xi\|_4^2+C \|J_{s,t}\xi\|_3^2    \|w_t\|_{1}^{4/3}+C \|J_{s,t}\xi\|_{1/2}^2   \|w_t\|_{4}^2.
  \end{eqnarray*}
Following similar arguments as that in the proof of (\ref{20-1}),
we get (\ref{22-2}).

  \end{proof}

\begin{lemma}
  \label{22-3}
  For any    $\kappa,T>0, n\geq 1,$    real numbers   $s_1,s_2,t$ with $\frac{T}{2} \leq s_1\leq s_2\leq t \leq T$
  and
  $H$-valued random variable $\xi\in \cF_{s_1}$,
it holds that
\begin{eqnarray*}
 &&  \mE \Big[\|J_{s_1,t}\xi-J_{s_2,t}\xi\|^{2n}\big|\cF_{s_1}\Big]
  \\ && \leq C_{\kappa,n,T}  |s_1-s_2|^{n/2 }\|\xi\|_1^{2n}
  \exp\{\kappa \|w_{s_1}\|^2 \}
  \mE \Big[ \sup_{r\in [s_1,t]}(1+\|w_r\|_1^{4n})\big|\cF_{s_1}\Big],
\end{eqnarray*}
where $C_{\kappa,n,T}$ is a constant  depending on $n,\kappa,T$ and $\nu,d,\aleph.$
  \end{lemma}

\begin{proof}
Denote $\rho_t^1=J_{s_1,t}\xi,\rho_t^2=J_{s_2,t}\xi.$
  Using It\^o's formula to $\langle \rho_t^1-\rho_t^2, \rho_t^1- \rho_t^2 \rangle^n,$ we get
  \begin{eqnarray*}
    && \dif \| \rho_t^1-\rho_t^2\|^{2n}=2n  \| \rho_t^1-\rho_t^2\|^{2n-2}
    \langle \rho_t^1-\rho_t^2, \dif \rho_t^1-\dif  \rho_t^2 \rangle
    \\ &&\quad\quad+ n  \| \rho_t^1-\rho_t^2\|^{2n-2}
    \sum_{j\in \cZ_0} |Dq_j(w_t)(\rho_t^1-\rho_t^2)|^2\dif t
    \\ &&\quad\quad+2n(n-1)\| \rho_t^1-\rho_t^2\|^{2n-4}
    \big( \sum_{j\in \cZ_0}|Dq_j(w_t)(\rho_t^1-\rho_t^2)|^2 \langle  \rho_t^1-\rho_t^2,e_j\rangle^2\big)\dif t.
  \end{eqnarray*}
  Observe that
  \begin{eqnarray*}
  && \langle B(\cK (\rho_t^1-\rho_t^2),w_t),  \rho_t^1- \rho_t^2\rangle
 \leq C \| \rho_t^1-\rho_t^2\| \| \rho_t^1-\rho_t^2\|_{1/2}
  \|w_t\|_1
  \\ &&\leq \frac{\nu}{8} \| \rho_t^1-\rho_t^2\|_{1}^2
  +C  \| \rho_t^1-\rho_t^2\|^2   \|w_t\|_1^{4/3}.
  \end{eqnarray*}
  Hence, we get
  \begin{eqnarray}
   \nonumber  && \mE \big[\| \rho_t^1-\rho_t^2\|^{2n}\big|\cF_{s_2}\big]
 \leq
    C_{n,T} \|J_{s_1,s_2}\xi-\xi\|^{2n} \mE\big[  \exp\{C_n \int_{s_2}^t \|w_s\|_1^{4/3}\dif s\}\big|\cF_{s_2}\big]
    \\ \label{29-4}  &&\leq  C_{\kappa,n,T}  \|J_{s_1,s_2}\xi-\xi\|^{2n} \exp\{\kappa \|w_{s_2}\|^2/2  \} .
  \end{eqnarray}
  So, it remains to  give an estimate of $\mE\big[ \|J_{s_1,s_2}\xi-\xi\|^{4n}\big|\cF_{s_1} \big]. $

Using It\^o's formula to $\varrho_t= \|J_{s_1,t}\xi\|^2,$ we get
\begin{eqnarray*}
 \|J_{s_1,s_2}\xi\|^2&=& \|\xi\|^2-2 \nu\int_{s_1}^{s_2} \|J_{s_1,r}\xi\|_1^2\dif r
+2 \int_{s_1}^{s_2} \langle \tilde B(w_r,J_{s_1,r}\xi),J_{s_1,r}\xi \rangle\dif r
 \\ && +2\sum_{j\in \cZ_0}\int_{s_1}^{s_2} \big(Dq_j(w_r)J_{s_1,r}\xi\big)  \langle J_{s_1,r}\xi, e_j\rangle \dif W_j(r)
 \\ && + \sum_{j\in \cZ_0}\int_{s_1}^{s_2} |Dq_j(w_r)J_{s_1,r}\xi|^2\dif r
\end{eqnarray*}
and
\begin{eqnarray*}
  && \langle J_{s_1,s_2}\xi,\xi\rangle=\|\xi\|^2+\int_{s_1}^{s_2}\langle \nu \Delta  J_{s_1,r}\xi, \xi\rangle \dif r
  +\int_{s_1}^{s_2} \langle \tilde B(w_r,J_{s_1,r}\xi),\xi \rangle\dif r
  \\ &&\quad \quad+ \sum_{j\in \cZ_0}\int_{s_1}^{s_2} \big(Dq_j(w_r)J_{s_1,r}\xi\big)  \langle \xi, e_j\rangle \dif W_j(r).
\end{eqnarray*}
Therefore,
\begin{align}
\nonumber  & \mE \big[\|J_{s_1,s_2}\xi-\xi\|^{4n}\big|\cF_{s_1} \big]
= \mE \Big[\big(\|J_{s_1,s_2}\xi\|^2  +\|\xi\|^2-2 \langle J_{s_1,s_2}\xi, \xi\rangle\big)^{2n} \big|\cF_{s_1} \Big]
\\ \label{AA-2} &  =\mE \Big[\big( I_1(s_1,s_2)+I_2(s_1,s_2)\big)^{2n}   \big|\cF_{s_1} \Big],
\end{align}
where
\begin{eqnarray*}
 I_1(s_1,s_2)&=& -2 \nu\int_{s_1}^{s_2} \|J_{s_1,r}\xi\|_1^2\dif r
 +2 \sum_{j\in \cZ_0}\int_{s_1}^{s_2} \big(Dq_j(w_r)J_{s_1,r}\xi\big)  \langle J_{s_1,r}\xi, e_j\rangle \dif W_j(r)
 \\ && +2 \int_{s_1}^{s_2} \langle \tilde B(w_r,J_{s_1,r}\xi),J_{s_1,r}\xi \rangle\dif r+\sum_{j\in \cZ_0}\int_{s_1}^{s_2} |Dq_j(w_r)J_{s_1,r}\xi|^2\dif r,
 \\   I_2(s_1,s_2)&=&
- 2 \int_{s_1}^{s_2}\langle  \nu \Delta  J_{s_1,r}\xi, \xi\rangle \dif r
  -2 \int_{s_1}^{s_2} \langle \tilde B(w_r,J_{s_1,r}\xi),\xi\rangle\dif r
 \\ && -2  \sum_{j\in \cZ_0}\int_{s_1}^{s_2} \big(Dq_j(w_r)J_{s_1,r}\xi\big)  \langle  \xi, e_j\rangle \dif W_j(r).
\end{eqnarray*}
With the help of
\begin{eqnarray*}
 &&  \langle \tilde B(w_r,J_{s_1,r}\xi),J_{s_1,r}\xi \rangle
\leq C \|J_{s_1,r}\xi\|   \|J_{s_1,r}\xi\|_{1/2}  \|w_r\|_{1},
\\ &&
 \langle \tilde B(w_r,J_{s_1,r}\xi),\xi\rangle \leq
 C \|w_r\|_{1/2}\|J_{s_1,r}\xi\|_1\|\xi\|
 +C\|w_r\|_{1}\|J_{s_1,r}\xi\|_{1/2}\|\xi\|,
\end{eqnarray*}
we get
\begin{eqnarray*}
&& |I_1(s_1,s_2)|+|I_2(s_1,s_2)|
\\ && \leq  C\sup_{r\in [s_1,s_2]}\|J_{s_1,r}\xi\|_1^2 \cdot  (1+
\sup_{r\in [s_1,s_2]}\|w_r\|_1)(s_2-s_1)
\\ &&\quad  + 2 \sum_{j\in \cZ_0}\Big| \int_{s_1}^{s_2} \big(Dq_j(w_r)J_{s_1,r}\xi\big)  \langle J_{s_1,r}\xi, e_j\rangle \dif W_j(r)\Big|
\\&& \quad +2 \sum_{j\in \cZ_0} \Big| \int_{s_1}^{s_2} \big(Dq_j(w_r)J_{s_1,r}\xi\big)  \langle  \xi, e_j\rangle \dif W_j(r)\Big|.
\end{eqnarray*}
Thus, with the help of Lemmas  \ref{19-1},  \ref{2-3}, H\"older's inequality   and Burkholder-Davis-Gundy  inequality,
we arrive at
\begin{align}
\label{29-3}
\begin{split}
 & \mE \big[\big( I_1(s_1,s_2)+I_2(s_1,s_2)\big)^{2n}    \big|\cF_{s_1}\big]
\\    & \leq C_{\kappa,T,n}  |s_2-s_1|^{ n }\|\xi\|_1^{4n}
  \exp\{\kappa \|w_{s_1}\|^2/2  \}
  \mE \Big[ \sup_{r\in [s_1,t]}(1+\|w_r\|_1^{4n})\big|\cF_{s_1}\Big].
  \end{split}
\end{align}

Combining (\ref{29-3}) with (\ref{AA-2})
(\ref{29-4}), we conclude that
\begin{eqnarray*}
&& \mE \big[\| \rho_t^1-\rho_t^2\|^{2n}\big|\cF_{s_1}\big] \leq  C_{\kappa,n,T}
\mE \big[  \|J_{s_1,s_2}\xi-\xi\|^{2n} \exp\{\kappa \|w_{s_2}\|^2/2  \}\big|\cF_{s_1} \big]
\\ &&\leq C_{\kappa,n,T} \Big(\mE \big[  \|J_{s_1,s_2}\xi-\xi\|^{4n} \big|\cF_{s_1} \big]\Big)^{1/2}
\Big(\mE \big[  \exp\{\kappa \|w_{s_2}\|^2  \}\big|\cF_{s_1} \big]\Big)^{1/2}
\\ &&\leq  C_{\kappa,n,T}|s_1-s_2|^{n/2 }  \exp\{\kappa \|w_{s_1}\|^2  \}\|\xi\|_1^{2n}\mE \Big[ \sup_{r\in [s_1,t]}(1+\|w_r\|_1^{4n})\big|\cF_{s_1}\Big].
\end{eqnarray*}
We complete the proof of this lemma.

\end{proof}

\begin{lemma}
  \label{19-3}
  For any  $T>0,$   there exists  a constant $\vartheta_2=\vartheta_2 (T)\in [T/2,T)$
 only depending on $T$ and $\nu,d,\aleph $
 such that
for any     $1\leq n\leq 100,$    real numbers $s_1,s_2,t$ with    $\vartheta_2  \leq s_1\leq s_2\leq  t\leq T$
and
$H$-valued random variable $\xi\in \cF_{s_1}$,
it holds that
\begin{eqnarray*}
 &&  \mE \Big[\|J_{s_1,t}\xi-J_{s_2,t}\xi\|_2^{n}\big|\cF_{s_1}\Big]
  \\ && \leq C_{T} |s_1-s_2|^{n/4 } \exp\{\eta_0 \|w_{s_1}\|^2/2  \}
   \big(\mE \sup_{r\in [s_1,T]}(1+\|w_r\|_{4}^{10n})\big|\cF_{s_1} \big)  \|\xi\|_3^{n},
\end{eqnarray*}
where $C_{T}$ is a constant   depending on $T,\nu,d,\aleph$
  and $\eta_0=\eta_0(\nu,d,\aleph)$ is  given in Lemma \ref{L:5.1}.
  \end{lemma}

\begin{proof}
Denote $\rho_t^1=J_{s_1,t}\xi,\rho_t^2=J_{s_2,t}\xi.$
  Using It\^o's formula to $\langle \rho_t^1-\rho_t^2,(-\Delta)^2 \rho_t^1-(-\Delta)^2 \rho_t^2 \rangle^n,$ we get
  \begin{align}
  \label{23-5}
  \begin{split}
 & \dif \| \rho_t^1-\rho_t^2\|_2^{2n}=2n  \| \rho_t^1-\rho_t^2\|_2^{2n-2}
    \langle (-\Delta)^2 \rho_t^1-(-\Delta)^2 \rho_t^2, \dif \rho_t^1-\dif  \rho_t^2 \rangle
    \\ &\quad\quad+ n  \| \rho_t^1-\rho_t^2\|_2^{2n-2}
    \sum_{j\in \cZ_0} |j|^4 |Dq_j(w_t)(\rho_t^1-\rho_t^2)|^2\dif t
    \\  &\quad\quad+2n(n-1)\| \rho_t^1-\rho_t^2\|_2^{2n-4}
    \sum_{j\in \cZ_0}|j|^8 |Dq_j(w_t)(\rho_t^1-\rho_t^2)|^2 \langle  \rho_t^1-\rho_t^2,e_j\rangle^2 \dif t.
    \end{split}
  \end{align}
  By (\ref{25-1}), Young's inequality  and  interpolating inequality(c.f. \cite[Property 1.1.4]{KS-book}),
   it holds that
  \begin{eqnarray*}
  && \big| \langle B(\cK (\rho_t^1-\rho_t^2),w_t), (-\Delta)^2 \rho_t^1-(-\Delta)^2 \rho_t^2\rangle\big|
 \\ && \leq C \| \rho_t^1-\rho_t^2\|_2 \| \rho_t^1-\rho_t^2\|_{5/2}
  \|w_t\|_1+ C\| \rho_t^1-\rho_t^2\|_2 \| \rho_t^1-\rho_t^2\|_{1/2}
  \|w_t\|_3
  \\ &&\leq \frac{\nu }{16}  \| \rho_t^1-\rho_t^2\|_{3}^2
  +C  \| \rho_t^1-\rho_t^2\|_2^2   \|w_t\|_1^{4/3}
  + C\|  \rho_t^1-\rho_t^2\|_{2}^{5/4}\| \rho_t^1-\rho_t^2\|^{3/4}
  \|w_t\|_3,
\\ &&\leq \frac{\nu }{8}  \| \rho_t^1-\rho_t^2\|_{3}^2
  +C  \| \rho_t^1-\rho_t^2\|_2^2   \|w_t\|_1^{4/3}
  + C\| \rho_t^1-\rho_t^2\|^{2}
  \|w_t\|_3^{8/3},
  \end{eqnarray*}
  and
  \begin{eqnarray*}
   && \big| \langle B(\cK w_t, \rho_t^1-\rho_t^2), (-\Delta)^2 \rho_t^1-(-\Delta)^2 \rho_t^2\rangle\big|
   \\ &&\leq C\|\rho_t^1-\rho_t^2\|_2\|\rho_t^1-\rho_t^2\|_3\|w_t\|_{1/2}
   +C\|\rho_t^1-\rho_t^2\|_2 \|\rho_t^1-\rho_t^2\|_1 \|w_t\|_{5/2}
   \\ &&\leq \frac{\nu }{16}\|\rho_t^1-\rho_t^2\|_3^2+
   C  \|\rho_t^1-\rho_t^2\|_2^2 \|w_t\|_{1/2}^2
   + C\|\rho_t^1-\rho_t^2\|_2^{3/2} \|\rho_t^1-\rho_t^2\|^{1/2}  \|w_t\|_{5/2}
      \\ &&\leq \frac{\nu }{8}\|\rho_t^1-\rho_t^2\|_3^2+
   C  \|\rho_t^1-\rho_t^2\|_2^2 \|w_t\|_{1/2}^2
   + C  \|\rho_t^1-\rho_t^2\|^{2}  \|w_t\|_{5/2}^4.
  \end{eqnarray*}
  Combining the above  estimates  with  (\ref{23-5}), for any $1\leq n\leq 100,$
  we conclude that
  \begin{eqnarray*}
    &&  \dif  \| \rho_t^1-\rho_t^2\|_2^{2n}\leq
     C_1  \| \rho_t^1-\rho_t^2\|_2^{2n}(\|w_t\|_{1/2}^2+\|w_t\|_{1}^{4/3}+1)\dif t
     \\ && ~+ C_1 \| \rho_t^1-\rho_t^2\|^{2n}
     (
  \|w_t\|_3^{4n}+1)\dif t
  \\ && ~+2n \sum_{j\in \cZ_0}  \| \rho_t^1-\rho_t^2\|_2^{2n-2}
    \Big\langle (-\Delta)^2 \rho_t^1-(-\Delta)^2 \rho_t^2, \big(Dq_j(w_t)(\rho_t^1-\rho_t^2)\big)e_j \dif W_j(t)  \Big\rangle,
  \end{eqnarray*}
  where  $C_1=C_1(\nu,d,\aleph)$.
  Therefore, by the above inequality,
  for  any $1\leq n\leq 100,$  it holds that
  \begin{eqnarray}
  \nonumber  && \mE \[ \| \rho_t^1-\rho_t^2\|_2^{2n}\exp\{-C_1\int_{s_2}^t
   (\|w_r\|_{1/2}^2+\|w_r\|_{1}^{4/3}+1) \dif r  \}\big|\cF_{s_2}\]
   \\    \nonumber &&\leq
   \| \rho_{s_2}^1-\rho_{s_2}^2\|_2^{2n}+
   C_1  \mE \Big(\int_{s_2}^t \| \rho_r^1-\rho_r^2\|^{2n}
     (
  \|w_r\|_3^{4n}+1)\dif r  \Big|\cF_{s_2} \Big)
  \\   \nonumber  && \leq \|J_{s_1,s_2}\xi-\xi\|_2^{2n}
  +C_1 \Big( \mE \int_{s_2}^t \| \rho_r^1-\rho_r^2\|^{4n}
     \dif r  \Big|\cF_{s_2} \Big)^{1/2}
\Big( \mE
     \sup_{r\in [s_2,t]}(1 +
  \|w_r\|_3^{8n})  \Big|\cF_{s_2} \Big)^{1/2}.
   \end{eqnarray}
  Furthermore, by  Lemma \ref{22-3} and the above inequality, one gets
  \begin{eqnarray}
  \nonumber  && \mE \[ \| \rho_t^1-\rho_t^2\|_2^{2n}\exp\{-C_1\int_{s_2}^t
   (\|w_r\|_{1/2}^2+\|w_r\|_{1}^{4/3}+1) \dif r  \}\big|\cF_{s_1}\]
  \\   \nonumber && \leq \mE \Big[\|J_{s_1,s_2}\xi-\xi\|_2^{2n}\big|\cF_{s_1}\Big]
  \\ \nonumber &&  +C_1 \Big( \mE \int_{s_2}^t \| \rho_r^1-\rho_r^2\|^{4n}
     \dif r  \Big|\cF_{s_1} \Big)^{1/2}
\Big( \mE
     \sup_{r\in [s_2,t]}(1 +
  \|w_r\|_3^{8n})  \Big|\cF_{s_1} \Big)^{1/2}
  \\    \nonumber &&  \leq   \mE \Big[\|J_{s_1,s_2}\xi-\xi\|_2^{2n}\big|\cF_{s_1}\Big]
   \\ \label{31-2}  &&+ C_T |s_1-s_2|^{n/2} \|\xi\|_3^{2n}  e^{\eta_0\|w_{s_1}\|^2/20 }  \Big(  \mE\sup_{r\in [s_1,t] }
     (1 +
  \|w_r\|_3^{8n})   \Big|\cF_{s_1} \Big).
   \end{eqnarray}
  So, it remains to estimate
 $
  \mE[ \|J_{s_1,s_2}\xi-\xi\|^{2n}_2\big|\cF_{s_1}].
$

Using It\^o's formula to $\varrho_t=\langle J_{s_1,t}\xi, (-\Delta)^2J_{s_1,t}\xi\rangle= \|J_{s_1,t}\xi\|_2^2,$ we get
\begin{eqnarray*}
 \|J_{s_1,s_2}\xi\|_2^2&=& \|\xi\|_2^2-2 \nu\int_{s_1}^{s_2} \|J_{s_1,r}\xi\|_3^2\dif r
 \\ &&+2 \int_{s_1}^{s_2}  \langle \tilde B(w_r,J_{s_1,r}\xi),(-\Delta)^2J_{s_1,r}\xi \rangle\dif r
 \\ && +2 \sum_{j\in \cZ_0}\int_{s_1}^{s_2}   |j|^4 \big(Dq_j(w_r)J_{s_1,r}\xi\big)  \langle J_{s_1,r}\xi, e_j\rangle \dif W_j(r)
 \\ &&+\sum_{j\in \cZ_0}\int_{s_1}^{s_2}   |j|^4 |Dq_j(w_r)J_{s_1,r}\xi|^2\dif r
\end{eqnarray*}
and
\begin{eqnarray*}
  && \langle J_{s_1,s_2}\xi,(-\Delta)^2\xi\rangle=\|\xi\|_2^2+\int_{s_1}^{s_2} \langle  \nu \Delta J_{s_1,r}\xi,  (-\Delta)^2 \xi\rangle \dif r
  +\int_{s_1}^{s_2}  \langle \tilde B(w_r,J_{s_1,r}\xi),(-\Delta)^2\xi \rangle\dif r
  \\ &&\quad \quad+ \sum_{j\in \cZ_0}\int_{s_1}^{s_2}  \big(Dq_j(w_r)J_{s_1,r}\xi\big)  \langle  (-\Delta)^2\xi, e_j\rangle \dif W_j(r).
\end{eqnarray*}
Therefore,
\begin{eqnarray}
\nonumber  && \mE \big[\|J_{s_1,s_2}\xi-\xi\|_2^{2n}\big|\cF_{s_1} \big]
\\ \nonumber && = \mE \Big[\big(\|J_{s_1,s_2}\xi\|_2^2  +\|\xi\|_2^2-2 \langle J_{s_1,s_2}\xi,(-\Delta)^2 \xi\rangle\big)^{n} \big|\cF_{s_1}\Big]
\\ \label{A-2} &&  =\mE \Big[\big( I_1(s_1,s_2)+I_2(s_1,s_2)\big)^n \big|\cF_{s_1}  \Big],
\end{eqnarray}
where
\begin{eqnarray*}
 I_1(s_1,s_2)&=& -2 \nu\int_{s_1}^{s_2} \|J_{s_1,r}\xi\|_3^2\dif r
 +2 \int_{s_1}^{s_2} \langle \tilde B(w_r,J_{s_1,r}\xi),(-\Delta)^2J_{s_1,r}\xi \rangle\dif r
 \\ &&
 +2 \sum_{j\in \cZ_0}|j|^4 \int_{s_1}^{s_2} \big(Dq_j(w_r)J_{s_1,r}\xi\big)  \langle J_{s_1,r}\xi, e_j\rangle \dif W_j(r)
 \\ &&+\sum_{j\in \cZ_0}\int_{s_1}^{s_2} |j|^4 |Dq_j(w_r)J_{s_1,r}\xi|^2\dif r,
 \\   I_2(s_1,s_2)&=&
 -2 \int_{s_1}^{s_2}\langle   \nu \Delta J_{s_1,r}\xi,  (-\Delta)^2 \xi\rangle \dif r
  -2\int_{s_1}^{s_2} \langle \tilde B(w_r,J_{s_1,r}\xi),(-\Delta)^2\xi\rangle\dif r
 \\ && - 2 \sum_{j\in \cZ_0}\int_{s_1}^{s_2} \big(Dq_j(w_r)J_{s_1,r}\xi\big)  \langle  (-\Delta)^2\xi, e_j\rangle \dif W_j(r).
\end{eqnarray*}
With the help of
\begin{eqnarray*}
 &&  \langle \tilde B(w_r,J_{s_1,r}\xi),(-\Delta)^2J_{s_1,r}\xi \rangle
  \\ &&\leq C \|J_{s_1,r}\xi\|_2   \|J_{s_1,r}\xi\|_3  \|w_r\|_{1/2}
  +C  \|J_{s_1,r}\xi\|_2 \|J_{s_1,r}\xi\|_1\|w_r\|_{5/2}
  \\ &&\quad+C \|J_{s_1,r}\xi\|_2   \|J_{s_1,r}\xi\|_{1/2}  \|w_r\|_{3}
  +C  \|J_{s_1,r}\xi\|_2 \|w_r\|_1\|J_{s_1,r}\xi \|_{5/2},
\end{eqnarray*}
we obtain
\begin{eqnarray*}
  && \big| I_1(s_1,s_2)+ I_2(s_1,s_2)\big|\leq (s_2-s_1) \sup_{s\in [s_1,s_2]}\|J_{s_1,r}\xi\|_3^2(1+\|w_r\|_3)
  \\ &&\quad +2 \sum_{j\in \cZ_0}|j|^4 \Big| \int_{s_1}^{s_2} \big(Dq_j(w_r)J_{s_1,r}\xi\big)  \langle J_{s_1,r}\xi, e_j\rangle \dif W_j(r)\Big|
  \\ &&\quad +  2 \sum_{j\in \cZ_0}\Big| \int_{s_1}^{s_2} \big(Dq_j(w_r)J_{s_1,r}\xi\big)  \langle  (-\Delta)^2\xi, e_j\rangle \dif   W_j(r)\Big|.
\end{eqnarray*}
Thus, with the help  of (\ref{A-2}),  H\"older's inequality,  Burkholder-Davis-Gundy  inequality and Lemmas  \ref{19-2},  \ref{2-3},
 for any $1\leq n\leq 100,$
we arrive at
\begin{align}
\nonumber & \mE \big[\|J_{s_1,s_2}\xi-\xi\|_2^{2n}\big|\cF_{s_1} \big] =\mE \big[\big( I_1(s_1,s_2)+I_2(s_1,s_2)\big)^n    \big|\cF_{s_1}\big]
\\ \nonumber &\leq (s_2-s_1)^n \Big( \mE \big[  \sup_{r\in [s_1,s_2]}\|J_{s_1,r}\xi\|_3^{4n}\big|\cF_{s_1}\big]\Big)^{1/2}
\\ \nonumber & \quad\quad \quad\quad\quad\quad \quad  \times  \Big(\mE \big[ \sup_{r\in [s_1,s_2]} (1+\|w_r\|_3)^{2n}  \big|\cF_{s_1}\big]\Big)^{1/2}
\\  \nonumber & \quad   +C \mE\[  \big(\int_{s_1}^{s_2}\|J_{s_1,r}\xi\|^4\dif r\big)^{n/2} \big| \cF_{s_1}\]
\\ \label{28-2}  & \leq C_{T}  \exp\{\eta_0 \|w_{s_1}\|^2/4 \}(s_2-s_1)^{n/2} \mE\big[\sup_{r\in [s_1,T]}(1+\|w_r\|_4^{10 n})\big|\cF_{s_1}\big] \|\xi\|_3^{2n}.
\end{align}

Combining (\ref{28-2})  with  (\ref{31-2}),
 also with the help of Lemma \ref{L:5.1} and H\"older's inequality,  for any $1\leq n\leq 100,$  we get
\begin{align}
\label{31-3}
\begin{split}
    &\mE \[ \| \rho_t^1-\rho_t^2\|_2^{2n}\exp\{-C_1\int_{s_2}^t
   (\|w_r\|_{1/2}^2+\|w_r\|_{1}^{4/3}+1) \dif r  \}\big|\cF_{s_1}\]
  \\ &  \leq
  C_{T}  \exp\{\eta_0 \|w_{s_1}\|^2/4 \}(s_2-s_1)^{n/2} \mE\big[\sup_{r\in [s_1,T]}(1+\|w_r\|_4^{10 n})\big|\cF_{s_1}\big] \|\xi\|_3^{2n}.
  \end{split}
\end{align}
Observe that for any $\vartheta_2\in [T/2,T)$ and
$\vartheta_2\leq s\leq t \leq T,$ it holds that
  \begin{eqnarray*}
  \begin{split}
   & \exp\big\{ C_1   \int_s^t\big(\|w_r\|_{1/2}^2+ \|w_r\|_1^{4/3}\big)\dif r \big\}
   \\ &\leq \exp\{\frac{\eta_0}{4}  \int_s^t \|w_r\|_1^2\dif r+  C_2   (t-s)+C_3    \sup_{r\in [s,t]}\|w_r\|^2 (T-\vartheta_2) \},
   \end{split}
  \end{eqnarray*}
  where  $C_2=C_2(\nu,d,\aleph),C_3=C_3(\nu,d,\aleph)$.
  Setting     $\vartheta_2$ close to $T$ enough(for example, $T-\vartheta_2 =\min\{ \frac{\eta_0}{4 C_3},\frac{T}{2}\}),$
  combining the above inequality with (\ref{31-3}) and Lemma \ref{L:5.1}, we get
  \begin{eqnarray*}
   && \mE \[ \| \rho_t^1-\rho_t^2\|_2^{n}\big|\cF_{s_1}\]\leq
   \Big(\mE\[ \exp\{C_1\int_{s_2}^t
   (\|w_r\|_{1/2}^2+\|w_r\|_{1}^{4/3}+1) \dif r  \}\big|\cF_{s_1}\] \Big)^{1/2}
   \\ && \times \Big( \mE \[ \| \rho_t^1-\rho_t^2\|_2^{2n}\exp\{-C_1\int_{s_2}^t
   (\|w_r\|_{1/2}^2+\|w_r\|_{1}^{4/3}+1) \dif r  \}\big|\cF_{s_1}\] \Big)^{1/2}
   \\ && \leq C_{T}  \exp\{\eta_0 \|w_{s_1}\|^2/2 \}(s_2-s_1)^{n/4} \mE\big[\sup_{r\in [s_1,T]}(1+\|w_r\|_4^{10 n})\big|\cF_{s_1}\big] \|\xi\|_3^{n},
  \end{eqnarray*}
which completes the proof.

\end{proof}

\begin{lemma}
  \label{19-5}
    For any  $T>0,$   let  $\vartheta=\max\{\vartheta_1,\vartheta_2\}\in [T/2,T)$,
    where  $\vartheta_1,\vartheta_2$ are given in Lemma \ref{19-2} and Lemma \ref{19-3},  respectively.
  For any    $1\leq n\leq 100,$
  real numbers   $s,t_1,t_2$  with    $\vartheta \leq s\leq t_1\leq  t_2\leq T$
  and
  $H$-valued random variable  $\xi \in \cF_s$,
it holds that
\begin{eqnarray*}
  && \mE\Big[ \|J_{s,t_1}\xi-J_{s,t_2}\xi\|_2^{2n}\big|\cF_s\Big]
  \\ && \leq C_{T}  |t_1-t_2|^{n/2}\|\xi\|_3^{2n} \exp\{\eta_0 \|w_s\|^2/2   \}\mE \Big(\sup_{r\in [s,T]}(1+\|w_r\|_4^{40  n}) \Big| \cF_s \Big) ,
\end{eqnarray*}
where $C_{T}$ is a constant   depending on $T,\nu,d,\aleph$
  and $\eta_0=\eta_0(\nu,d,\aleph)$ is  given in Lemma \ref{L:5.1}.
  \end{lemma}
\begin{proof}
By (\ref{28-2}),
for any $1\leq n\leq 100,$ it holds that
\begin{eqnarray*}
  && \mE \big[\|J_{t_1,t_2}\xi-\xi\|_2^{2n}\big| \cF_{t_1} \big]
  \\ && \leq
 C_{T}  \exp\{\eta_0 \|w_{t_1}\|^2/4 \}(t_2-t_1)^{n/2} \mE\big[\sup_{r\in [t_1,T]}(1+\|w_r\|_4^{10 n})\big|\cF_{t_1}\big] \|\xi\|_3^{2n}.
\end{eqnarray*}
Thus, by the above inequality,  H\"older's inequality and Lemma  \ref{19-2},   for any $1\leq n\leq 100,$   we conclude that
\begin{eqnarray*}
 && \mE \big[\|J_{s,t_2}\xi-J_{s,t_1}\xi\|_2^{2n}  \big| \cF_s \big]
\\ &&=\mE \Big[ \mE \big[\|J_{t_1,t_2}J_{s,t_1}\xi-J_{s,t_1}\xi\|_2^{2n}  \big|\cF_{t_1}\big]\Big| \cF_s  \Big]
\\ &&\leq C_{T}  \mE \Big[\exp\{\eta_0  \|w_{t_1}\|^2/4 \}\|J_{s,t_1}\xi\|_3^{2n}\sup_{r\in [t_1,T]}(1+\|w_r\|_4^{10n}) \Big| \cF_s \Big] (t_2-t_1)^{n/2}
\\ &&\leq C_{T} (t_2-t_1)^{n/2}
\Big(\mE \exp\{\eta_0  \|w_{t_1}\|^2 \}\big|\cF_{s_1}\Big)^{1/4}
\\ && \quad\quad \times \Big(\mE \|J_{s,t_1}\xi\|_3^{4n}\big|\cF_{s}\Big)^{1/2}
\Big(\mE \sup_{r\in [t_1,T]}(1+\|w_r\|_4^{40n}) \Big| \cF_s \Big)^{1/4}
\\ &&\leq  C_{T}  |t_1-t_2|^{n/2}\|\xi\|_3^{2n} \exp\{\eta_0 \|w_s\|^2 /2  \}\mE \Big(\sup_{r\in [s,T]}(1+\|w_r\|_4^{40  n}) \Big| \cF_s \Big)   .
\end{eqnarray*}
The proof is complete.
\end{proof}

\textbf{Now, we are in a position to prove Lemma \ref{18-2}.}

\begin{proof}
For any  $T>0,$  let  $\vartheta=\max\{\vartheta_1,\vartheta_2\}\in [T/2,T)$
and $\lambda \in [\vartheta,T]$,
    where  $\vartheta_1,\vartheta_2$ are given in Lemma \ref{19-2} and Lemma \ref{19-3},  respectively.
    By Kolmogorov's continuity theorem, we only need to prove
    \begin{align}
     \label{28-4}
      \begin{split}& \mE |f_i(u_1)-f_i(u_2)|^{2n} \leq  C_{T} |u_1-u_2|^{n/2}\exp\{\eta_0 \|w_0\|^2 \},
      \\   & \quad \quad \forall i\in \{1,2,3,4,5,6,7,8\}
      \text{ and }u_1=(t_1,s_1),u_2=(t_2,s_2)\in D_{\lambda,T}
      \end{split}
    \end{align}
 for $n=10, $
 here  $C_{T}$ is a constant   depending on $T,\nu,d,\aleph$  and $\eta_0=\eta_0(\nu,d,\aleph)$ is  given in Lemma \ref{L:5.1}.
  We only prove (\ref{28-4}) for $i=3$ and $i=4$. The proof of the other      cases are similar.

First, we consider the case of $i=3.$
For $n=10$ and any $ u_1=(t_1,s_1),u_2=(t_2,s_2)\in D_{\lambda,T}$ with $s_1\leq s_2,$ we have
\begin{eqnarray}
 \nonumber  && \mE \| J_{s_1,T}\big(\Delta J_{\lambda,t_1}\tilde B(w_\lambda,e_k)\big) -J_{s_2,T}\big(\Delta J_{\lambda,t_2}\tilde B(w_\lambda,e_k)\big) \|^{2n}
\\ \nonumber  && \leq
C \mE \| J_{s_1,T}\big(\Delta J_{\lambda,t_1}\tilde B(w_\lambda,e_k)\big) -J_{s_2,T}\big(\Delta J_{\lambda,t_1}\tilde B(w_\lambda,e_k)\big) \|^{2n}
\\  \nonumber  && +C \mE \| J_{s_2,T}\big(\Delta J_{\lambda,t_1}\tilde B(w_\lambda,e_k)\big) -J_{s_2,T}\big(\Delta J_{\lambda,t_2}\tilde B(w_\lambda,e_k)\big) \|^{2n}
.
\\  \label{21-6} && :=\mathcal{J}_1+\mathcal J_2.
\end{eqnarray}
For the first term $\mathcal J_1$ in (\ref{21-6}),    by H\"older's inequality,
Lemmas  \ref{22-3},  \ref{19-2},  \ref{L:5.1}, we obtain
\begin{eqnarray*}
 \nonumber \mathcal J_1 & = & C \mE \Big( \mE \big(\| J_{s_1,T}\big(\Delta J_{\lambda,t_1}\tilde B(w_\lambda,e_k)\big) -J_{s_2,T}\big(\Delta J_{\lambda,t_1}\tilde B(w_\lambda,e_k)\big) \|^{2n}  \big|\cF_{s_1}\big)\Big)
\\ \nonumber  &\leq & C_T  |u_1-u_2|^{n/2}\mE\Big[ \big\|\Delta J_{\lambda,t_1}\tilde B(w_\lambda,e_k)\big\|_{1}^{2n}
 \exp\{\eta_0  \|w_{s_1}\|^2/4 \}
   \sup_{r\in [s_1,T]}(1+\|w_r\|_1^{4n})\Big]
   \\ \nonumber  &\leq&
    C_T  |u_1-u_2|^{n/2}\Big[ \mE \big\|\Delta J_{\lambda,t_1}\tilde B(w_\lambda,e_k)\big\|_{1}^{4n}\Big]^{1/2}
   \\ && \quad\quad\quad\quad \times \Big[\mE
 \exp\{ \eta_0 \|w_{s_1}\|^2/2  \}
   \sup_{r\in [s_1,T]}(1+\|w_r\|_1^{8n})\Big]^{1/2}.
\\ &&  \leq C_{T} |u_1-u_2|^{n/2}\exp\{\eta_0  \|w_0\|^2\}.
\end{eqnarray*}
For the second  term  $\mathcal J_2$ in (\ref{21-6}),
by Lemma \ref{2-3}, Lemma \ref{19-5} and Lemma \ref{L:5.1},
\begin{align*}
& \mathcal J_2
\\ &
=C \mE\Big(\mE \big(\| J_{s_2,T}\big(\Delta J_{\lambda,t_1}\tilde B(w_\lambda,e_k)\big) -J_{s_2,T}\big(\Delta J_{\lambda,t_2}\tilde B(w_\lambda,e_k)\big) \|^{2n} \big|\cF_{s_2}\big)\Big)
\\ &\leq C_{T}\mE\Big(\exp\{\eta_0  \|w_{s_2}\|^2/4 \}\|\Delta J_{\lambda,t_1}\tilde B(w_\lambda,e_k)-\Delta J_{\lambda,t_2}\tilde B(w_\lambda,e_k)\|^{2n}
\Big)
\\ &\leq C_{T} \Big( \mE \exp\{\eta_0  \|w_{s_2} \|^2/2 \}\Big)^{1/2}
\Big(\mE \|\Delta J_{\lambda,t_1}\tilde B(w_\lambda,e_k)-\Delta J_{\lambda,t_2}\tilde B(w_\lambda,e_k)\|^{4n}\Big)^{1/2}
\\ &\leq C_{T} |u_1-u_2|^{n/2}\exp\{\eta_0  \|w_0\|^2/4 \}
\\ & \quad\quad \times \Big[ \mE\| \tilde B(w_\lambda,e_k)\|_3^{4n}  \exp\{\eta_0  \|w_\lambda\|^2/2   \}\sup_{r\in [\lambda,T]}(1+\|w_r\|_4^{80 n})  \Big]^{1/2}
\\ &\leq C_{T} |u_1-u_2|^{n/2}\exp\{\eta_0  \|w_0\|^2 \}.
\end{align*}
Combining the above estimates of $\mathcal J_1,\mathcal J_2$ with  (\ref{21-6}),  for $n=10,$   one arrives at that
\begin{eqnarray*}
  \mE \|f_3(u_1)-f_3(u_2)\|^{2n}\leq C_{T} |u_1-u_2|^{n/2}\exp\{\eta_0 \|w_0\|^2\}.
\end{eqnarray*}

Now, we consider the case $i=4.$
For  $ n= 10 $ and any  $u_1=(t_1,s_1),u_2=(t_2,s_2)\in D_{\lambda,T}$ with $s_1\leq s_2,$ we have
\begin{eqnarray}
\nonumber  &&  \mE \| J_{s_1,T}\tilde B(w_{t_1},J_{\lambda,t_1}\tilde B(w_\lambda,e_k))
  -J_{s_2,T}\tilde B(w_{t_2},J_{\lambda,t_2}\tilde B(w_\lambda,e_k))\|^{2n}
  \\ \nonumber  &&\leq C \mE \| J_{s_1,T}\tilde B(w_{t_1},J_{\lambda,t_1}\tilde B(w_\lambda,e_k))
  -J_{s_2,T}\tilde B(w_{t_1},J_{\lambda,t_1}\tilde B(w_\lambda,e_k))\|^{2n}
  \\ \nonumber  &&+ C \mE \| J_{s_2,T}\tilde B(w_{t_1},J_{\lambda,t_1}\tilde B(w_\lambda,e_k))-J_{s_2,T}\tilde B(w_{t_2},J_{\lambda,t_2}\tilde B(w_\lambda,e_k))\|^{2n}
  \\  &&:=C \cI_1+C \cI_2. \label{21-10}
\end{eqnarray}
For the first term in (\ref{21-10}),  by Lemmas
\ref{22-3},  \ref{19-2},  \ref{L:5.1} and H\"older's inequality,    we obtain
\begin{eqnarray*}
   &&\cI_1 = \mE\Big(\mE\big(\| J_{s_1,T}\tilde B(w_{t_1},J_{\lambda,t_1}\tilde B(w_\lambda,e_k))
  -J_{s_2,T}\tilde B(w_{t_1},J_{\lambda,t_1}\tilde B(w_\lambda,e_k))\|^{2n} \big|\cF_{s_1}\big)\Big)
  \\ &&
  \leq C_{T}  |s_1-s_2|^{n/2}  \mE \Big[  \|\tilde B(w_{t_1},J_{\lambda,t_1}\tilde B(w_\lambda,e_k)) \|_1^{2n}
  \exp\{\eta_0  \|w_{s_1}\|^2/4  \}
\sup_{r\in [s_1,T]}(1+\|w_r\|_1^{4n}) \Big]
\\ && \leq C_{T}  |s_1-s_2|^{n/2}  \mE \Big[  \|w_{t_1}\|_2^{2n}\|J_{\lambda,t_1}\tilde B(w_\lambda,e_k)) \|_2^{2n}
  \exp\{\eta_0  \|w_{s_1}\|^2/4  \}
\sup_{r\in [s_1,T]}(1+\|w_r\|_1^{4n}) \Big]
\\ && \leq C_{T} |u_1-u_2|^{n/2}\exp\{\eta_0 \|w_0\|^2\}.
\end{eqnarray*}
Now we consider   the second  term in (\ref{21-10}).    By
Lemmas  \ref{2-3},  \ref{25-2}, we obtain
\begin{eqnarray*}
 &&\cI_2 = \mE \Big(\mE \big[\| J_{s_2,T}\tilde B(w_{t_1},J_{\lambda,t_1}\tilde B(w_\lambda,e_k))-J_{s_2,T}\tilde B(w_{t_2},J_{\lambda,t_2}\tilde B(w_\lambda,e_k))\|^{2n}\big|\cF_{s_2}\Big)
\\ &&\leq C_{T}\mE \Big[\exp\{\eta_0  \|w_{s_2}\|^2/4 \}
\|\tilde B(w_{t_1},J_{\lambda,t_1}\tilde B(w_\lambda,e_k))- \tilde B(w_{t_2},J_{\lambda,t_2}\tilde B(w_\lambda,e_k))\|^{2n}\Big]
\\ &&\leq   C_{T}\mE \Big[\exp\{\eta_0  \|w_{s_2}\|^2/4 \}
\|w_{t_1}\|_1^{2n}\|J_{\lambda,t_1}\tilde B(w_\lambda,e_k))- J_{\lambda,t_2}\tilde B(w_\lambda,e_k))\|_1^{2n}\Big]
\\ &&\quad+C_{T}  \mE \Big[\exp\{\eta_0  \|w_{s_2}\|^2/4 \}
\|w_{t_1}-w_{t_2}\|_1^{2n}\|J_{\lambda,t_2}\tilde B(w_\lambda,e_k))\|_1^{2n}\Big]
\\ &&:=\cI_{21}+\cI_{22}.
\end{eqnarray*}
First, we consider   the term $\cI_{21}.$ By
  Lemmas  \ref{19-5},   \ref{L:5.1} and H\"older's inequality,  we get
\begin{eqnarray*}
&& \cI_{21}
\\ && \leq   C_{T}\Big(\mE e^{\eta_0 \|w_{s_2}\|^2/2 }
\|w_{t_1}\|_1^{4n}\Big)^{1/2}\Big(\mE\|J_{\lambda,t_1}\tilde B(w_\lambda,e_k))- J_{\lambda,t_2}\tilde B(w_\lambda,e_k))\|_1^{4n}\Big)^{1/2}
\\ &&\leq   C_{T}\Big(\mE \exp\{\eta_0  \|w_{s_2}\|^2/2 \}
\|w_{t_1}\|_1^{4n}\Big)^{1/2}
\\ && \quad\quad \times  \Big(
|u_1-u_2|^{n}
\mE \big(\exp\{\eta_0  \|w_\lambda\|^2/2  \} \sup_{r\in [\lambda,T]}(1+\|w_r\|_4^{80n})\|\tilde B(w_\lambda,e_k)\|_3^{4n}\big)   \Big)^{1/2}
\\ && \leq C_{T} |u_1-u_2|^{n/2}\exp\{\eta_0 \|w_0\|^2\}.
\end{eqnarray*}
In the next, we consider   the term $\cI_{22}.$ By  Lemmas  \ref{19-1}, \ref{L:5.1} and H\"older's inequality, it holds that
 \begin{eqnarray*}
 &&\cI_{22}\leq    C_{T} \Big( \mE e^{ \eta_0  \|w_{s_2}\|^2 }\Big)^{1/4}
\Big(\mE \|J_{\lambda,t_2}\tilde B(w_\lambda,e_k))\|_1^{8n}\Big)^{1/4}
\Big(\mE \|w_{t_1}-w_{t_2}\|_1^{4n}\Big)^{1/2}
\\ &&\leq
C_{T} \Big( \mE e^{ \eta_0 \|w_{s_2}\|^2 }\Big)^{1/4}
\Big(\mE e^{\eta_0  \|w_\lambda\|^2/4 } \|\tilde B(w_\lambda,e_k)\|_1^{8n}\Big)^{1/4}
\Big(\mE \|w_{t_1}-w_{t_2}\|_1^{4n}\Big)^{1/2}
\\ && \leq C_{T} |u_1-u_2|^{n/2}\exp\{\eta_0 \|w_0\|^2\}.
\end{eqnarray*}
Combining the above estimates of $\cI_1,\cI_{21}, \cI_{22}$ with  (\ref{21-10}),    for $n=10,$ one arrives at that
\begin{eqnarray*}
  \mE \|f_4(u_1)-f_4(u_2)\|^{2n}\leq C_{T} |u_1-u_2|^{n/2}\exp\{\eta_0 \|w_0\|^2\}.
\end{eqnarray*}
We complete the proof of  (\ref{28-4}) for $i=4.$
\end{proof}

\section{Proof of Proposition   \ref{1-2}. }
\label{B}

We will prove Proposition  \ref{1-2} by contradiction.

Suppose that   Proposition  \ref{1-2} { were}  not true,  then   there exist  sequences $\{w_0^{(k)} \} \subseteq  B_H(\mathfrak{R}), \{\eps_k\}\subseteq (0,1)$  and a positive number $\delta_0$  such that
\begin{eqnarray}
\label{D-1}
\lim_{k\rightarrow \infty} \mP(  X^{w_0^{(k)},\alpha,N}<\eps_k)\geq \delta_0>0  \text{ and  } \lim_{k\rightarrow \infty}\eps_k=0.
\end{eqnarray}
Our  aim is to find something
which contradicts (\ref{D-1}).

Since $H$ is a Hilbert space, there exists a subsequence $\{w_0^{(n_k)},k\geq 1\}$ of  $\{w_0^{(k)},k\geq 1\}$
such that    $w_0^{(n_k)}$ converges weakly to  some  $ w_0^{(0)}\in H.$
We still denote this subsequence  by $\{w_0^{(k)},k\geq 1\}.$
Let $w_t^{(k)}$ denote the solution of equation (\ref{0.1}) with $w_t|_{t=0}=w_0^{(k)}(k\geq 0)$.
For any $\xi \in H,$ in the following  equation
    \begin{eqnarray}
    \label{9-2}
    \begin{split}
	 & J_{s,t}\xi=\xi+  \int_s^t \big[ \nu \Delta J_{s,r}\xi +\tilde{B}(w_r,J_{s,r}\xi)\big]\dif r
\\ & \quad\quad\quad\quad +\sum_{j\in \cZ_0}  \int_s^t \big(Dq_j(w_r)J_{s,r}\xi\big)  e_j \dif W_j(r) ,
  \end{split}
\end{eqnarray}
when $w_t$ is replaced by $w_t^{(k)}$, we denote its solution by $J_{s,t}^{(k)}\xi$.
We denote  the adjoint of  $J_{s,t}^{(k)}$ by $K_{s,t}^{(k)}.$

Recall that  for any $M\in \mN,$ ${H}_M=\{e_j: j\in \mZ_*^2  \text{ and } |j|\leq M\},$
${P}_M$  denotes the orthogonal projections from $H$ onto ${H}_M$
and   ${Q}_M u:=u-{P}_Mu, \forall u\in H.$
As before,  $C$  denotes   a constant depending
     $ \nu,\aleph,d,$
   $C_{\mathfrak{R}}$  denotes  a constant depending on $\mathfrak{R}$
   and   $ \nu,\aleph,d.$
   The values of these  constants may change from line to line.

\begin{lemma}
\label{E-5}
For any $t\geq 0$, $k\in \mN$ and $M>\max\{|j|: j\in \cZ_0\},$ one has
  \begin{eqnarray}
  \label{E-1}
  \begin{split}
    \|Q_Mw_t^{(k)}\|^2& \leq e^{-\nu M^2 t }\|Q_Mw_0^{(k)}\|^2+
   \\ &
 \frac{C }{M^{1/2} }\Big(\int_0^t \|w_s^{(k)}\|_1^{4/3}\dif s\Big)^{3/4} \sup_{s\in [0,t]} \|w_s^{(k)}\|^{3}
 \end{split}
  \end{eqnarray}
  and
     \begin{eqnarray}
     \label{E-2}
     \begin{split}
      & \mE \sup_{t\in [0,1] }\|P_M w_t^{(k)}-P_M w_t^{(0)} \| \leq
C_{\mathfrak{R}} \|P_M w_0^{(k)}-P_M w_0^{(0)}\|^{1/2}+\frac{C_{\mathfrak{R}}, }{M^{1/8} } .
     \end{split}
   \end{eqnarray}

\end{lemma}
\begin{proof}
  First, we give a proof of  (\ref{E-1}).
   By (\ref{44-1}),
  one has
  \begin{eqnarray*}
\big| \langle  B (\cK w_t^{(k)},w_t^{(k)}),Q_M w_t^{(k)}\rangle \big| & \leq &
C\|Q_M  w_t^{(k)}\|_1\|w_t^{(k)}\|_{1/2}\|w_t^{(k)}\|
\\ & \leq & \frac{\nu}{4} \|Q_M w_t^{(k)}\|_1^2+C \|w_t^{(k)}\|_1 \|w_t^{(k)}\|^3.
  \end{eqnarray*}
Thus, in view of the equation (\ref{0.1}),  we obtain
\begin{eqnarray*}
 \dif \|Q_Mw_t^{(k)}\|^2 &=&  -2 \nu \|Q_Mw_t^{(k)}\|_1^2+ \langle  B (\cK w_t^{(k)},w_t^{(k)}),Q_M w_t^{(k)}\rangle
 \\ &\leq & -\nu M^2  \|Q_M w_t^{(k)}\|^2\dif t+C \|w_t^{(k)}\|_1 \|w_t^{(k)}\|^3\dif t.
\end{eqnarray*}
 It follows that
\begin{eqnarray*}
&& \|Q_Mw_t^{(k)}\|^2  \leq  e^{-\nu M^2 t} \|Q_Mw_0^{(k)}\|^2
+C \int_0^t e^{-\nu M^2 (t-s) }  \|w_s^{(k)}\|_1 \|w_s^{(k)}\|^3\dif s
\\ & & \leq   e^{-\nu M^2 t} \|Q_Mw_0^{(k)}\|^2
\\ &&\quad +
C \big(\int_0^t e^{-4\nu M^2 (t-s) }  \dif s \big)^{1/4}
\big(\int_0^t \|w_s^{(k)}\|_1^{4/3}\dif s\big)^{3/4} \sup_{s\in [0,t]} \|w_s^{(k)}\|^{3}.
\end{eqnarray*}
This completes the   proof of  (\ref{E-1}).

Next, we will prove  (\ref{E-2}).
One easily sees that
\begin{eqnarray}
\label{E-3}
\begin{split}
& \dif  \|P_M w_t^{(k)}-P_M w_t^{(0)}\|^2
\\ & = -2\nu \|P_M w_t^{(k)}-P_M w_t^{(0)}\|_1^2\dif t
  \\ &\quad  +  2 \langle P_M w_t^{(k)}-P_M w_t^{(0)},
      B(\cK w_t^{(k)},w_t^{(k)})- B(\cK w_t^{(0)},w_t^{(0)})   \rangle
      \\& \quad +\sum_{j\in \cZ_0}|q_j(w_t^{(k)})-q_j(w_t^{(0)})|^2\dif t
      \\ & \quad +2\sum_{j\in \cZ_0}(q_j(w_t^{(k)})-q_j(w_t^{(0)}))\langle
      e_j, P_M w_t^{(k)}-P_M w_t^{(0)} \rangle \dif W_j(t).
      \end{split}
\end{eqnarray}
Clearly, we have
   \begin{eqnarray*}
     && \langle P_M w_t^{(k)}-P_M w_t^{(0)},
      B(\cK w_t^{(k)},w_t^{(k)})- B(\cK w_t^{(0)},w_t^{(0)})   \rangle
      \\ && = \langle P_M w_t^{(k)}-P_M w_t^{(0)},
      B(\cK w_t^{(k)}-\cK w_t^{(0)},w_t^{(k)})   \rangle
     \\ && \quad  + \langle P_M w_t^{(k)}-P_M w_t^{(0)},
    B(\cK w_t^{(0)},  w_t^{(k)}-w_t^{(0)})   \rangle
    \\ && :=I_1+I_2.
   \end{eqnarray*}
   For the term $I_1$, we have
   \begin{eqnarray*}
     I_1 &\leq &  C \|P_M w_t^{(k)}-P_M w_t^{(0)} \|_{1/2}
     \| P_M w_t^{(k)}- P_Mw_t^{(0)} \| \|w_t^{(k)}\|_1
     \\ && \quad \quad +C \|P_M w_t^{(k)}-P_M w_t^{(0)} \|_{1}
     \| Q_M w_t^{(k)}- Q_M w_t^{(0)} \|  \|w_t^{(k)}\|_{1/2}
     \\ &\leq &\frac{\nu}{6}\|P_M w_t^{(k)}-P_M w_t^{(0)} \|_{1}^2+C
     \|P_M w_t^{(k)}-P_M w_t^{(0)} \|^2  \|w_t^{(k)}\|_1^{4/3}
     \\ && \quad \quad +C \| Q_M w_t^{(k)}- Q_M w_t^{(0)} \|^2
     \|w_t^{(k)}\|^2_{1/2}.
   \end{eqnarray*}
    Obviously,
   \begin{eqnarray*}
     I_2& \leq&
      C\|P_M w_t^{(k)}-P_M w_t^{(0)} \|_{1}
     \| Q_M w_t^{(k)}- Q_M w_t^{(0)} \| \|w_t^{(0)}\|_{1/2}
     \\ &\leq & \frac{\nu}{6}\|P_M w_t^{(k)}-P_M w_t^{(0)} \|_{1}^2+C
     \| Q_M w_t^{(k)}- Q_M w_t^{(0)} \|^2 \|w_t^{(0)}\|_{1/2}^2.
   \end{eqnarray*}
 Combining the estimates of $I_1$,$I_2$ with  (\ref{E-3}),
 also in view  of
 \begin{eqnarray*}
 && \sum_{j\in \cZ_0}|q_j(w_t^{(k)})-q_j(w_t^{(0)})|^2\leq d \aleph^2 \|w_t^{(k)}-w_t^{(0)}\|^2
 \\ && =d \aleph^2 \|P_M w_t^{(k)}-P_Mw_t^{(0)}\|^2+d \aleph^2 \|Q_M w_t^{(k)}-Q_Mw_t^{(0)}\|^2,
 \end{eqnarray*}
  we obtain
  \begin{eqnarray}
  \label{14-1}
  \begin{split}
  & \|P_M w_t^{(k)}-P_M w_t^{(0)}\|^2 e^{-C\int_0^t\|w_s^{(k)}\|_1^{4/3} \dif s-Ct }
 \\ & \leq
  \|P_M w_0^{(k)}-P_M w_0^{(0)}\|^2+\mathscr{M}_t
    \\   &  \quad +  \int_0^t \| Q_M w_s^{(k)}- Q_M w_s^{(0)} \|^2 (\|w_s^{(k)}\|_{1/2}^2 +\|w_s^{(0)}\|_{1/2}^2+C)\dif s,
\end{split}
  \end{eqnarray}
  where $$\mathscr{M}_t=\int_0^t 2e^{-\int_0^s\|w_r^{(k)}\|_1^{4/3} \dif r-Cs}\sum_{j\in \cZ_0}(q_j(w_s^{(k)})-q_j(w_s^{(0)}))\langle
      e_j, P_M w_s^{(k)}-P_M w_s^{(0)} \rangle \dif W_j(s).$$
Thus, by (\ref{E-1}) and Lemma  \ref{L:5.1},
also  taking into account the fact that  $\|w_0^{(k)}\|\leq \mathfrak{R},$
it  holds that
 \begin{eqnarray*}
 && \mE\big[ \|P_M w_t^{(k)}-P_M w_t^{(0)}\|^2 e^{-C\int_0^t\|w_s^{(k)}\|_1^{4/3} \dif s-Ct }\big]
 \\ &&\leq  C_{\mathfrak{R}} \|P_M w_0^{(k)}-P_M w_0^{(0)}\|^2+\frac{C_{\mathfrak{R}} }{M^{1/2} }, \quad \forall t\in [0,1].
 \end{eqnarray*}
 By the above inequality, (\ref{14-1}), (\ref{E-1}) and  Burkholder-Davis-Gundy inequality, one arrives at that
 \begin{eqnarray*}
  && \mE  \sup_{t\in [0,1 ]}\Big(\|P_M w_t^{(k)}-P_M w_t^{(0)}\|^2 e^{-C\int_0^t\|w_s^{(k)}\|_1^{4/3} \dif s-Ct }\Big)
  \\ &&\leq  \|P_M w_0^{(k)}-P_M w_0^{(0)}\|^2
 + \mE  \int_0^1  \| Q_M w_s^{(k)}- Q_M w_s^{(0)} \|^2 (\|w_s^{(k)}\|_{1/2}^2 +\|w_s^{(0)}\|_{1/2}^2+C)\dif s
 \\&&+\mE \sup_{t\in [0,1]}\Big|\int_0^t 2e^{-C \int_0^s\|w_r^{(k)}\|_1^{4/3} \dif r-Cs}\sum_{j\in \cZ_0}(q_j(w_s^{(k)})-q_j(w_s^{(0)}))\langle
      e_j, P_M w_s^{(k)}-P_M w_s^{(0)} \rangle \dif W_j(s) \Big|
      \\    && \leq  C_{\mathfrak{R}} \|P_M w_0^{(k)}-P_M w_0^{(0)}\|^2
 + \frac{C_{\mathfrak{R}}}{M^{1/2}}
 \\&&+C_{\mathfrak{R}} \mE \sum_{j\in \cZ_0} \Big(\int_0^1  \Big| e^{-C \int_0^s\|w_r^{(k)}\|_1^{4/3} \dif r-Cs}(q_j(w_s^{(k)})-q_j(w_s^{(0)}))\langle
      e_j, P_M w_s^{(k)}-P_M w_s^{(0)} \rangle\Big|^2  \dif s \Big)^{1/2}
         \\    && \leq C_{\mathfrak{R}}  \|P_M w_0^{(k)}-P_M w_0^{(0)}\|^2
 + \frac{C_{\mathfrak{R}} }{M^{1/2}}
 \\ &&   \quad +C_{\mathfrak{R}}  \Big( \mE  \int_0^1
     \| P_M w_s^{(k)}-P_M w_s^{(0)} \|^2e^{-2C \int_0^s\|w_r^{(k)}\|_1^{4/3} \dif r-2Cs}     \dif s \Big)^{1/2}
     \\ && \leq C_{\mathfrak{R}}  \|P_M w_0^{(k)}-P_M w_0^{(0)}\|
 + \frac{C_{\mathfrak{R}}}{M^{1/4}}.
 \end{eqnarray*}
 (\ref{E-2}) is directly obtained by the above inequality, Cauchy-Schwartz inequality and Lemma \ref{L:5.1}.

  \end{proof}

 \begin{lemma}
  \label{9-3}
For any $\frac{1}{2}\leq s\leq t\leq 1,$  $k,M\in \mN,\xi\in H$ and $M>\max\{|j|: j\in \cZ_0\},$   one has
   \begin{eqnarray*}
    && \mE \|J_{s,t}^{(k)}\xi-J_{s,t}^{(0)}\xi\|
    \leq  C_{\mathfrak{R}} \|P_M w_0^{(k)}-P_M w_0^{(0)}\|^{1/16 }\|\xi\|+\frac{C_{\mathfrak{R}} }{M^{1/64} }\|\xi\|.
    \end{eqnarray*}
  \end{lemma}
  \begin{proof}
    By the equation (\ref{9-2}), one easily  sees that
    \begin{eqnarray}
      \nonumber && \dif \|J_{s,t}^{(k)}\xi-J_{s,t}^{(0)}\xi\|^2
    \leq   -
     2  \nu\|J_{s,t}^{(k)}\xi-J_{s,t}^{(0)}\xi\|^2_1\dif t
      \\  \nonumber &&\quad +2 \langle J_{s,t}^{(k)}\xi-J_{s,t}^{(0)}\xi,
      {B}(\cK w_t^{(k)},J_{s,t}^{(k)}\xi)
      -{B}(\cK w_t^{(0)},J_{s,t}^{(0)}\xi) \rangle \dif t
      \\ \nonumber &&\quad+  2 \langle J_{s,t}^{(k)}\xi-J_{s,t}^{(0)}\xi,
      {B}(\cK J_{s,t}^{(k)}\xi, w_t^{(k)})
      -{B}(\cK J_{s,t}^{(0)}\xi, w_t^{(0)}) \rangle \dif t
      \\ \nonumber && \quad+\sum_{j\in \cZ_0 }
     \big| Dq_j(w_t^{(k)} )J_{s,t}^{(k)}\xi- Dq_j(w_t^{(0)} )J_{s,t}^{(0)}\xi\big|^2\dif t
      \\ \nonumber && \quad+ 2 \sum_{j\in \cZ_0 }
      \Big(Dq_j(w_t^{(k)} )J_{s,t}^{(k)}\xi- Dq_j(w_t^{(0)} )J_{s,t}^{(0)}\xi \Big)
      \langle e_j, J_{s,t}^{(k)}\xi-J_{s,t}^{(0)}\xi \rangle\dif W_j(t)
      \\ \label{E-4} &&:=2I_1(t)\dif t+2I_2(t)\dif t+2I_3(t)\dif t+2I_4(t)\dif t+2\sum_{j\in \cZ_0}I_{5j}(t) \dif W_j(t).
    \end{eqnarray}
 For the terms $I_2(t),I_3(t)$, we have
    \begin{eqnarray*}
      I_2(t)&= &
       \langle J_{s,t}^{(k)}\xi-J_{s,t}^{(0)}\xi,
      {B}(\cK w_t^{(k)}-\cK w_t^{(0)},J_{s,t}^{(k)}\xi)
       \rangle
       \\ &\leq &
      C  \|w_t^{(k)}- w_t^{(0)}\| \| J_{s,t}^{(k)}\xi-J_{s,t}^{(0)}\xi\|_1 \|J_{s,t}^{(k)}\xi\|_{1/2}
    \\ &\leq & C  \|w_t^{(k)}- w_t^{(0)}\|^2  \|J_{s,t}^{(k)}\xi\|_{1/2}^2+\frac{\nu}{6}
    \|J_{s,t}^{(k)}\xi-J_{s,t}^{(0)}\xi\|_1^2
    \end{eqnarray*}
    and
    \begin{eqnarray*}
      I_3(t) &= &
      \langle J_{s,t}^{(k)}\xi-J_{s,t}^{(0)}\xi,
      {B}(\cK J_{s,t}^{(k)}\xi-\cK J_{s,t}^{(0)}\xi, w_t^{(k)})
      +{B}(\cK J_{s,t}^{(0)}\xi, w_t^{(k)}- w_t^{(0)}) \rangle
      \\ &\leq & C \|w_t^{(k)}\|_1 \| J_{s,t}^{(k)}\xi- J_{s,t}^{(0)}\xi\|^{3/2}
      \| J_{s,t}^{(k)}\xi- J_{s,t}^{(0)}\xi\|_1^{1/2}
      \\ && +C\|w_t^{(k)} - w_t^{(0)}\|\| J_{s,t}^{(k)}\xi- J_{s,t}^{(0)}\xi\|_1\|J_{s,t}^{(0)}\xi \|_{1/2}
    \\ &\leq & \frac{\nu}{6}\| J_{s,t}^{(k)}\xi- J_{s,t}^{(0)}\xi\|_1^{2}
    +C \| J_{s,t}^{(k)}\xi- J_{s,t}^{(0)}\xi\|^{2} \|w_t^{(k)}\|_1^{4/3}
    \\ && +C \|w_t^{(k)}- w_t^{(0)}\|^2\|J_{s,t}^{(0)}\xi\|_{1/2}^2.
    \end{eqnarray*}
    For the term $I_4(t),$ it holds that
    \begin{eqnarray*}
      I_4(t)\leq C \|w_t^{(k)}- w_t^{(0)}\|^2\|J_{s,t}^{(k)}\xi\|^2
      +C\|J_{s,t}^{(k)}\xi-J_{s,t}^{(0)}\xi\|^2.
    \end{eqnarray*}
  Combining the above estimates of $I_2(t),I_3(t),I_4(t)$ with (\ref{E-4}),
  for some $C\geq 1,$ one gets
  \begin{eqnarray*}
    && \dif   \Big(\|J_{s,t}^{(k)}\xi-J_{s,t}^{(0)}\xi\|^2 \exp\big\{-C\int_s^t \|w_r^{(k)}\|_1^{4/3}
    \dif r -C(t-s)\big\}
    \Big)
    \\ &&\leq  C\|w_t^{(k)}-w_t^{(0)}\|^2 \big(\|J_{s,t}^{(k)}\xi\|_{1/2}^2+\|J_{s,t}^{(0)}\xi\|_{1/2}^2 \big)\dif t
    \\ && \quad +2\sum_{j\in \cZ_0} I_{5,j}(t) \exp\big\{-C\int_s^t \|w_r^{(k)}\|_1^{4/3}
    \dif r -C(t-s)\big\} \dif W_j(t),
  \end{eqnarray*}
  which   implies
  \begin{eqnarray*}
    && \|J_{s,t}^{(k)}\xi-J_{s,t}^{(0)}\xi\|^2 \exp\big\{-C\int_s^t \|w_r^{(k)}\|_1^{4/3}
    \dif r -C(t-s)\big\}
    \\ &&\leq  C\int_s^t \|w_r^{(k)}-w_r^{(0)}\|^2 \big(\|J_{s,r}^{(k)}\xi\|_{1/2}^2+\|J_{s,r}^{(0)}\xi\|_{1/2}^2 \big)\dif r
    \\ && \quad +2\sum_{j\in \cZ_0}  \int_s^t I_{5,j}(r) \exp\big\{-C\int_s^r  \|w_u^{(k)}\|_1^{4/3}
    \dif u -C(r-s)\big\} \dif W_j(r).
  \end{eqnarray*}
 Taking expectation in the above,
  by  Lemma  \ref{E-5} and  H\"older's inequality,   we obtain
    \begin{align}
    \nonumber & \mE \Big[\|J_{s,t}^{(k)}\xi-J_{s,t}^{(0)}\xi\|^2\exp\{-C\int_s^t \|w_r^{(k)}\|_1^{4/3}
    \dif r -C(t-s) \}\Big]
    \\   \nonumber & \leq
   C \mE
    \int_s^t\|w_r^{(k)}-w_r^{(0)}\|^2 \big(\|J_{s,r}^{(k)}\xi\|_{1/2}^2+\|J_{s,r}^{(0)}\xi\|_{1/2}^2 \big) \dif r.
    \\   \nonumber &\leq   C   \mE \Big[ \sup_{r\in [s,t]}\|w_r^{(k)}-w_r^{(0)}\|^{1/4}
    \sup_{r\in [s,t]} \big( \|J_{s,r}^{(k)}\xi\|+\|J_{s,r}^{(0)}\xi\|\big)
    \\   \nonumber & \quad\quad \quad\quad  \times
    \int_s^t \|w_r^{(k)}-w_r^{(0)}\|^{7/4} \big(\|J_{s,r}^{(k)}\xi\|_{1}+\|J_{s,r}^{(0)}\xi\|_{1}  \big) \dif r\Big]
    \\   \nonumber &\leq   C   \big( \mE \sup_{r\in [s,t]}\|w_r^{(k)}-w_r^{(0)}\|\big)^{1/4}
    \cdot\Big( \mE \sup_{r\in [s,t]} \big( \|J_{s,r}^{(k)}\xi\|+\|J_{s,r}^{(0)}\xi\|\big)^{12 }\Big)^{1/12}
    \\   \nonumber & \quad\times \Big(\mE \int_s^t \|w_r^{(k)}-w_r^{(0)}\|^{\frac{7}{4}\cdot \frac{3}{2}} \big(\|J_{s,r}^{(k)}\xi\|_{1}^{3/2}+\|J_{s,r}^{(0)}\xi\|_{1}^{3/2}  \big) \dif r   \Big)^{2/3}
    \\   \nonumber &\leq   C   \big( \mE \sup_{r\in [s,t]}\|w_r^{(k)}-w_r^{(0)}\|\big)^{1/4}
    \cdot\Big( \mE \sup_{r\in [s,t]} \big( \|J_{s,r}^{(k)}\xi\|+\|J_{s,r}^{(0)}\xi\|\big)^{12 }\Big)^{1/12}
    \\  \label{p1112-5}  & \quad\times \Big(\mE \Big[ \sup_{r\in [s,t]}
    (\|w_r^{(k)}\|^{\frac{21}{8}}+\|w_r^{(0)}\|^{\frac{21}{8}})
    \int_s^t \big(\|J_{s,r}^{(k)}\xi\|_{1}^{3/2}+\|J_{s,r}^{(0)}\xi\|_{1}^{3/2}  \big) \dif r \Big]  \Big)^{2/3}.
    \end{align}
    Using   Lemmas \ref{L:5.1}, \ref{2-3} and  H\"older's inequality,
    with regard to the last line of (\ref{p1112-5}),
   it holds that
    \begin{eqnarray*}
    &&  \mE \Big[ \sup_{r\in [s,t]}
    (\|w_r^{(k)}\|^{\frac{21}{8}}+\|w_r^{(0)}\|^{\frac{21}{8}})
    \int_s^t \big(\|J_{s,r}^{(k)}\xi\|_{1}^{3/2}+\|J_{s,r}^{(0)}\xi\|_{1}^{3/2}  \big) \dif r \Big]
    \\ &&  \leq C \Big[\mE  \sup_{r\in [s,t]}
    (\|w_r^{(k)}\|^{\frac{21}{2}}+\|w_r^{(0)}\|^{\frac{21}{2}})\Big]^{1/4}
   \Big[ \mE \int_s^t (\|J_{s,r}^{(k)}\xi\|_{1}^{2}+\|J_{s,r}^{(0)}\xi\|_{1}^{2} ) \dif r  \Big]^{3/4}
   \\ &&\leq C_{\mathfrak{R}}.
    \end{eqnarray*}
    Therefore,  the  inequality (\ref{p1112-5}) implies
    \begin{eqnarray*}
    && \mE \Big[\|J_{s,t}^{(k)}\xi-J_{s,t}^{(0)}\xi\|^2\exp\{-C\int_s^t \|w_r^{(k)}\|_1^{4/3}
    \dif r -C(t-s) \}\Big]
    \\ &&\leq C_{\mathfrak{R}}  \|P_M w_0^{(k)}-P_M w_0^{(0)}\|^{1/8}\|\xi\|+\frac{C_{\mathfrak{R}} \|\xi\| }{M^{1/32} }.
    \end{eqnarray*}
In the last step, we have used (\ref{E-1})(\ref{E-2}) and Lemma  \ref{2-3}.
The result of this lemma can be   directly obtained by the above inequality and Cauchy-Schwartz inequality. The proof is complete.

\end{proof}

\textbf{Now we are in a position to complete the   proof of   Proposition  \ref{1-2}.}

For any $j\in \cZ_0$ and $\xi\in H$ with $\|\xi \|=1,$ observe  that
\begin{eqnarray*}
  &&
 \Big|  \int_{1/2}^1 q_j^2(w_r^{(k)})\langle \xi ,
 J_{r,1}^{(k)} e_j\rangle^2 \dif r- \int_{1/2}^1  q_j^2(w_r^{(0)})\langle \xi,
 J_{r,1}^{(0)} e_j\rangle^2 \dif r \Big|
 \\ && \leq   \int_{1/2}^1 \big|  q_j^2(w_r^{(k)})\langle \xi ,
 J_{r,1}^{(k)} e_j\rangle^2 -  q_j^2(w_r^{(0)})\langle \xi,
 J_{r,1}^{(k)} e_j\rangle^2\big|\dif r
 \\ && \quad\quad     +   \int_{1/2}^1 \big|  q_j^2(w_r^{(0)})\langle \xi,
 J_{r,1}^{(k)} e_j\rangle^2 -  q_j^2(w_r^{(0)})\langle \xi,
 J_{r,1}^{(0)} e_j\rangle^2\big|\dif r
 \\ &&\leq 2\aleph   \int_{1/2}^1 |q_j(w_r^{(k)})-q_j(w_r^{(0)}) |
 \cdot \|J_{r,1}^{(k)}e_j\|^2 \dif r
 \\ &&  \quad\quad     + \aleph^2  \int_{1/2}^1 \|J_{r,1}^{(k)} e_j-J_{r,1}^{(0)} e_j \|\cdot (\|J_{r,1}^{(k)} e_j\|+\|J_{r,1}^{(0)} e_j \|)\dif r
 \\ &&\leq  2\aleph^2   \int_{1/2}^1  \|w_r^{(k)}-w_r^{(0)} \|
 \|J_{r,1}^{(k)}e_j\|^2 \dif r
  \\ &&  \quad\quad     +  \aleph^2  \int_{1/2}^1 \|J_{r,1}^{(k)} e_j-J_{r,1}^{(0)} e_j \|^{1/2}\cdot (\|J_{r,1}^{(k)} e_j\|+\|J_{r,1}^{(0)} e_j \|)^{3/2}\dif r
  \\ &&:=X_{k,j}.
\end{eqnarray*}
Thus, for any $\eps>0,$ it holds  that
\begin{eqnarray}
 \nonumber && \mP\Big(\inf_{\xi\in \cS_{\alpha,N}}\sum_{j\in \cZ_0}\int_{0}^{1} q_j^2(w_r^{(k)}) \langle \xi,
 J_{r,1}^{(k)} e_j\rangle^2\dif r<\eps_k\Big)
 \\  \nonumber &&\leq \mP\Big(\inf_{\xi\in \cS_{\alpha,N} }\sum_{j\in \cZ_0}\int_{1/2}^1 q_j^2(w_r^{(k)})  \langle \xi,
 J_{r,1}^{(k)} e_j\rangle^2\dif  r<\eps_k
 \Big)
 \\  \nonumber &&\leq \mP\Big(\inf_{\xi\in \cS_{\alpha,N} }\sum_{j\in \cZ_0}\int_{1/2}^1 q_j^2(w_r^{(0)})  \langle \xi,
 J_{r,1}^{(0)} e_j\rangle^2\dif  r<\eps_k+\sum_{j\in \cZ_0} X_{k,j}
 \Big)
 \\ \nonumber   && \leq \mP\Big(\inf_{\xi\in \cS_{\alpha,N} }\sum_{j\in \cZ_0}\int_{1/2}^1 q_j^2(w_r^{(0)})  \langle \xi,
 J_{r,1}^{(0)} e_j\rangle^2\dif  r<\eps_k+\eps
 \Big)
 \\ \label{p1111-1}  && \quad\quad\quad\quad \quad\quad \quad\quad \quad\quad +\mP\Big(\sum_{j\in \cZ_0} X_{k,j} > \eps \Big).
 \end{eqnarray}
  In view of  Lemmas  \ref{E-5}--\ref{9-3} and Lemmas  \ref{L:5.1},  \ref{2-3},
  for any $\eps>0,k\in \mN$  and $M>\max\{|j|: j\in \cZ_0\},$
  we obtain
  \begin{eqnarray*}
 && \mP\Big(\sum_{j\in \cZ_0} X_{k,j} > \eps \Big)
 \leq \sum_{j\in \cZ_0} \frac{\mE X_{k,j}}{\eps}
 \\ && \leq  \frac{2\aleph^2 }{\eps} \sum_{j\in \cZ_0}
  \int_{1/2}^1 \mE  \Big( \|w_r^{(k)}-w_r^{(0)} \|^{1/2}
   (\|w_r^{(k)}\|+\|w_r^{(0)} \|)^{1/2}
 \|J_{r,1}^{(k)}e_j\|^2\Big)  \dif r
 \\ && \quad + \frac{2\aleph^2 }{\eps}
 \int_{1/2}^1 \mE \Big(\|J_{r,1}^{(k)} e_j-J_{r,1}^{(0)} e_j \|^{1/2}\cdot (\|J_{r,1}^{(k)} e_j\|+\|J_{r,1}^{(0)} e_j \|)^{3/2}\Big) \dif r
  \\ && \leq  \frac{C\aleph^2 }{\eps} \sum_{j\in \cZ_0}
  \int_{1/2}^1 \big(\mE  \|w_r^{(k)}-w_r^{(0)} \|\big)^{1/2}
   (\mE \|w_r^{(k)}\|^2+\mE \|w_r^{(0)} \|^2)^{1/4}
 (\mE \|J_{r,1}^{(k)}e_j\|^8)^{1/4}  \dif r
 \\ && \quad + \frac{C \aleph^2 }{\eps}
  \sum_{j\in \cZ_0} \int_{1/2}^1 (\mE\|J_{r,1}^{(k)} e_j-J_{r,1}^{(0)} e_j \|)^{1/2}\cdot (\mE \|J_{r,1}^{(k)} e_j\|^3+\mE \|J_{r,1}^{(0)} e_j \|^3)^{1/2}\Big) \dif r
 \\ &&\leq   \frac{C_\mathfrak{R}  \aleph^2 }{\eps}\big(\|P_M w_0^{(k)}-P_M w_0^{(0)}\|^{1/16}+\frac{C_{\mathfrak{R}} }{M^{1/64} }\big)^{1/2}.
  \end{eqnarray*}
 Combining the above inequality with (\ref{p1111-1}),
 the following
 \begin{eqnarray*}
   && \mP\Big(\inf_{\xi\in \cS_{\alpha,N}}\sum_{j\in \cZ_0}\int_{0}^{1} q_j^2(w_r^{(k)}) \langle \xi,
 J_{r,1}^{(k)} e_j\rangle^2\dif r<\eps_k\Big)
 \\ &&\leq \mP\Big(\inf_{\xi\in \cS_{\alpha,N} }\sum_{j\in \cZ_0}\int_{1/2}^1 q_j^2(w_r^{(0)})  \langle \xi,
 J_{r,1}^{(0)} e_j\rangle^2\dif  r<\eps_k+\eps
 \Big)
 \\ && \quad\quad+   \frac{C_\mathfrak{R}  \aleph^2 }{\eps}\big(\|P_M w_0^{(k)}-P_M w_0^{(0)}\|^{1/16}+\frac{C_{\mathfrak{R}} }{M^{1/64} }\big)^{1/2}
 \end{eqnarray*}
 holds for  any $\eps>0,k\in \mN$  and $M>\max\{|j|: j\in \cZ_0\}.$
 Letting $k\rightarrow \infty$ in the above,  with the help of  our assumption (\ref{D-1}),  for any  $\eps>0$ and  $M>\max\{|j|: j\in \cZ_0\},$  one arrives at
 \begin{eqnarray}
\nonumber
   \delta_0 & \leq  &  \mP\Big(\inf_{\xi\in \cS_{\alpha,N} }\sum_{j\in \cZ_0}\int_{1/2}^1 q_j^2(w_r^{(0)})  \langle \xi,
 J_{r,1}^{(0)} e_j\rangle^2\dif  r\leq   \eps  \Big)+
  \frac{C_\mathfrak{R}  \aleph^2 }{\eps M^{1/128}}.
 \end{eqnarray}
 Letting $M\rightarrow \infty$ in the above,  it holds that
  \begin{eqnarray}
\label{B-3}
   && \delta_0 \leq
   \mP\Big(\inf_{\xi\in \cS_{\alpha,N} }\sum_{j\in \cZ_0}\int_{1/2}^1 q_j^2(w_r^{(0)})  \langle \xi,
 J_{r,1}^{(0)} e_j\rangle^2\dif  r\leq   \eps  \Big),~ \forall \eps>0.
 \end{eqnarray}
 On the other way,  (\ref{28-1})  implies that
\begin{eqnarray*}
  \lim_{\eps \rightarrow 0}\mP\Big(\inf_{\xi\in \cS_{\alpha,N} }\sum_{j\in \cZ_0}\int_{1/2}^1q_j^2(w_r^{(0)})  \langle \xi,
 J_{r,T}^{(0)} e_j\rangle^2\dif  r\leq    \eps  \Big)=0.
\end{eqnarray*}
This conflicts  with (\ref{B-3}).
The proof is  complete.
\end{appendix}

\begin{acks}[Acknowledgments]
We  would like to thank Yong Liu, Ziyu Liu and  Tianyi Pan    for  useful  discussions and suggestions.
\end{acks}

\begin{funding}
This work was  partially supported by  National Key R and  D Program of China (No. 2020YFA0712700),
NSFC (Grant Nos.
12090010, 12090014, 12471128),
the science and technology innovation Program of Hunan Province (No. 2022RC1189)  and
Key Laboratory of Random Complex Structures and Data Science, Academy
of Mathematics and Systems Science, Chinese Academy of Sciences (Grant No. 2008DP173182).
\end{funding}

\end{document}